\documentclass[11pt]{article}
\usepackage{amsmath,amsthm,amssymb,amscd,fancybox,ifthen}
\usepackage[all]{xy}

\newcommand\dbarb{\bar{\partial}_b}

\DeclareMathOperator\eo{eo}
\DeclareMathOperator\ooee{oe}

\newcommand\Dom{\operatorname{Dom}}

\newcommand\cdeg{\operatorname{c-deg}}

\newcommand\sign{\operatorname{sig}}
\newcommand\SW{\operatorname{SW}}

\newcommand\Span{\operatorname{span}}
\newcommand\spnc{\operatorname{Spin}_{\bbC}}

\newcommand\Tr{\operatorname{tr}}

\newcommand\hn{(1,-\frac 12)}
\newcommand\cT{\mathcal T}

\newcommand\cH{\mathcal H}
\newcommand\cI{\mathcal I}
\newcommand\cJ{\mathcal J}
\newcommand\cO{\mathcal O}
\newcommand\cU{\mathcal U}
\newcommand\cV{\mathcal V}
\newcommand\sO{\mathfrak O}
\newcommand\fT{\mathfrak T}

\newcommand\cK{\mathcal K}

\newcommand\cR{\mathcal R}
\newcommand\cP{\mathcal P}
\newcommand{\ccD}{\mathfrak D}
\newcommand\ho{\mathfrak H}
\newcommand\px{\tilde x}
\newcommand\pz{\tilde z}
\newcommand\bc{\boldsymbol c}
\newcommand\bzero{\boldsymbol 0}

\newcommand\tcP{\widetilde{\cP}}

\newcommand\cW{\mathcal W}

\newcommand\dbar{\bar{\pa}}
\newcommand\dbnc{\dbar\rho\rfloor}

\newcommand\bz{\bar{z}}
\newcommand\bj{\bar{j}}
\newcommand\sypr{\pi}

\newcommand\scp{\mathfrak p}
\newcommand\sd{\mathfrak d}

\newcommand\tQ{\widetilde{Q}}

\newcommand\tX{\widetilde{X}}
\newcommand\pX{X'}
\newcommand\thX{\widetilde{\pX}}
\newcommand\hX{\widehat{X}}
\newcommand\tK{\widetilde{K}}
\newcommand\sh[1]{\overset{\sqcap}{#1}}

\newcommand\esym[2]{{}^{#1}\sigma^{#2}}

\newcommand{\Spn}{S\mspace{-10mu}/ }

\newcommand\Ker{\operatorname{ker}}
\newcommand\coker{\operatorname{coker}}

\newcommand\cL{\mathcal{L}}
\newcommand\cS{\mathcal{S}}

\newcommand\cD{\mathcal{D}}

\newcommand\bomega{\bar{\omega}}

\newcommand\ha{\frac12}
\newcommand\oqu{\frac14}
\renewcommand\Re{\operatorname{Re}}

\newcommand\bbC{\mathbb C}

\newcommand\bbQ{\mathbb Q}
\newcommand\bbN{\mathbb N}

\newcommand\bbR{\mathbb R}
\newcommand\bbZ{\mathbb Z}

\newcommand\bbT{\mathbb T}
\newcommand\pa{\partial}

\newcommand\restrictedto{\upharpoonright}

\newcommand\subsubset{\subset\!\subset}

\newcommand\CI{{\mathcal C}^{\infty}}

\newcommand\CmI{{\mathcal C}^{-\infty}}

\newcommand\Id{\operatorname{Id}}

\DeclareMathOperator{\range}{range}
\DeclareMathOperator{\odd}{o}
\DeclareMathOperator{\even}{e}
\DeclareMathOperator{\Ind}{Ind}

\DeclareMathOperator{\Rind}{R-Ind}
\DeclareMathOperator{\tind}{t-Ind}

\DeclareMathOperator{\Hom}{Hom}

\newtheorem{theorem}{Theorem}
\newtheorem{proposition}{Proposition}
\newtheorem{corollary}{Corollary}
\newtheorem{lemma}{Lemma}

\theoremstyle{definition}
\newtheorem{definition}{Definition}

\theoremstyle{remark}
\newtheorem{remark}{Remark}

\begin{document}

\title{Subelliptic $\spnc$ Dirac operators, III\\
The Atiyah-Weinstein conjecture} 

\author{Charles L. Epstein\footnote{ Research partially supported by NSF grant DMS02-03795 and the Francis
J. Carey term chair.  E-mail: cle@math.upenn.edu} \\ Department of
Mathematics\\ University of Pennsylvania}

\date{Revised submitted version: January 11, 2006}

\maketitle

\centerline{\Large{\it This paper is dedicated to my wife Jane}}
 \centerline{\Large{\it for her enduring love and support.}}
\medskip
\begin{abstract}  In this paper we extend the results obtained
  in~\cite{Epstein4,Epstein3} to manifolds with $\spnc$-structures defined,
  near the boundary, by an almost complex structure. We show that on such a
  manifold with a strictly pseudoconvex boundary, there are modified
  $\dbar$-Neumann boundary conditions defined by projection operators,
  $\cR^{\eo}_+,$ which give subelliptic Fredholm problems for the $\spnc$-Dirac
  operator, $\eth^{\eo}_{+}$. We introduce a generalization of Fredholm pairs
  to the ``tame'' category. In this context, we show that the index of the
  graph closure of $(\eth^{\eo}_+,\cR^{\eo}_+)$ equals the relative index, on
  the boundary, between $\cR^{\eo}_+$ and the Calderon projector,
  $\cP^{\eo}_+.$ Using the relative index formalism, and in particular, the
  comparison operator, $\cT^{\eo}_{+},$ introduced in~\cite{Epstein4,
  Epstein3}, we prove a trace formula for the relative index that generalizes
  the classical formula for the index of an elliptic operator. Let $(X_0,J_0)$
  and $(X_1,J_1)$ be strictly pseudoconvex, almost complex manifolds, with
  $\phi:bX_1\to bX_0,$ a contact diffeomorphism. Let $\cS_0, \cS_1$ denote
  generalized Szeg\H o projectors on $bX_0, bX_1,$ respectively, and
  $\cR_0^{\eo},$ $\cR_1^{\eo},$ the subelliptic boundary conditions they
  define. If $\overline{X_1}$ is the manifold $X_1$ with its orientation
  reversed, then the glued manifold $X=X_0\amalg_{\phi}\overline{X_1}$ has a
  canonical $\spnc$-structure and Dirac operator, $\eth_X^{\eo}.$ Applying
  these results and those of our previous papers we obtain a formula for the
  relative index, $\Rind(\cS_0,\phi^*\cS_1),$
\begin{equation}
\Rind(\cS_0,\phi^*\cS_1)=\Ind(\eth^{\even}_X)-
\Ind(\eth^{\even}_{X_0},\cR^{\even}_{0})+\Ind(\eth^{\even}_{X_1},\cR^{\even}_{1}).
\label{frm1}
\end{equation}
For the special case that $X_0$ and $X_1$ are strictly pseudoconvex complex
manifolds and $\cS_0$ and $\cS_1$ are the classical Szeg\H o projectors defined
by the complex structures this formula implies that
\begin{equation}
\Rind(\cS_0,\phi^*\cS_1)=\Ind(\eth^{\even}_X)-\chi'_{\cO}(X_0)+\chi'_{\cO}(X_1),
\label{7.15.1}
\end{equation}
which is essentially the formula conjectured by Atiyah and Weinstein,
see~\cite{Weinstein}. We show that, for the case of embeddable CR-structures on
a compact, contact 3-manifold, this formula specializes to show that the
boundedness conjecture for relative indices from~\cite{Epstein} reduces to a
conjecture of Stipsicz concerning the Euler numbers and signatures of Stein
surfaces with a given contact boundary, see~\cite{stipsicz}.
\end{abstract}

\section*{Introduction}
Let $X$ be an even dimensional manifold with a Spin${}_{\bbC}$-structure,
 see~\cite{LawsonMichelsohn}.  A compatible choice of metric, $g,$ and
 connection $\nabla^{\Spn},$ define a Spin${}_{\bbC}$-Dirac operator, $\eth$
 which acts on sections of the bundle of complex spinors, $\Spn.$ This bundle
 splits as a direct sum $\Spn=\Spn^{\even}\oplus\Spn^{\odd}.$ If $X$ has a
 boundary, then the kernels and cokernels of $\eth^{\eo}$ are generally
 infinite dimensional. To obtain a Fredholm operator we need to impose boundary
 conditions. In this instance, there are no local boundary conditions for
 $\eth^{\eo}$ that define elliptic problems.  In our earlier
 papers,~\cite{Epstein4,Epstein3}, we analyzed \emph{subelliptic} boundary
 conditions for $\eth^{\eo}$ obtained by modifying the classical
 $\dbar$-Neumann and dual $\dbar$-Neumann conditions for $X,$ under the
 assumption that the $\spnc$-structure near to the boundary of $X$ is that
 defined by an integrable almost complex structure, with the boundary of
 $X$ either strictly pseudoconvex or pseudoconcave. The boundary conditions
 considered in our previous papers have natural generalizations to almost
 complex manifolds with strictly pseudoconvex or pseudoconcave boundary.

A notable feature of our analysis is that, properly understood, we show that
the natural generality for Kohn's classic analysis of the $\dbar$-Neumann
problem is that of an almost complex manifold with a strictly pseudoconvex
contact boundary. Indeed it is quite clear that analogous results hold true for
almost complex manifolds with contact boundary satisfying the obvious
generalizations of the conditions $Z(q),$ for a $q$ between $0$ and $n,$
see~\cite{FollandKohn1}. The principal difference between the integrable and
non-integrable cases is that in the latter case one must consider all form
degrees at once because, in general, $\eth^2$ does not preserve form degree.

Before going into the details of the geometric setup we briefly describe the
philosophy behind our analysis. There are three principles:
\begin{enumerate}
\item On an almost complex manifold the $\spnc$-Dirac operator, $\eth,$ is the proper
  replacement for $\dbar+\dbar^*.$
\item Indices can be computed using trace formul{\ae}.
\item The index of a boundary value problem should be expressed as a relative
  index between projectors on the boundary.
\end{enumerate}
The first item is a well known principle that I learned from
reading~\cite{duistermaat}. Technically, the main point here is that $\eth^2$
differs from a metric Laplacian by an operator of order zero. 
As to the second item, this is a basic principle in the analysis of elliptic
operators as well. It allows one to take advantage of the remarkable invariance
properties of the trace. The last item is not entirely new, but our
applications require a substantial generalization of the notion of Fredholm
pairs.  In an appendix we define \emph{tame Fredholm pairs} and prove
generalizations of many standard results. Using this approach we reduce the
Atiyah-Weinstein conjecture to Bojarski's formula, which expresses the index of
a Dirac operator on a compact manifold as a relative index of a pair of
Calderon projectors defined on a separating hypersurface. That Bojarski's
formula would be central to the proof of formula~\eqref{frm1} was suggested by
Weinstein in~\cite{Weinstein}.

The Atiyah-Weinstein conjecture, first enunciated in the 1970s, was a
conjectured formula for the index of a class of elliptic Fourier integral
operators defined by contact transformations between co-sphere bundles of
compact manifolds. We close this introduction with a short summary of the
evolution of this conjecture and the prior results.  In the original conjecture
one began with a contact diffeomorphism between co-sphere bundles:
$\phi:S^*M_1\to S^*M_0.$ This contact transformation defines a class of
elliptic Fourier integral operators. There are a variety of ways to describe an
operator from this class; we use an approach that makes the closest contact
with the analysis in this paper.

Let $(M,g)$ be a smooth Riemannian manifold; it is possible to define complex
structures on a neighborhood of the zero section in $T^*M$ so that the zero
section and fibers of $\pi:T^*M\to M$ are totally real, see~\cite{LS},
\cite{GuilleminStenzel1,GuilleminStenzel2}. For each $\epsilon>0,$ let
$B^*_\epsilon M$ denote the co-ball bundle of radius $\epsilon,$ and let
$\Omega^{n,0} B^*_{\epsilon}M$ denote the space of holomorphic $(n,0)$-forms on
$B^*_\epsilon M$ with tempered growth at the boundary. For small enough
$\epsilon>0,$ the push-forward defines maps
\begin{equation}
G_\epsilon:\Omega^{n,0} B_\epsilon^*M\longrightarrow \CmI(M),
\end{equation}
such that forms smooth up to the boundary map to $\CI(M).$ Boutet de Monvel and
Guillemin conjectured, and Epstein and Melrose proved that there is an
$\epsilon_0>0$ so that, if $\epsilon<\epsilon_0,$ then $G_\epsilon$ is an
isomorphism, see~\cite{EpsteinMelrose2}. With $S^*_\epsilon M=b B^*_\epsilon
M,$ we let $\Omega_b^{n,0} S^*_\epsilon M$ denote the distributional boundary
values of elements of $\Omega^{n,0} B_\epsilon^*M.$ One can again define a
push-forward map
\begin{equation}
G_{b\epsilon}:\Omega_b^{n,0} S_\epsilon^*M\longrightarrow \CmI(M).
\end{equation}
In his thesis, Raul Tataru showed that, for small enough $\epsilon,$ this map
is also an isomorphism, see~\cite{tataru}.  As the canonical bundle is
holomorphically trivial for $\epsilon$ sufficiently small, it suffices to work
with holomorphic functions (instead of $(n,0)$-forms).

Let $M_0$ and $M_1$ be compact manifolds and $\phi: S^* M_1\to S^* M_0$ a
contact diffeomorphism. Such a transformation canonically defines a contact
diffeomorphism $\phi_\epsilon: S_\epsilon^* M_1\to S^*_\epsilon M_0$ for all
$\epsilon>0.$ For sufficiently small positive $\epsilon,$ we define the
operator:
\begin{equation}
F^{\phi}_\epsilon f= G_{b\epsilon}^{1}\phi_\epsilon^*[ G_{b\epsilon}^{0}]^{-1}f.
\end{equation}
This is an elliptic Fourier integral operator, with canonical relation
essentially the graph of $\phi.$ The original Atiyah-Weinstein conjecture
(circa 1975) was a formula for the index of this operator as the index of the
$\spnc$-Dirac operator on the compact $\spnc$-manifold $B^*_\epsilon
M_0\amalg_{\phi} \overline{B^*_\epsilon M_1}.$ Here $\overline{X}$ denotes a
reversal of the orientation of the oriented manifold $X.$ If we let
$\cS^j_\epsilon$ denote the Szeg\H o projectors onto the boundary values of
holomorphic functions on $B^*_\epsilon M_j,\, j=0,1,$ then, using the
Epstein-Melrose-Tataru result, Zelditch observed that the index of
$F^{\phi}_\epsilon$ could be computed as the relative index between the Szeg\H
o projectors, $\cS^0_\epsilon,$ and $[\phi^{-1}]^* \cS^1_\epsilon \phi^*,$ defined
on $S^*_\epsilon M_0,$ i.e.
\begin{equation}
\Ind(F^{\phi}_{\epsilon})=\Rind(\cS^0_\epsilon, [\phi^{-1}]^* \cS^1_\epsilon \phi^*).
\end{equation}
Weinstein subsequently generalized the conjecture to allow for contact
transforms $\phi :bX_1\to bX_0,$ where $X_0, X_1$ are strictly pseudoconvex
complex manifolds with boundary, see~\cite{Weinstein}. In this paper Weinstein
suggests a variety of possible formul{\ae} depending upon whether or not the
$X_j$ are Stein manifolds.

Several earlier papers treat special cases of this conjecture (including the
original conjectured formula). In~\cite{EpsteinMelrose}, Epstein and Melrose
consider operators defined by contact transformations $\phi:Y\to Y,$ for $Y$ an
arbitrary compact, contact manifold. If $\cS$ is any generalized Szeg\H o
projector defined on $Y,$ then they show that $\Rind(\cS,[\phi^{-1}]^*
\cS\phi^*)$ depends only on the contact isotopy class of $\phi.$ In light of
its topological character, Epstein and Melrose call this relative index the
\emph{contact degree} of $\phi,$ denoted $\cdeg(\phi).$ It equals the index of
the $\spnc$-Dirac operator on the mapping torus $Z_\phi=Y\times
[0,1]/(y,0)\sim(\phi(y),1).$ Generalized Szeg\H o projectors were originally
introduced by Boutet de Monvel and Guillemin, in the context of the Hermite
calculus, see~\cite{BoutetdeMonvel-Guillemin1}. A discussion of generalized
Szeg\H o projectors and their relative indices, in the Heisenberg calculus, can
be found in~\cite{EpsteinMelrose}.

Leichtnam, Nest and Tsygan consider the case of contact transformations $\phi
:S^*M_1\to S^*M_0$ and obtain a cohomological formula for the index of
$F^{\phi}_\epsilon,$ see~\cite{LeichtNestTsy}. The approaches of these two
papers are quite different: Epstein and Melrose express the relative index as a
spectral flow, which they compute by using the extended Heisenberg calculus to
deform, through Fredholm operators, to the $\spnc$-Dirac operator on $Z_\phi.$
Leichtnam, Nest and Tsygan use the deformation theory of Lie algebroids and the
general algebraic index theorem from~\cite{NestTsygan} to obtain their formula
for the index of $F^{\phi}_\epsilon.$ In this paper we also make extensive
usage of the extended Heisenberg calculus, but the outline of the argument here
is quite different here from that in~\cite{EpsteinMelrose}.

One of our primary motivations for studying this problem was to find a formula
for the relative index between pairs of Szeg\H o projectors, $\cS_0,\cS_1,$
defined by embeddable, strictly pseudoconvex CR-structures on a compact,
3-dimensional contact manifold $(Y,H).$ In~\cite{Epstein} we conjectured that,
among small embeddable deformations, the relative index, $\Rind(\cS_0,\cS_1)$
should assume finitely many distinct values. It is shown that the relative
index conjecture implies that the set of small embeddable perturbations of an
embeddable CR-structure on $(Y,H)$ is closed in the $\CI$-topology.

Suppose that $j_0, j_1$ are embeddable CR-structures on $(Y,H),$ which bound
the strictly pseudoconvex, complex surfaces $(X_0,J_0),$ $(X_1,J_1),$
respectively. In this situation our general formula,~\eqref{7.15.1} takes a
very explicit form:
\begin{equation}
\begin{split}
\Rind(\cS_0,\cS_1)=&\dim H^{0,1}(X_0,J_0)-\dim H^{0,1}(X_1,J_1)+\\
&\frac{\sign[X_0]-\sign[X_1]+\chi[X_0]-\chi[X_1]}{4}.
\end{split}
\label{7.15.2}
\end{equation}
Here $\sign[M]$ is the signature of the oriented 4-manifold $M$ and $\chi(M)$
is its Euler characteristic. In~\cite{stipsicz}, Stipsicz conjectures that,
among Stein manifolds $(X,J)$ with $(Y,H)$ as boundary, the characteristic
numbers $\sign[X],\chi[X]$ assume only finitely many values. Whenever
Stipsicz's conjecture is true it implies a strengthened form of the relative
index conjecture: the function $\cS_1\mapsto\Rind(\cS_0,\cS_1)$ is bounded from
above throughout the entire deformation space of embeddable CR-structures on
$(Y,H).$ Many cases of Stipsicz's conjecture are proved
in~\cite{OzbagciStipsicz,stipsicz}. As a second consequence of~\eqref{7.15.2}
we show that, if $\dim M_j=2,$ then  $\Ind(F^\phi_\epsilon)=0.$

{\small \centerline{Acknowledgments} Boundary conditions similar to those
considered in this paper, as well as the idea of finding a geometric formula
for the relative index were first suggested to me by Laszlo Lempert. I would
like to thank Richard Melrose for our years of collaboration on problems in
microlocal analysis and index theory; it provided many of the tools needed to
do the current work. I would also like to thank Alan Weinstein for very useful
comments on an early version of this paper.  I am very grateful to
John Etnyre for references to the work of Ozbagci and Stipsicz and our many
discussions about contact manifolds and complex geometry, and to Julius Shaneson
for providing the proof of Lemma~\ref{lem10}. I would like to thank the referee
for many suggestions that improved the exposition and for simplifying the
proof of Proposition~\ref{prp10}.}

\section{Outline of Results}
Let $X$ be an even dimensional manifold with a $\spnc$-structure and let
$\Spn\to X$ denote the bundle of complex spinors. A choice of metric on $X$ and
compatible connection, $\nabla^{\Spn},$ on the bundle $\Spn$ define the
$\spnc$-Dirac operator, $\eth:$
\begin{equation}
\eth \sigma=\sum_{j=0}^{\dim X} \bc(\omega_j)\cdot\nabla^{\Spn}_{V_j} \sigma,
\end{equation}
with $\{V_j\}$ a local framing for the tangent bundle and $\{\omega_j\}$ the
dual coframe. Here $\bc(\omega)\cdot$ denotes the Clifford action of $T^*X$ on
$\Spn.$  It is customary to
split $\eth$ into its chiral parts: $\eth=\eth^{\even}+\eth^{\odd},$ where
$$\eth^{\eo}:\CI(X;\Spn^{\eo})\longrightarrow\CI(X;\Spn^{\ooee}).$$ The
operators $\eth^{\odd}$ and $\eth^{\even}$ are formal adjoints.

An almost complex structure on $X$ defines a Spin${}_{\bbC}$-structure, and
bundle of complex spinors $\Spn,$ see~\cite{duistermaat}.  The bundle of
complex spinors is canonically identified with $\oplus_{q\geq 0}\Lambda^{0,q}.$
We use the notation
\begin{equation}
\Lambda^{\even}=\bigoplus\limits_{q=0}^{\lfloor\frac{n}{2}\rfloor}
\Lambda^{0,2q}\quad
\Lambda^{\odd}=\bigoplus\limits_{q=0}^{\lfloor\frac{n-1}{2}\rfloor}\Lambda^{0,2q+1}.
\end{equation} 
These bundles are in turn canonically identified with the bundles of even and
odd spinors, $\Spn^{\eo},$ which are defined as the $\pm 1$-eigenspaces of the
orientation class. 
A metric $g$ on $X$ is compatible with the almost complex structure, if for
every $x\in X$ and $V,W\in T_xX,$ we have:
\begin{equation}
g_x(J_xV,J_xW)=g_x(V,Y).
\label{eqn4}
\end{equation}

Let $X$ be a compact manifold with a co-oriented contact structure $H\subset
TbX,$on its boundary. Let $\theta$ denote a globally defined contact form in
the given co-orienta-tion class. An almost complex structure $J$ defined in a
neighborhood of $bX$ is compatible with the contact structure if, for every
$x\in bX,$ we have
\begin{equation}
\begin{split}
&J_xH_x\subset H_x,\text{ and for all }V,W\in H_x\\
&d\theta_x(J_xV,W)+d\theta_x(V,J_x W)=0\\
&d\theta_x(V,J_x V)>0,\text{ if }V\neq 0.
\end{split}
\end{equation}
We usually assume that $g\restrictedto_{H\times H}=d\theta(\cdot,J\cdot).$ If
the almost complex in not integrable, then $\eth^2$ does not preserve the
grading of $\Spn$ defined by the $(0,q)$-types.

As noted, the almost complex structure defines the bundles $T^{1,0}X, T^{0,1}X$
as well as the form bundles $\Lambda^{0,q}X.$ This in turn defines the
$\dbar$-operator.  The bundles $\Lambda^{0,q}$ have a splitting at the boundary
into almost complex normal and tangential parts, so that a section $s$
satisfies:
\begin{equation}
s\restrictedto_{bX}=s^t+\dbar\rho\wedge s^n,\text{ where }\dbnc s^t=\dbnc
s^n=0.
\end{equation}
Here $\rho$ is defining function for $bX.$ The $\dbar$-Neumann condition for
sections $s\in\CI(X;\Lambda^{0,q})$ is the requirement that
\begin{equation}
\dbnc [s]_{bX}=0,
\end{equation}
i.e., $s^n=0.$ As before this does not impose any requirement on forms of
degree $(0,0).$ 

The contact structure on $bX$ defines the class of generalized Szeg\H o
projectors acting on scalar functions, see~\cite{Epstein3,EpsteinMelrose} for the
definition. Using the identifications of $\Spn^{\eo}$ with $\Lambda^{0,\eo},$ a
generalized Szeg\H o projector, $\cS,$ defines a modified (strictly
pseudoconvex) $\dbar$-Neumann condition as follows:
\begin{equation}
\begin{split}
&\cR \sigma^{00}\overset{d}{=}\cS[\sigma^{00}]_{bX}=0\\
&\cR \sigma^{01}\overset{d}{=}(\Id-\cS)[\dbnc \sigma^{01}]_{bX}=0\\
&\cR \sigma^{0q}\overset{d}{=}[\dbnc\sigma^{0q}]_{bX}=0,\text{ for }q>1.
\end{split}
\end{equation}
We choose the defining function so that $s^t$ and $\dbar\rho\wedge s^n$ are
orthogonal, hence the mapping $\sigma\to \cR\sigma$ is a self adjoint
projection operator.  Following the practice in~\cite{Epstein4,Epstein3} we use
$\cR^{\eo}$ to denote the restrictions of this projector to the subbundles of
even and odd spinors. 

We follow the conventions for the $\spnc$-structure and Dirac operator on an
almost complex manifold given in~\cite{duistermaat}. Lemma 5.5
in~\cite{duistermaat} states that the principal symbol of $\eth_X$ agrees with
that of the Dolbeault-Dirac operator $\dbar+\dbar^*,$ which implies that
$(\eth_X^{\eo},\cR^{\eo})$ are formally adjoint operators. It is a consequence
of our analysis that, as unbounded operators on $L^2,$
\begin{equation}
(\eth_X^{\eo},\cR^{\eo})^*=\overline{(\eth_X^{\ooee},\cR^{\ooee})}.
\end{equation}
The almost complex structure is only needed to define the boundary condition.
Hence we assume that $X$ is a $\spnc$-manifold, where the $\spnc$-structure is
defined by an almost complex structure $J$ defined in a neighborhood of the
boundary.

In this paper we begin by showing that the analytic results obtained in our
earlier papers remain true in the almost complex case. As noted above, this
shows that integrability is not needed for the validity of Kohn's estimates for
the $\dbar$-Neumann problem. By working with $\spnc$-structures we are able to
fashion a much more flexible framework for studying index problems than that
presented in~\cite{Epstein4,Epstein3}. As before, we compare the projector
$\cR$ defining the subelliptic boundary conditions with the Calderon projector
for $\eth,$ and show that these projectors are, in a certain sense, relatively
Fredholm. These projectors are not relatively Fredholm in the usual sense of
say Fredholm pairs in a Hilbert space, used in the study of elliptic boundary
value problems. We circumvent this problem by extending the theory of Fredholm
pairs to that of \emph{tame Fredholm pairs}. We then use our analytic results
to obtain a formula for a parametrix for these subelliptic boundary value
problems that is precise enough to prove, among other things, higher norm
estimates. The extended Heisenberg calculus introduced
in~\cite{EpsteinMelrose3} remains at the center of our work. The basics of this
calculus are outlined in~\cite{Epstein3}.

If $\cR^{\eo}$ are  projectors defining  modified $\dbar$-Neumann
conditions and $\cP^{\eo}$ are the Calderon projectors, then we show that the
comparison operators,
\begin{equation}
\cT^{\eo}=\cR^{\eo}\cP^{\eo}+(\Id-\cR^{\eo})(\Id-\cP^{\eo})
\end{equation}
are graded elliptic elements of the extended Heisenberg calculus. As such there are
parametrices $\cU^{\eo}$ that satisfy
\begin{equation}
\cT^{\eo}\cU^{\eo}=\Id-K_1^{\eo},\quad\cU^{\eo}\cT^{\eo}=\Id-K_2^{\eo},
\end{equation}
where $K_1^{\eo}, K_2^{\eo}$ are smoothing operators. We define Hilbert spaces,
$\cH_{\cU^{\eo}}$  to be the closures of $\CI(bX;\Spn^{\eo}\restrictedto_{bX})$
    with respect to the inner products
\begin{equation}
\langle\sigma,\sigma\rangle_{\cU^{\eo}}=\langle\sigma,\sigma\rangle_{L^2}+
\langle\cU^{\eo}\sigma,\cU^{\eo}\sigma\rangle_{L^2}.
\label{eqn1}
\end{equation}
The operators $\cR^{\eo}\cP^{\eo}$ are Fredholm from $\range\cP^{\eo}\cap L^2$
to $\range\cR^{\eo}\cap\cH_{\cU^{\eo}}.$ As usual, we let
$\Rind(\cP^{\eo},\cR^{\eo})$ denote the indices of these restrictions; we show
that
\begin{equation}
\Ind(\eth^{\eo},\cR^{\eo})=\Rind(\cP^{\eo},\cR^{\eo}).
\label{eqn2}
\end{equation}

Using the standard formalism for computing indices we show that
\begin{equation}
\Rind(\cP^{\eo},\cR^{\eo})=\Tr\cR^{\eo} K_1^{\eo}\cR^{\eo}-\Tr\cP^{\eo}
K_2^{\eo}\cP^{\eo}.
\label{eqn3}
\end{equation}
There is some subtlety in the interpretation of this formula in that
$\cR^{\eo} K_1^{\eo}\cR^{\eo}$ act on $\cH_{\cU^{\eo}}.$ But, as is also used
implicitly in the elliptic case, we show that the computation of the trace does
not depend on the topology of the underlying Hilbert space. Among other things,
this formula allows us to prove that the indices of the boundary problems
$(\eth^{\eo},\cR^{\eo})$ depend continuously on the data defining the
boundary condition and the $\spnc$-structure, allowing us to employ
deformation arguments.

To obtain the gluing formula we use the invertible double construction
introduced in~\cite{BBW}. Using this construction, we are able to express the
relative index between two generalized Szeg\H o projectors as the index of the
$\spnc$-Dirac operators on a compact manifold with corrections coming from
boundary value problems on the ends. Let $X_0, X_1$ be $\spnc$-manifolds with
contact boundaries.  Assume that the $\spnc$-structures are defined in
neighborhoods of the boundaries by compatible almost complex structures, such
that $bX_0$ is contact isomorphic to $bX_1,$ let $\phi: bX_1\to bX_0$ denote a
contact diffeomorphism. If $\overline{X_1}$ denotes $X_1$ with its orientation
reversed, then $\tX_{01}=X_0\amalg_\phi\overline{X_1}$ is a compact manifold with a
canonical $\spnc$-structure and Dirac operator, $\eth^{\eo}_{\tX_{01}}.$ Even if $X_0$
and $X_1$ have globally defined almost complex structures, the manifold $\tX_{01},$ in
general, does not. In case $X_0=X_1,$ as $\spnc$-manifolds, then this is the
invertible double introduced in~\cite{BBW}, where they show that $\eth_{\tX_{01}}$ is
an invertible operator.

Let $\cS_0,\cS_1$ be generalized Szeg\H o projectors on $bX_0, bX_1,$
respectively. If $\cR^{\even}_0,$ $\cR^{\even}_1$ are the subelliptic boundary
conditions they define, then the main result of this paper is the following
formula:
\begin{equation}
\Rind(\cS_0,\cS_1)=\Ind(\eth^{\even}_{\tX_{01}})-\Ind(\eth^{\even}_{X_0},\cR^{\even}_{0})+
\Ind(\eth^{\even}_{X_1},\cR^{\even}_{1}).
\label{eqn6}
\end{equation}
As detailed in the introduction, such a formula was conjectured, in a more
restricted case, by Atiyah and Weinstein, see~\cite{Weinstein}. Our approach
differs a little from that conjectured by Weinstein, in that $\tX_{01}$ is
constructed using the extended double construction rather than the
stabilization of the almost complex structure on the glued space described
in~\cite{Weinstein}. A result of Cannas da Silva implies that the stable almost
complex structure on $\tX_{01}$ defines a $\spnc$-structure, which very likely agrees
with that used here, see~\cite{GGK}. Our formula is very much in the spirit
suggested by Atiyah and Weinstein, though we have not found it necessary to
restrict to $X_0, X_1$ to be Stein manifolds (or even complex manifolds), nor
have we required the use of ``pseudoconcave caps'' in the non-Stein case. It is
quite likely that there are other formul{\ae} involving the pseudoconcave caps
and they will be considered in a subsequent publication.

In the case that $X_0$ is isotopic to $X_1$ through $\spnc$-structures
compatible with the contact structure on $Y,$ then $\tX_{01},$ with its canonical
$\spnc$-structure, is isotopic to the invertible double of $X_0\simeq X_1.$
In~\cite{BBW} it is shown that in this case, $\eth^{\eo}_{\tX_{01}}$ are invertible
operators and hence $\Ind(\eth^{\eo}_{\tX_{01}})=0.$ Thus~\eqref{eqn6} states that
\begin{equation}
\Rind(\cS_0,\cS_1)=\Ind(\eth^{\even}_{X_1},\cR^{\even}_{1})-
\Ind(\eth^{\even}_{X_0},\cR^{\even}_{0}).
\label{eqn7}
\end{equation}

If $X_0\simeq X_1$ are complex manifolds with strictly pseudoconvex boundaries,
and the complex structures are isotopic as above (through compatible almost
complex structures), and the Szeg\H o projectors are those defined by the
complex structure, then formula (77) in~\cite{Epstein4} implies that
$\Ind(\eth^{\even}_{X_j},\cR^{\even}_{j})=\chi_{\cO}'(X_j)$ and therefore:
\begin{equation}
\Rind(\cS_0,\cS_1)=\chi_{\cO}'(X_1)-\chi_{\cO}'(X_0).
\label{eqn8}
\end{equation}
When $\dim_{\bbC} X_j=2,$  this formula becomes:
\begin{equation}
\Rind(\cS_0,\cS_1)=\dim H^{0,1}(X_0)-\dim H^{0,1}(X_1),
\label{eqn9}
\end{equation}
which has applications to the relative index conjecture in~\cite{Epstein}.  In
the case that $\dim_{\bbC}X_j=1,$ a very similar formula was obtained by Segal
and Wilson, see~\cite{SegalWilson,Kang}.  A detailed analysis of the complex
2-dimensional case is given in Section~\ref{s.3d}, where we
prove~\eqref{7.15.2}. 

In Section~\ref{sec11} we show how these results can be extended to allow for
vector bundle coefficients. An interesting consequence of this analysis is a
proof, which makes no mention of K-theory, that the index of a classically
elliptic operator on a compact manifold $M$ equals that of a $\spnc$-Dirac
operator on the glued space $B^*M\amalg_{S^*M} \overline{B^*M}.$ Hence, using
relative indices and the extended Heisenberg calculus, along with Getzler's
rescaling argument we obtain an entirely analytic proof of the Atiyah-Singer
formula.

\begin{remark} In this paper we restrict our attention to the pseudoconvex
  case. There are analogous results for other cases with non-degenerate
  $d\theta(\cdot,J\cdot).$ We will return to these in a later publication. The
  subscript $+$ sometimes refers to the fact that the underlying manifold is
  pseudoconvex. Sometimes, however, we use $\pm$ to designate the two sides of
  a separating hypersurface.
\end{remark}

\section{The symbol of the Dirac Operator and its inverse}\label{ss.2}
In this section we show that, under appropriate geometric hypotheses, the
results of Sections 2--5 of~\cite{Epstein3} remain valid, with small
modifications, for the $\spnc$-Dirac operator on an almost complex manifold,
with strictly pseudoconvex boundary. As noted above the $\spnc$-structure only
need be defined by an almost complex structure near the boundary. This easily
implies that the operators $\cT^{\eo}_{+}$ are elliptic elements of the
extended Heisenberg calculus. To simplify the exposition we treat only the
pseudoconvex case. The results in the pseudoconcave case are entirely
analogous. For simplicity we also omit vector bundle coefficients. There is no
essential difference if they are included; the modifications necessary to treat
this case are outlined in Section~\ref{sec11}.

Let $X$ be a manifold with boundary, $Y.$ We suppose that $(Y,H)$ is a contact
manifold and $X$ has an almost complex structure $J,$ defined near the
boundary, compatible with the contact structure, with respect to which the
boundary is strictly pseudoconvex, see~\cite{ABKLR}. We let $g$ denote a metric
on $X$ compatible with the almost complex structure: for every $x\in X, V,W\in
T_xX$ we have
\begin{equation}
g_x(J_xV,J_xW)=g_x(V,W).
\label{eqn1.1}
\end{equation}
We suppose that $\rho$ is a defining function for the boundary of $X$ that is
negative on $X.$ Let $\dbar$ denote the (possibly non-integrable)
$\dbar$-operator defined by $J.$ We assume that $JH\subset H,$ and that the one
form,
\begin{equation}
\theta=\frac{i}{2}\dbar\rho\restrictedto_{TbX},
\end{equation}
is a contact form for $H.$ The quadratic form defined on $H\times H$ by
\begin{equation}
\cL(V,W)=d\theta(V,JW)
\end{equation}
is assumed to be positive definite. In the almost complex category this is the
statement that $bX$ is strictly pseudoconvex.

Let $T$ denote the Reeb vector field: $\theta(T)=1, i_Td\theta=0.$ For
simplicity we assume that 
\begin{equation}
g\restrictedto_{H\times H}=\cL\text{ and }g(T,V)=0,\quad\forall V\in H.
\label{eqn1.2}
\end{equation}
Note that~\eqref{eqn1.1}
and~\eqref{eqn1.2} imply that $J$ is compatible with $d\theta$ in that, for all
$V,W\in H$ we have
\begin{equation}
d\theta(JV,JW)=d\theta(V,W)\text{ and }d\theta(V,JV)>0\text{ if }V\neq 0.
\end{equation}
\begin{definition} Let $X$ be a $\spnc$-manifold with almost complex structure
  $J,$ defined near $bX.$ If the $\spnc$-structure near $bX$ is that specified
  by $J,$ then the quadruple $(X,J,g,\rho)$
  satisfying~\eqref{eqn1.1}--\eqref{eqn1.2} defines a \emph{normalized
  strictly pseudoconvex $\spnc$-manifold}.
\end{definition}

On an almost complex manifold with compatible metric there is a
$\spnc$-struc-ture so that the bundle of complex spinors $\Spn\to X$ is a
complex Clifford module. As noted above, if the $\spnc$-structure is defined by
an almost complex structure, then $\Spn\simeq\oplus\Lambda^{0,q}.$ Under this
isomorphism, the Clifford action of a real one-form $\xi$ is given by
\begin{equation}
\bc(\xi)\cdot\sigma\overset{d}{=}(\xi-iJ\xi)\wedge\sigma-\xi\rfloor\sigma.
\end{equation}
It is extended to the complexified Clifford algebra complex linearly.  We
largely follow the treatment of $\spnc$-geometry given in~\cite{duistermaat},
though with some modifications to make easier comparisons with the results of
our earlier papers.

There is a compatible connection $\nabla^{\Spn}$ on $\Spn$ and a formally
self adjoint $\spnc$-Dirac operator defined on sections of $\Spn$ by
\begin{equation}
\eth\sigma=\frac{1}{2}\sum_{j=1}^{2n}\bc(\omega_j)\cdot\nabla^{\Spn}_{V_j}\sigma,
\end{equation}
with $\{V_j\}$ a local framing for the tangent bundle and $\{\omega_j\}$ the
dual coframe.  Here we differ slightly from~\cite{duistermaat} by including the
factor $\frac 12$ in the definition of $\eth.$ This is so that, in the case
that $J$ is integrable, the leading order part of $\eth$ is $\dbar+\dbar^*$
(rather than $2(\dbar+\dbar^*)),$ which makes for a more direct comparison with
results in~\cite{Epstein4,Epstein3}.

The spinor bundle splits into even and odd components $\Spn^{\eo},$ and the
Dirac operator splits into  even and odd parts, $\eth^{\eo},$ where
\begin{equation}
\eth^{\eo}:\CI(X;\Spn^{\eo})\longrightarrow \CI(X;\Spn^{\ooee}).
\end{equation}
Note that, in each fiber, Clifford multiplication by a non-zero co-vector gives an
isomorphism $\Spn^{\eo}\leftrightarrow \Spn^{\ooee}.$

Fix a point $p$ on the boundary of $X$ and let $(x_1,\dots,x_{2n})$ denote
normal coordinates centered at $p.$ This means that
\begin{enumerate}
\item $p\leftrightarrow (0,\dots,0)$
\item The Hermitian metric tensor $g_{i\bj}$ in these coordinates satisfies
\begin{equation}
g_{i\bj}=\frac{1}{2}\delta_{i\bj}+O(|x|^2).
\label{7.2.1}
\end{equation}
\end{enumerate}
If $V\in T_pX$ is a unit vector, then $V^{0,1}=\frac 12(V+iJV),$  and
\begin{equation}
\langle V^{0,1},V^{0,1}\rangle_g=\frac 12.
\end{equation}

Without loss of generality we may also assume that the coordinates are ``almost
complex'' and adapted to the contact geometry \emph{at $p$}: that is the vectors
$\{\pa_{x_j}\}\subset T_pX$ satisfy
\begin{equation}
\begin{split}
&J_p\pa_{x_j}=\pa_{x_{j+n}}\text{ for }j=1,\dots,n\\
&\{\pa_{x_2},\dots,\pa_{x_{2n}}\}\in T_pbX\\
&\{\pa_{x_2},\dots,\pa_{x_n},\pa_{x_{n+2}},\dots,\pa_{x_{2n}}\}\in H_p.
\end{split}
\label{3.1.2}
\end{equation}
We let 
$$z_j=x_j+ix_{j+n}.$$
As $d\rho\restrictedto_{bX}=0,$ equation~\eqref{3.1.2} implies that
\begin{equation}
\rho(z)=-\frac{2}{\alpha}\Re z_1+\langle az,z\rangle+\Re (bz,z)+O(|z|^3).
\label{7.2.2}
\end{equation}
In this equation $\alpha>0,$ $a$ and $b$ are $n\times n$ complex matrices, with
$a=a^*,$ $b=b^t,$ and
\begin{equation}
\langle w,z\rangle=\sum_{j=1}^n w_j\bz_j\text{ and }(w,z)=\sum_{j=1}^n w_jz_j.
\end{equation}
With these normalizations we have the following formul{\ae} for the contact
form at $p:$
\begin{lemma} Under the assumptions above
\begin{equation}
\theta_p=-\frac{1}{2\alpha}dx_{n+1}\text{ and
}d\theta_p=\sum\limits_{j=2}^{n}dx_j\wedge dx_{j+n}.
\label{eqn3.1.3}
\end{equation}
\end{lemma}
\begin{proof} The formula for $\theta_p$ follows from~\eqref{7.2.2}.  The
  normality of the coordinates,~\eqref{eqn1.2} and~\eqref{3.1.2} imply that, for
  a one-form $\phi_p$ we have
\begin{equation}
d\theta_p=\sum\limits_{j=2}^{n}dx_j\wedge dx_{j+n}+\theta_p\wedge\phi_p.
\end{equation}
The assumption that the Reeb vector field is orthogonal to $H_p$
and~\eqref{3.1.2} imply that $\pa_{x_{n+1}}$ is a multiple of the Reeb vector
field. Hence $\phi_p=0.$
\end{proof}

For symbolic calculations the following notation proves very useful: a term
which is a symbol of order at most $k$ vanishing, at $p,$ to order $l$ is
denoted by $\sO_{k}(|x|^l).$ As we work with a variety of operator calculi, it
is sometimes necessary to be specific as to the sense in which the order should
be taken. The notation $\sO^{C}_j$ refers to terms of order at most $j$ in the
sense of the symbol class $C.$ If no symbol class is specified, then the order
is with respect to the classical, radial scaling. If no rate of vanishing is
specified, it should be understood to be $O(1).$

If $\{f_j\}$ is an orthonormal frame for $TX,$ then the Laplace
operator on the spinor bundle is defined by
\begin{equation}
\Delta=\sum_{j=1}^{2n}\nabla^{\Spn}_{f_j}\circ\nabla^{\Spn}_{f_j}-
\nabla^{\Spn}_{\nabla^g_{f_j}f_j}.
\end{equation}
Here $\nabla^g$ is the Levi-Civita connection on $TX.$ As explained
in~\cite{duistermaat}, the reason for using the $\spnc$-Dirac operator
as a replacement for $\dbar+\dbar^*$ is because of its very close connection to
the Laplace operator.
\begin{proposition}\label{prp1.1} Let $(X,g,J)$ be a Hermitian, almost complex manifold and
  $\eth$ the $\spnc$-Dirac operator defined by this data. Then
\begin{equation}
\eth^2=\frac{1}{2}\Delta+R,
\label{5.23.3}
\end{equation}
where $R:\Spn\to\Spn$ is an endomorphism.
\end{proposition}
After changing to the normalizations used here, e.g. $\langle
V^{0,1},V^{0,1}\rangle_g=\frac 12,$ this is Theorem 6.1 in~\cite{duistermaat}. Using this
result we can compute the symbols of $\eth$ and $\eth^2$ at $p.$ Recall that
\begin{equation}
\nabla^g\pa_{x_k}=O(|x|).
\end{equation}
We can choose a local orthonormal framing for $\Spn,$ $\{\sigma_J\}$
($J=(j_1,\dots,j_q)$ with $1\leq j_1<\dots<j_q\leq n$) so that
\begin{equation}
\sigma_J-d\bz^J=O(|x|)\text{ and }\nabla^{\Spn}\sigma_J=O(|x|)
\label{5.23.4}
\end{equation}
as well.

With respect to this choice of frame, the symbol of $\eth,$ in a geodesic normal
coordinate system, is
\begin{equation}
\sigma(\eth)(x,\xi)=d_1(x,\xi)+d_0(x).
\label{5.23.1}
\end{equation}
Because the connection coefficients vanish at $p$ we obtain:
\begin{equation}
d_1(x,\xi)=d_1(0,\xi)+\sO_1(|x|),\quad
d_0(z)=\sO_0(|x|).
\label{7.2.3}
\end{equation}
The linear polynomial $d_1(0,\xi)$ is the symbol of $\dbar+\dbar^*$ on $\bbC^n$
with respect to the flat metric.  This is slightly different from the K\"ahler
case where $d_1(x,\xi)-d_1(0,\xi)=\sO_{1}(|x|^2).$ First order vanishing is
sufficient for our applications, we only needed the quadratic vanishing to
obtain the formula for the symbol of $\eth^2,$ obtained here from
Proposition~\ref{prp1.1}. 

Proposition~\ref{prp1.1} implies that
\begin{equation}
\sigma(\eth^2)(x,\xi)=\sigma(\Delta-R)(x,\xi)=
\Delta_2(x,\xi)+\Delta_1(x,\xi)+\Delta_0(x),
\end{equation}
where $\Delta_j$ is a polynomial in $\xi$ of degree $j$ and
\begin{equation}
\begin{split}
\Delta_2(x,\xi)=\Delta_2&(0,\xi)+\sO_2(|x|^2)\\
\Delta_1(x,\xi)=\sO_1(|x|), &\quad \Delta_0(x,\xi)=\sO_0(1).
\end{split}
\label{7.2.4}
\end{equation}
Because we are working in geodesic normal coordinates,  the principal symbol at
$p$ is
\begin{equation}
\Delta_2(0,\xi)=\frac{1}{2}|\xi|^2\otimes\Id.
\label{5.23.2}
\end{equation}
Here $\Id$ is the identity homomorphism on the appropriate bundle.  These
formul{\ae} are justified in Section~\ref{sec11}, where we explain the
modifications needed to include vector bundle coefficients.

The manifold $X$ can be included into a larger manifold $\tX$ (the invertible
double) in such a way that its $\spnc$-structure and Dirac operator extend
smoothly to $\tX$ and such that the extended operators $\eth^{\eo}$ are
invertible. We return to this construction in Section~\ref{sec6}. Let $Q^{\eo}$
denote the inverses of $\eth^{\eo}$ extended to $\tX.$ These are classical
pseudodifferential operators of order $-1.$

We set $\tX\setminus Y=\tX_{+}\amalg \tX_{-},$ where $\tX_+=X;$ note that
$\rho<0$ on $\tX_+,$ and $\rho>0$ on $\tX_-.$ Let $r_{\pm}$ denote the
operations of restriction of a section of $\Spn^{\eo},$ defined on $\tX$ to
$\tX_{\pm},$ and $\gamma_\epsilon$ the operation of restriction of a smooth
section of $\Spn^{\eo}$ to $Y_{\epsilon}=\{\rho^{-1}(\epsilon)\}.$ 
Define the operators
\begin{equation}
\tK^{\eo}_{\pm}\overset{d}{=} r_{\pm} Q^{\eo}\gamma^*_0:\CI(Y;\Spn^{\ooee}\restrictedto_{Y})
\longrightarrow \CI(\tX_{\pm};\Spn^{\eo}).
\end{equation}
Here $\gamma^*_0$ is the formal adjoint of $\gamma_0.$ We recall that, along
$Y,$ the symbol $\sigma_1(\eth^{\eo},d\rho)$ defines an isomorphism
\begin{equation}
\sigma_1(\eth^{\eo},d\rho):\Spn^{\eo}\restrictedto_{Y}\longrightarrow
\Spn^{\ooee}\restrictedto_{Y}.
\end{equation}
Composing, we get the usual Poisson operators
\begin{equation}
\cK^{\eo}_{\pm}=\frac{\mp}{i\|d\rho\|}\tK^{\eo}_{\pm}\circ\sigma_1(\eth^{\eo},d\rho):
\CI(Y;\Spn^{\eo}\restrictedto_{Y})
\longrightarrow \CI(\tX_{\pm};\Spn^{\eo}),
\end{equation}
which map sections of $\Spn^{\eo}\restrictedto_{Y}$ into the nullspaces of
$\eth^{\eo}_{\pm}.$ The factor $\mp$ is inserted because $\rho<0$ on $X.$

The Calderon projectors are defined by
\begin{equation}
\cP^{\eo}_{\pm} s\overset{d}{=}\lim_{\mp\epsilon\to 0^+}
\gamma_{\epsilon}\cK^{\eo}_{\pm}s\text{ for }s\in\CI(Y;\Spn^{\eo}\restrictedto_{Y}).
\end{equation}
The fundamental result of Seeley is that $\cP^{\eo}_{\pm}$ are classical
pseudodifferential operators of order $0.$ The ranges of these operators are the
boundary values of elements of $\Ker\eth^{\eo}_{\pm}.$ Seeley gave a
prescription for computing the symbols of these operators using contour
integrals, which we do not repeat here, as we shall be computing these symbols in
detail in the following sections, see~\cite{Seeley0}.
\begin{remark}[{\bf Notational remark}] Unlike in~\cite{Epstein4,Epstein3}, the
  notation $\cP^{\eo}_+$ and $\cP^{\eo}_-$ refers to the Calderon projectors
  defined on the two sides of a separating hypersurface in a single manifold
  $\tX,$ with an invertible $\spnc$-Dirac operator. This is the more standard
  usage; in this case we have the identities $\cP^{\eo}_++\cP^{\eo}_-=\Id.$ In
  our earlier papers $\cP^{\eo}_+$ are the Calderon projectors on a
  pseudoconvex manifold, and $\cP^{\eo}_-,$ the Calderon projectors on a
  pseudoconcave manifold.
\end{remark}

Given the formul{\ae} above for $\sigma(\eth)$ and $\sigma(\eth^2)$ the
computation of the symbol of $Q^{\eo}$ proceeds exactly as in the K\"ahler
case. As we only need the principal symbol, it suffices to do the computations
in the  fiber over a fixed point $p\in bX.$ Set
\begin{equation}
\sigma(Q^{\eo})=q=q_{-1}+q_{-2}+\dots
\end{equation}
We summarize the results of these calculations in the following proposition:
\begin{proposition} Let $(X,J,g,\rho)$ define a normalized strictly
  pseudoconvex $\spnc$-manifold. For $p\in bX,$ let $(x_1,\dots,x_{2n})$
  denote boundary adapted, geodesic normal coordinates centered at $p.$ The symbols of $Q^{\eo}$
  at $p$ are given by
\begin{equation}
\begin{split}
&q_{-1}(\xi)=\frac{2d_1(\xi)}{|\xi|^2}\\
&q_{-2}=\sO_{-2}(|z|)
\end{split}
\label{5.13.1}
\end{equation}
Here $\xi$ are the coordinates on $T^*_pX$ defined by $\{dx_j\},$ $|\xi|$ is
the standard Euclidean norm, and $d_1(\xi)$ is the symbol of $\dbar+\dbar^*$ on
$\bbC^n$ with respect to the flat metric. For $k\geq 2$ we have:
\begin{equation}
q_{-2k}=\sum_{j=0}^{l_k}\frac{\sO_{2j}(1)}{|\xi|^{2(k+j)}},\quad
q_{-(2k-1)}=\sum_{j=0}^{l_k'}\frac{\sO_{2j+1}(1)}{|\xi|^{2(k+j)}}.
\label{7.6.6}
\end{equation}
The terms in the numerators of~\eqref{7.6.6} are polynomials in $\xi$ of the
indicated degrees.
\end{proposition}

In order to compute the symbol of the Calderon projector, we  introduce
boundary adapted coordinates, $(t,x_2,\dots, x_{2n})$ where
\begin{equation}
t=-\frac{\alpha}{2}\rho(z)=x_1+O(|x|^2).
\end{equation}
Note that $t$ is positive on a pseudoconvex manifold and $dt$ is an inward
pointing, unit co-vector.

We need to use the change of coordinates formula to
express the symbol in the new variables. From~\cite{Hormander3} we obtain the
following prescription: Let $w=\phi(x)$ be a diffeomorphism and $c(x,\xi)$ the
symbol of a classical pseudodifferential operator $C$. Let $(w,\eta)$ be linear
coordinates in the cotangent space, then $c_\phi(w,\eta),$ the symbol of $C$
in the new coordinates, is given by
\begin{equation}
c_{\phi}(\phi(x),\eta)\sim
\sum_{k=0}^{\infty}\sum_{\theta\in\cI_k}
\frac{(-i)^k\pa_\xi^\theta c(x,d\phi(x)^t\eta)\pa_{\px}^\theta e^{i\langle
    \Phi_x(\px),\eta\rangle}}
{\theta!}\bigg|_{x=\px},
\label{7.6.1}
\end{equation}
where
\begin{equation}
\Phi_x(\px)=\phi(\px)-\phi(x)-d\phi(x)(\px-x).
\label{7.6.2}
\end{equation}
Here $\cI_k$ are multi-indices of length $k.$ Our symbols are matrix valued,
e.g. $q_{-2}$ is really $(q_{-2})_{pq}.$ As the change of variables applies
component by component, we suppress these indices in the computations that
follow.

In the case at hand, we are interested in evaluating this expression at $z=x=0,$
where we have $d\phi(0)=\Id$ and
$$\Phi_0(\px)=(-\frac{\alpha}{2}[\langle
  a\pz,\pz\rangle+\Re(b\pz,\pz)+O(|\pz|^3)],0,\dots,0).$$ This is exactly as in
  the K\"ahler case, but for two small modifications: In the K\"ahler case
  $\alpha=1$ and $a=\Id.$ These differences slightly modify the symbolic
  result, but not the invertibility of the symbols of $\cT^{\eo}_{+}.$ As
  before, only the $k=2$ term is of importance. It is given by
\begin{equation}
-\frac{i\xi_1}{2}\Tr[\pa^2_{\xi_j\xi_k} q(0,\xi)\pa^2_{x_jx_k}\phi(0)].
\end{equation}

To compute this term we need to compute the Hessians of $q_{-1}$ and
$\phi(x)$ at $x=0.$ 
We define the $2n\times 2n$ real matrices  $A,B$ so that
\begin{equation}
\langle az,z\rangle=\langle Ax,x\rangle\text{ and }\Re(bz,z)=\langle B x,x\rangle;
\end{equation}
if $a=a^0+ia^1$ and $b=b^0+ib^1,$ then
\begin{equation}
A=\left(\begin{matrix} a^0 & -a^1\\
a^1 & a^0\end{matrix}\right)\quad
B=\left(\begin{matrix} b^0 & -b^1\\
-b^1 & -b^0\end{matrix}\right).
\end{equation}
Here $a^{0t}=a^{0}$ $a^{1t}=-a^{1},$ and $b^{0t}=b^0,$ $b^{1t}=b^1.$
With these definitions we see that
\begin{equation}
\pa^2_{x_jx_k}\phi(0)=-\alpha(A+B).
\end{equation}

As before we compute:
\begin{equation}
\frac{\pa^2q_{-1}}{\pa\xi_k\pa\xi_j}=
-4\frac{d_1\Id+\xi\otimes\pa_\xi d_1^t+\pa_\xi d_1\otimes\xi^t}{|\xi|^4}+
16d_1\frac{\xi\otimes\xi^t}{|\xi|^6}.
\end{equation}
Here $\xi$ and $\pa_\xi d_1$ are regarded as column vectors.
The principal part of the $k=2$ term is
\begin{equation}
q_{-2}^{c}(\xi)=i\xi_1\alpha\Tr\left[(A+B)\left(-2\frac{\Id d_1
+\xi\otimes\pa_\xi
    d_1^t+\pa_\xi 
d_1\otimes\xi^t}{|\xi|^4}+ 8d_1\frac{\xi\otimes\xi^t}{|\xi|^6}\right)\right]
\label{7.6.8}
\end{equation}
Observe that $q_{-2}^c$ depends linearly on $A$ and $B$. It is shown in
Proposition 6 of~\cite{Epstein3} that the contribution, along the contact
direction, of a matrix with the symmetries of $B$ vanishes.  Because $q_{-2}$
vanishes at $0$ and because the order of a symbol is preserved under a change
of variables we see that the symbol of $Q^{\eo}$ at $p$ is 
\begin{equation}
q(0,\xi)=\frac{2d_1(\xi)}{|\xi|^2}+q_{-2}^c(\xi)+\sO_{-3}(1).
\label{7.9.7}
\end{equation}
As before the $\sO_{-3}$-term contributes nothing to the extended Heisenberg
principal symbol of the Calderon projector. Only the term 
\begin{equation}
q_{-2}^{cA}(\xi)=
2i\xi_1\alpha\left[-\frac{\Tr Ad_1}{|\xi|^4}+4\frac{d_1\langle
    A\xi,\xi\rangle}{|\xi|^6}-
2\frac{\langle A\xi,\pa_\xi d_1\rangle}{|\xi|^4}\right].
\label{5.13.2}
\end{equation}
makes a contribution.  To find the contribution of $q_{-2}^{cA}$ to the symbol
of the Calderon projector, we need to compute the contour integral
\begin{equation}
\scp_{-2\pm}^{c}(p,\xi')=
\frac{1}{2\pi}\int\limits_{\Gamma_{\pm}(\xi')}q_{-2}^{cA}(\xi)d\xi_1.
\end{equation}
Let $\xi=(\xi_1,\xi').$ As this term is lower order, in the classical sense, we
only need to compute it for $\xi'$ along the contact line. We do this
computation in the next section.

\section{The symbol of the Calderon projector}\label{s.3}
We are now prepared to compute the symbol of the Calderon projector; it is
expressed as 1-variable contour integral in the symbol of $Q^{\eo}.$ If
$q(t,x',\xi_1,\xi')$ is the symbol of $Q^{\eo}$ in the boundary adapted
coordinates, then the symbol of the Calderon projector is
\begin{equation}
\scp_{\pm}(x',\xi')=\frac{1}{2\pi}\int\limits_{\Gamma_{\pm}(\xi_1)}q(0,x',\xi_1,\xi')d\xi_1
\circ\sigma_1(\eth^{\eo},\mp idt).
\label{7.9.1}
\end{equation}
Here we recall that $q(0,x',\xi_1,\xi')$ is a meromorphic function of $\xi_1.$
For each fixed $\xi',$ the poles of $q$ lie on the imaginary axis. For $t>0,$ we
take $\Gamma_+(\xi_1)$ to be a contour enclosing the poles of
$q(0,x',\cdot,\xi')$ in the upper half plane, for $t<0,$ $\Gamma_-(\xi_1)$ is a
contour enclosing the poles of $q(0,x',\cdot,\xi')$ in the lower half plane. In
a moment we use a residue computation to evaluate these integrals. For this
purpose we note that the contour $\Gamma_+(\xi_1)$ is positively oriented,
while $\Gamma_-(\xi_1)$ is negatively oriented.

The Calderon projector is a classical pseudodifferential operator of order $0$ and
therefore its symbol has an asymptotic expansion of the form
\begin{equation}
\scp_{\pm}=\scp_{0\pm}+\scp_{-1\pm}+\dots
\end{equation}
The contact line, $L_p,$ is defined in $T_p^*Y$ by the equations
\begin{equation}
\xi_2=\dots=\xi_{n}=\xi_{n+2}=\dots=\xi_{2n}=0,
\end{equation}
and $\xi_{n+1}$ is a coordinate along the contact line. Because
$t=-\frac{\alpha}{2}\rho,$ the positive contact direction is given by
$\xi_{n+1}<0.$ As before we have the principal symbols of $\cP^{\eo}_{\pm}$
away from the contact line:
\begin{proposition} If $\tX$ is an invertible double, containing $X$ as an open
  set, and  $p\in bX$ with coordinates normalized at $p$ as above, then
\begin{equation}
\scp^{\eo}_{0\pm}(0,\xi')=
\frac{d_1^{\ooee}(\pm i|\xi'|,\xi')}{|\xi'|}\circ\sigma_1(\eth^{\eo},\mp i dt).
\label{7.9.3}
\end{equation}
\end{proposition}

Along the contact directions we need to evaluate higher order terms; as shown
in~\cite{Epstein3}, the error terms in~\eqref{7.9.7} contribute terms that lift to
have Heisenberg order less than $-2.$
To finish our discussion of the symbol of the Calderon projector we need to
compute the symbol along the contact direction. This entails computing the
contribution from $q_{-2}^{cA}.$ As before, the terms arising from the
holomorphic Hessian of $\rho$ do not contribute anything to the symbol of the
Calderon projector. However, the terms arising from $\pa^{2}_{z_j\bz_k}$ still
need to be computed. To do these computations, we need to have an explicit
formula for the principal symbol $d_1(\xi)$ of $\eth$ at $p.$ For the purposes
of these and our subsequent computations, it is useful to use the chiral
operators $\eth^{\eo}.$ As we are working in a geodesic normal coordinate system, we
only need to find the symbols of $\eth^{\eo}$ for $\bbC^n$ with the flat
metric. Let $\sigma$ denote a section of $\Lambda^{\eo}.$ We split $\sigma$
into its normal and tangential parts at $p:$
\begin{equation}
\sigma=\sigma^t+\frac{d\bz_1}{\sqrt{2}}\wedge\sigma^n,\quad
i_{\pa_{\bz_1}}\sigma^t
=i_{\pa_{\bz_1}}\sigma^n=0.
\label{7.12.1}
\end{equation}

With this splitting we see that
\begin{equation}
\begin{split}
\eth^{\even}\sigma&=\sqrt{2}\left(\begin{matrix}
\pa_{\bz_1}\otimes\Id_{n} & \ccD_t\\
-\ccD_t & -\pa_{z_1}\otimes\Id_{n}\end{matrix}\right)\left(\begin{matrix} \sigma^t\\
\sigma^n\end{matrix}\right)\\
\eth^{\odd}\sigma&=\sqrt{2}\left(\begin{matrix}
-\pa_{z_1}\otimes\Id_{n} & -\ccD_t\\
\ccD_t & \pa_{\bz_1}\otimes\Id_{n}\end{matrix}\right)\left(\begin{matrix} \sigma^n\\
\sigma^t\end{matrix}\right),
\end{split}
\label{7.12.2}
\end{equation}
where $\Id_{n}$ is the identity matrix acting on the normal, or tangential
parts of $\Lambda^{\eo}\restrictedto_{bX}$ and 
\begin{equation}
\ccD_t=\sum_{j=2}^n[\pa_{z_j} e_j-\pa_{\bz_j}\epsilon_j]\text{ with }
e_j= i_{\sqrt{2}\pa_{\bz_j}}\text{ and }
\epsilon_j=\frac{d\bz_j}{\sqrt{2}}\wedge.
\label{7.12.3}
\end{equation}
It is now a simple matter to compute $d_1^{\eo}(\xi):$
\begin{equation}
\begin{split}
&d^{\even}_1(\xi)=\frac{1}{\sqrt{2}}\left(\begin{matrix}
(i\xi_1-\xi_{n+1})\otimes\Id_{n} & \sd(\xi'')\\
-\sd(\xi'') & -(i\xi_1+\xi_{n+1}) \otimes\Id_{n}\end{matrix}\right)\\
&d^{\odd}_1(\xi)=
\frac{1}{\sqrt{2}}\left(\begin{matrix}
-(i\xi_1+\xi_{n+1})\otimes\Id_{n} & -\sd(\xi'')\\
\sd(\xi'') & (i\xi_1-\xi_{n+1})\otimes\Id_{n}\end{matrix}\right)
\end{split}
\label{7.12.4}
\end{equation}
where $\xi''=(\xi_2,\dots,\xi_n,\xi_{n+2},\dots,\xi_{2n})$ and
\begin{equation}
\sd(\xi'')=\sum_{j=2}^n[(i\xi_j+\xi_{n+j})e_j-(i\xi_j-\xi_{n+j})\epsilon_j].
\end{equation}
As $\epsilon_j^*=e_j$ we see that $\sd(\xi'')$ is a self adjoint symbol.

The principal symbols of $\cT^{\eo}_{+}$ have the same block structure as in
the K\"ahler case.  The symbol $q_{-2}^c$ produces a term that lifts to have
Heisenberg order $-2$ and therefore, in the pseudoconvex case, we only need to
compute the $(2,2)$ block arising from this term.

We start with the nontrivial term of order $-1.$
\begin{lemma} If $X$ is either pseudoconvex or pseudoconcave we have that
\begin{equation}
\frac{1}{2\pi}\int\limits_{\Gamma_{\pm}(\xi')}\frac{2i\xi_1\alpha\Tr A
  d_1(\xi_1,\xi')d\xi_1}{|\xi|^4}=-\frac{i\alpha\Tr A\pa_{\xi_1} d_1}{2|\xi'|}
\label{7.12.7}
\end{equation}
\end{lemma}
\begin{remark} As $d_1$ is a linear polynomial, $\pa_{\xi_1}d_1$ is a constant
  matrix. 
\end{remark}
\begin{proof} See Lemma 1 in~\cite{Epstein3}.
\end{proof} 

We complete the computation by evaluating the contribution from the other terms
in $q_{-2}^{cA}$ along the contact line.
\begin{proposition}
For $\xi'$ along
the positive (negative) contact line we have
\begin{equation}
\frac{1}{2\pi}\int\limits_{\Gamma_{\pm}(\xi')}[q_{-2}^{cA}(p,\xi)]d\xi_1=
-\frac{\alpha (a^0_{11}-\ha\Tr A)}{|\xi'|}\pa_{\xi_1} d_1.
\label{7.12.8}
\end{equation}
If $\xi_{n+1}<0,$ then we use $\Gamma_+(\xi'),$ whereas
if $\xi_{n+1}>0,$ then we use $\Gamma_-(\xi').$
\end{proposition}
\begin{proof}
To prove this result we need to evaluate the contour integral with
$$\xi'=\xi'_c=(0,\dots,0,\xi_{n+1},0,\dots,0),$$ recalling that the positive
contact line corresponds to $\xi_{n+1}<0.$ Hence, along the positive contact
line $|\xi'|=-\xi_{n+1}.$ We first compute the integrand along $\xi'_c.$
\begin{lemma} For $\xi'$ along the contact line we have
\begin{equation}
\left[\frac{2d_1^{\even}(\xi)\langle A\xi,\xi\rangle-|\xi|^2
\langle A\xi,\pa_\xi d_1^{\even}\rangle}{|\xi|^6}\right]=
\frac{a^0_{11}}{|\xi|^4}d_1^{\even}(\xi)
\label{7.12.9}
\end{equation}
\begin{equation}
\left[\frac{2d_1^{\odd}(\xi)\langle A\xi,\xi\rangle-|\xi|^2
\langle A\xi,\pa_\xi d_1^{\odd}\rangle}{|\xi|^6}\right]=
\frac{a^0_{11}}{|\xi|^4}d_1^{\odd}(\xi)
\label{7.12.10}
\end{equation}
\end{lemma} 
\begin{proof}
As $a^1_{11}=0$ we observe that along the contact line
\begin{equation}
\begin{split}
&\langle A\xi,\xi\rangle=a^0_{11}(\xi_1^2+\xi_{n+1}^2)\\
&\langle A\xi,\pa_{\xi}d^{\even}_1\rangle=
a^{0}_{11}\left(\begin{matrix}
(i\xi_1-\xi_{n+1})\otimes\Id & 0\\
0&-(i\xi_1+\xi_{n+1})\otimes\Id\end{matrix}\right)=a^0_{11}d_1^{\even}\\
&\langle A\xi,\pa_{\xi}d^{\odd}_1\rangle=
a^{0}_{11}\left(\begin{matrix}
-(i\xi_1+\xi_{n+1})\otimes\Id & 0\\
0&(i\xi_1-\xi_{n+1})\otimes\Id\end{matrix}\right)=a^0_{11}d_1^{\odd}
\end{split}
\end{equation}
As $\xi_1^2+\xi_{n+1}^2=|\xi|^2$ for $\xi'$ along the contact line these
formul{\ae} easily imply~\eqref{7.12.9} and~\eqref{7.12.10}.
\end{proof}

The proposition is an easy consequence of these formul{\ae}.
\end{proof}

For subsequent calculations we set
\begin{equation}
\beta\overset{d}{=}\ha\Tr
A-a^0_{11}=\sum_{j=2}^n[\pa^2_{x_j}\rho+\pa^2_{y_j}\rho]_{x=p}.
\end{equation}
As a corollary, we have a formula for the $-1$ order term in the symbol of the
Calderon projector
\begin{corollary}
In the normalizations defined above, along the contact directions, we have
\begin{equation}
\scp_{-1}^{\eo}(0,\xi')
=-\frac{i\alpha \beta\pa_{\xi_1}d_1^{\ooee}}{|\xi'|}\circ\sigma_1(\eth^{\eo},\mp
idt).
\label{5.13.4}
\end{equation}
\end{corollary}
\begin{remark} In the K\"ahler case $\alpha=1$ and $\beta=n-1.$ The values of
  these numbers turn out to be unimportant. It is only important that
  $\alpha>0$ and that they depend smoothly on local geometric data, which they
  obviously do.
\end{remark}

We have shown that the order $-1$ term in the symbol of the Calderon projector,
along the appropriate half of the contact line, is given by the right hand side
of equation~\eqref{5.13.4}. It is determined by the principal symbol of
$Q^{\eo}$ and does not depend on the higher order geometry of $bX.$ As all
other terms in the symbol of $Q^{\eo}$ contribute terms that lift to have
Heisenberg order less than $-2,$ these computations allow us to find the
principal symbols of $\cT^{\eo}_{+}$ and extend the main results
of~\cite{Epstein3} to the pseudoconvex almost complex category. As noted above,
the off diagonal blocks have Heisenberg order $-1,$ so the classical terms of
order less than zero cannot contribute to their principal parts.

We now give formul\ae\ for the chiral forms of the subelliptic boundary
conditions defined in~\cite{Epstein4} as well as the isomorphisms
$\sigma_1(\eth^{\eo},\mp i dt).$ Let $\cS$ be a generalized Szeg\H o
projector. 
\begin{lemma}\label{lemm4} According to the splittings of sections of $\Lambda^{\eo}$
  given in~\eqref{7.12.1}, the subelliptic boundary conditions, defined by the
  generalized Szeg\H o projector $\cS,$  on even (odd)
  forms are given by $\cR^{\eo}_+\sigma\restrictedto_{bX}=0$ where
\begin{equation}
\cR^{\even}_+\sigma\restrictedto_{bX}=
\left(\begin{matrix}\begin{matrix} \cS & 0 \\
0 & \bzero\\
\end{matrix} &
\begin{matrix} \bzero
\end{matrix}\\
\begin{matrix}\bzero
\end{matrix}&  
\begin{matrix}\Id \end{matrix} \end{matrix}\right)
\left[\begin{matrix} \sigma^{t}\\\phantom{\sigma}\\
\sigma^{n}
\end{matrix}\right]_{bX}
\quad
\cR^{\odd}_+\sigma\restrictedto_{bX}=
\left(\begin{matrix}\begin{matrix} 1-\cS & 0 \\
0 & \Id\\
\end{matrix} &
\begin{matrix} \bzero
\end{matrix}\\
\begin{matrix}\bzero
\end{matrix}&  
\begin{matrix}\bzero \end{matrix} \end{matrix}\right)
\left[\begin{matrix} \sigma^{n}\\\phantom{\sigma}\\
\sigma^{t}
\end{matrix}\right]_{bX}
\end{equation}
\end{lemma}

\begin{lemma}
The isomorphisms at the boundary between $\Lambda^{\eo}$ and 
$\Lambda^{\ooee}$ are given by
\begin{equation}
\sigma_1(\eth^{\eo}_{\pm},\mp i dt)\sigma^t=\frac{\pm }{\sqrt{2}}\sigma^t,\quad
\sigma_1(\eth^{\eo}_{\pm},\mp i dt)\sigma^n=\frac{\mp}{\sqrt{2}}\sigma^n.
\label{7.13.3}
\end{equation}
\end{lemma}

Thus far, we have succeeded in computing the symbols of the Calderon projectors
to high enough order to compute the principal symbols of $\cT^{\eo}_{+}$ as
elements of the extended Heisenberg calculus. The computations have been carried
out in a coordinate system adapted to the boundary. This suffices to examine
the classical parts of the symbols.  Recall that the positive contact direction 
$L^+,$  is given at $p$ by $\xi''=0,\xi_{n+1}<0.$ As before we obtain:
\begin{proposition}\label{prp8} If $(X,J,g,\rho)$ is a normalized strictly
  pseudoconvex $\spnc$-manifold, then, on the complement of
  the positive contact direction, the classical symbols
  ${}^R\sigma_0(\cT^{\eo}_+)$ are given by
\begin{equation}
\begin{split}
{}^R\sigma_0(\cT^{\even}_+)(0,\xi')&=
\frac{1}{2|\xi'|}\left(\begin{matrix} (|\xi'|+\xi_{n+1})\Id & -\sd(\xi'')\\
\sd(\xi'') & (|\xi'|+\xi_{n+1})\Id \end{matrix}\right)\\
{}^R\sigma_0(\cT^{\odd}_+)(0,\xi')&=
\frac{1}{2|\xi'|}\left(\begin{matrix} (|\xi'|+\xi_{n+1})\Id & \sd(\xi'')\\
-\sd(\xi'') & (|\xi'|+\xi_{n+1})\Id \end{matrix}\right)
\end{split}
\label{7.13.1}
\end{equation}
These symbols are invertible on the complement of $L^+.$
\end{proposition}
 \begin{proof}See Proposition 8 in~\cite{Epstein3}.
\end{proof}

\section{The Heisenberg symbols of $\cT^{\eo}_{+}$}
To compute the Heisenberg symbols of $\cT^{\eo}_{+}$ we change coordinates, one
last time, to get Darboux coordinates at $p.$ Up to this point we have used the
coordinates $(\xi_2,\dots,\xi_{2n})$ for $T^*_pbX,$ which are defined by the
coframe $dx_2,\dots,dx_{2n},$ with $dx_{n+1}$ the contact direction. Recall
that the contact form $\theta,$ defined by the complex structure and defining
function $\rho,$ is given by $\theta=\frac{i}{2}\dbar\rho.$ The symplectic form
on $H$ is defined by $d\theta.$ At $p$ we have
\begin{equation}
\theta_p=-\frac{1}{2\alpha}dx_{n+1},\quad d\theta_p=\sum_{j=2}^{n}dx_j\wedge
dx_{j+n}.
\label{7.14.1}
\end{equation}
By comparison with equation (5) in~\cite{Epstein3}, we see that properly
normalized coordinates for $T^*_pbX$ (i.e., Darboux coordinates) are obtained
by setting 
\begin{equation}
\eta_0=-2\alpha\xi_{n+1},\quad \eta_j=\xi_{j+1},\quad
\eta_{j+n-1}=\xi_{j+n+1}\text{ for }j=1,\dots,n-1.
\end{equation}
As usual we let $\eta'=(\eta_1,\dots,\eta_{2(n-1)});$ whence $\xi''=\eta'.$

As a first step in lifting the symbols of the Calderon projectors to the
extended Heisenberg compactification, we re-express them, through order $-1$ in
the $\xi$-coordinates:
\begin{equation}
\scp^{\even}_+(\xi')=
\frac{1}{2|\xi'|}\left[
\left(\begin{matrix} (|\xi'|-\xi_{n+1})\Id & \sd(\xi'')\\\sd(\xi'')&
(|\xi'|+\xi_{n+1})\Id\end{matrix}\right)-\alpha\beta\left(\begin{matrix}\Id &\bzero\\
\bzero&\Id\end{matrix}\right)\right]
\label{7.14.5}
\end{equation}

\begin{equation}
\scp^{\odd}_+(\xi')=
\frac{1}{2|\xi'|}\left[
\left(\begin{matrix} (|\xi'|+\xi_{n+1})\Id & \sd(\xi'')\\\sd(\xi'')&
(|\xi'|-\xi_{n+1})\Id\end{matrix}\right)-\alpha\beta\left(\begin{matrix}\Id &\bzero\\
\bzero&\Id\end{matrix}\right)\right]
\label{7.14.6}
\end{equation}

Various identity and zero matrices appear in these symbolic
computations. Precisely which matrix is needed depends on the dimension, the
parity, etc. We do not encumber our notation with these distinctions.

In order to compute ${}^H\sigma(\cT^{\eo}_{+}),$ we represent the Heisenberg
symbols as model operators and use operator composition.  To that end we need
to quantize $\sd(\eta')$ as well as the terms coming from the diagonals
in~\eqref{7.14.5}--\eqref{7.14.6}.  For the pseudoconvex side, we need to
consider the symbols on positive Heisenberg face, where the function
$|\xi'|+\xi_{n+1}$ vanishes.

We express the various terms in the symbol $\scp^{\eo}_{+},$ near the positive contact line
as sums of Heisenberg homogeneous terms
\begin{equation}
\begin{split}
|\xi'|=&\frac{\eta_0}{2\alpha}(1+\sO_{-2}^H)\\
 |\xi'|-\xi_{n+1}=\frac{\eta_0}{\alpha}(1+\sO_{-2}^H),\quad&
|\xi'|+\xi_{n+1}=\frac{\alpha|\eta'|^2}{\eta_0}(1+\sO_{-2}^H)\\
\sd(\xi'')=\sum_{j=1}^{n-1}[(i\eta_j+\eta_{n+j-1})e_j&-(i\eta_j-\eta_{n+j-1})\epsilon_j].
\end{split}
\label{7.14.9}
\end{equation}
Recall that the notation $\sO_j^H$ denotes a term of Heisenberg order at most
$j.$ To find the model operators, we use the
quantization rule, equation (20) in~\cite{Epstein3} (with the $+$ sign), obtaining
\begin{equation}
\begin{split}
\eta_j-i\eta_{n+j-1}&\leftrightarrow C_j\overset{d}{=}(w_j-\pa_{w_j})\\
\eta_j+i\eta_{n+j-1}&\leftrightarrow C_j^*\overset{d}{=}(w_j+\pa_{w_j})\\
|\eta'|^2&\leftrightarrow \ho\overset{d}{=}\sum_{j=1}^{n-1}w_j^2-\pa_{w_j}^2.
\end{split}
\label{7.14.10}
\end{equation}
The following standard identities are  useful
\begin{equation}
\sum_{j=1}^{n-1}C_j^*C_j-(n-1)=\ho=\sum_{j=1}^{n-1}C_jC_j^*+(n-1)
\label{7.14.11}
\end{equation}
We let $\cD_+$ denote the model operator defined, using the $+$ quantization,
by $\sd(\xi''),$ it is given by
\begin{equation}
\cD_+=i\sum_{j=1}^{n-1}[C_j e_j-C_j^*\epsilon_j].
\end{equation}
This operator can be split into even and odd parts, $\cD^{\eo}_+$ and these
chiral forms of the operator are what appear in the model operators below.

With these preliminaries, we can compute the model operators for $\cP^{\even}_+$
and $\Id-\cP^{\even}_+$ in the positive contact direction. They are:
\begin{equation}
{}^{eH}\sigma(\cP^{\even}_+)(+)=\left(\begin{matrix} \Id
  &\frac{\alpha\cD_+^{\odd}}{\eta_0}
\\\frac{\alpha\cD_+^{\even}}{\eta_0}
  &\frac{\alpha^2\ho-\alpha^2\beta\eta_0}{\eta_0^2}\end{matrix}\right)\quad
{}^{eH}\sigma(\Id-\cP^{\even}_+)(+)=\left(\begin{matrix} \frac{\alpha^2\ho+\alpha^2\beta\eta_0}{\eta_0^2} 
&-\frac{\alpha\cD_+^{\odd}}{\eta_0}\\-\frac{\alpha\cD_+^{\even}}{\eta_0} &\Id\end{matrix}\right).
\label{7.14.12}
\end{equation}
The denominators involving $\eta_0$ are meant to remind the reader of the
Heisenberg orders of the various blocks: $\eta_0^{-1}$ indicates a term of
Heisenberg order $-1$ and $\eta_0^{-2}$ a term of order $-2.$ Similar
computations give the model operators in the odd case:
\begin{equation}
{}^{eH}\sigma(\cP^{\odd}_+)(+)=\left(\begin{matrix} \frac{\alpha^2\ho-\alpha^2\beta\eta_0}{\eta_0^2}
  &\frac{\alpha\cD_+^{\odd}}{\eta_0}\\\frac{\alpha\cD_+^{\even}}{\eta_0} 
  &\Id\end{matrix}\right)\quad
{}^{eH}\sigma(\Id-\cP^{\odd}_+)(+)=\left(\begin{matrix} \Id
&-\frac{\alpha\cD_+^{\odd}}{\eta_0}\\-\frac{\alpha\cD_+^{\even}}{\eta_0}
  &\frac{\alpha^2\ho+\alpha^2\beta\eta_0}{\eta_0^2} \end{matrix}\right). 
\label{7.14.13}
\end{equation}
Let $\sypr_0'={}^{eH}\sigma(+)(\cS);$ it is a self adjoint rank one projection
defined by a compatible almost complex structure on $H.$ We use the ${}'$ to
distinguish this rank one projection, from the rank one projection $\pi_0$
defined by the CR-structure on the fiber of cotangent bundle at $p.$ The model
operators for $\cR_+^{\eo}$ in the positive contact direction are:
\begin{equation}
{}^{eH}\sigma(\cR^{\even}_{+})(+)=
\left(\begin{matrix}\begin{matrix} \sypr'_0 & 0 \\
0 & \bzero\\
\end{matrix} &
\begin{matrix} \bzero
\end{matrix}\\
\begin{matrix}\bzero
\end{matrix}&  
\begin{matrix}\Id \end{matrix} \end{matrix}\right),\quad
{}^{eH}\sigma(\cR_+^{\odd})(+)=
\left(\begin{matrix}\begin{matrix} 1-\sypr'_0 & 0 \\
0 & \Id
\end{matrix} &
\begin{matrix} \bzero
\end{matrix}\\
\begin{matrix}\bzero
\end{matrix}&  
\begin{matrix}\bzero \end{matrix} \end{matrix}\right).
\label{5.17.17}
\end{equation}
We can now compute the model operators for $\cT^{\eo}_+$
on the upper Heisenberg face.
\begin{proposition} If $(X,J,g,\rho)$ is normalized strictly pseudoconvex
  $\spnc$-manifold, then, at $p\in bX,$ the model operators for
  $\cT^{\eo}_+,$ in the positive contact direction, are given by
\begin{equation}
{}^{eH}\sigma(\cT^{\even}_+)(+)=
\left(\begin{matrix}\begin{matrix} \sypr'_0 & 0 \\
0 & \bzero\\
\end{matrix} &\begin{matrix}-\left[\begin{matrix} 1-2\sypr'_0 & 0 \\
0 &\Id
\end{matrix}\right]\frac{\alpha\cD_+^{\odd}}{\eta_0}
\end{matrix}
\\
\begin{matrix} \frac{\alpha\cD_+^{\even}}{\eta_0}
\end{matrix}&  
\begin{matrix} \frac{\alpha^2\ho-\alpha^2\beta\eta_0}{\eta_0^2} \end{matrix} \end{matrix}\right)
\label{7.14.21}
\end{equation}
\begin{equation}
{}^{eH}\sigma(\cT^{\odd}_+)(+)=
\left(\begin{matrix}\begin{matrix} \sypr'_0 & 0 \\
0 & \bzero\\
\end{matrix} &\begin{matrix}\left[\begin{matrix} 1-2\sypr'_0 & 0 \\
0 &\Id
\end{matrix}\right]\frac{\alpha\cD_+^{\odd}}{\eta_0}
\end{matrix}
\\
\begin{matrix} -\frac{\alpha\cD_+^{\even}}{\eta_0}
\end{matrix}&  
\begin{matrix} \frac{\alpha^2\ho+\alpha^2\beta\eta_0}{\eta_0^2} \end{matrix} \end{matrix}\right).
\label{7.14.22}
\end{equation}
\end{proposition}
\begin{proof} Observe that the Heisenberg orders of the blocks
  in~\eqref{7.14.21} and~\eqref{7.14.22} are
\begin{equation}
\left(\begin{matrix} 0 & -1 \\-1 & -2\end{matrix}\right).
\end{equation}
Proposition 6 in~\cite{Epstein3} shows that all other terms in the symbol of
the Calderon projector lead to diagonal terms of Heisenberg order at most $-4,$
and off diagonal terms of order at most $-2.$ This, along with the computations
above, completes the proof of the proposition.
\end{proof}

This brings us to the generalization, in the non-K\"ahler case, of Theorem 1
in~\cite{Epstein3}:
\begin{theorem}\label{thm1} Let $(X,J,g,\rho)$ be a normalized strictly pseudoconvex
  $\spnc$-manifold, and $\cS$ a generalized Szeg\H o projector, defined by a
  compatible deformation of the almost complex structure on $H$ induced by the
  embedding of $bX$ as the boundary of $X.$ The comparison operators,
  $\cT^{\eo}_{+},$ are elliptic elements of the extended Heisenberg calculus,
  with parametrices having Heisenberg orders
\begin{equation}
\left(\begin{matrix} 0 & 1\\ 1 & 1\end{matrix}\right).
\end{equation}
\end{theorem}
\begin{proof} The proof is identical to the proof of Theorem 1
  in~\cite{Epstein3}: we need to show that the principal symbols of $\cT^{\eo}_{+}$
  are invertible, which is done in the next section.
\end{proof}

\section{Invertibility of the model operators} 
In this section we produce inverses for the model operators
$\sigma^H(\cT^{\eo}_+)(+).$ We begin by writing down inverses for the model
operators using the projector compatible with the CR-structure induced at $p$
by $J.$ We denote this projector by $\pi_0$ to distinguish it from
$\pi_0'=\sigma_0^H(p)(\cS).$ \emph{In this section,} we let
$\esym{eH}{}(\cT^{\eo}_{+})(+),$ denote the model operators with this projector
to distinguish it from $\esym{eH}{}(\cT^{\prime\eo}_{+})(+),$ the model
operators with $\pi_0'.$ As before, the inverse in the general case is a finite
rank perturbation of this case. For the computations in this section we recall
that $\alpha$ is a positive number.

The operators  $\{C_j\}$ are called the creation operators and the operators
$\{C_j^*\}$ the annihilation operators.  They satisfy the commutation relations
\begin{equation}
[C_j,C_k]=[C_j^*,C_k^*]=0,\quad
[C_j,C_k^*]=-2\delta_{jk}
\label{7.16.4}
\end{equation}
The operators $\cD_{\pm}$ act on sums of the form
\begin{equation}
\omega=\sum_{k=0}^{n-1}\sum_{I\in\cI_k'} f_I\bomega^I,
\end{equation}
here $\cI_k'$ are increasing multi-indices of length $k.$ We refer to the terms
with $|I|=k$ as the terms of degree $k.$ For an increasing $k$-multi-index
$I=1\leq i_1<i_2<\dots<i_k\leq n-1,$ $\bomega^I$ is defined by
\begin{equation}
\bomega^I=\frac{1}{2^{\frac k2}}d\bz_{i_1}\wedge\dots\wedge d\bz_{i_k}.
\end{equation}

The projector $\pi_0$ and the operator $\cD_+$ satisfy the following relation:
\begin{lemma}\label{lem7} Let $\sypr_0$  be the symbol of the
  generalized Szeg\H o projector compatible with the CR-structure defined on
  the fiber of $T_pbX$ by the almost complex structure, then
\begin{equation}
\left[\begin{matrix} \sypr_0 & 0 \\
0 &\bzero
\end{matrix}\right]\cD^{\odd}_+=0
\end{equation}
\end{lemma}
\begin{proof} See Lemma 7 in~\cite{Epstein3}
\end{proof}

This lemma  simplifies the analysis of the  model operators for
$\cT^{\eo}_{+}.$  The following lemma is useful in finding their inverses.
\begin{lemma} Let $\Pi_q$ denote projection onto the terms of degree $q,$
\begin{equation}
\Pi_q\omega=\sum_{I\in\cI_q'} f_I\bomega^I.
\end{equation}
The operators $\cD_{+}$ satisfies the identity
\begin{equation}
\cD_+^2=\sum_{j=1}^{n-1}C_jC_j^*\otimes\Id+\sum_{q=0}^{n-1}2q\Pi_q
\label{7.16.9}
\end{equation}
\end{lemma}
\begin{proof}
See Lemma 9 in~\cite{Epstein3}.
\end{proof}

As before
$\esym{eH}{}(\cT^{\eo}_{+})(+)$ are Fredholm elements (in the graded
sense), in the isotropic algebra. Notice that this is a purely symbolic
statement in the isotropic algebra. The blocks have isotropic orders
\begin{equation}
\left(\begin{matrix} 0 & 1\\ 1& 2\end{matrix} \right).
\end{equation}
The leading order part in the isotropic algebra is
independent of the choice of generalized Szeg\H o projector.  In the former
case we can think of the operator as defining a map from
$H^1(\bbR^{n-1};E_1)\oplus H^{2}(\bbR^{n-1}; E_2)$ to
$H^1(\bbR^{n-1};F_1)\oplus H^{0}(\bbR^{n-1}; F_2)$ for appropriate vector
bundles $E_1, E_2, F_1, F_2.$ It is as maps between these spaces that the model
operators are Fredholm.
\begin{proposition}\label{prop12} The model operators,
  $\esym{eH}{}(\cT^{\eo}_{+})(+),$ are graded Fredholm elements in the
  isotropic algebra.
\end{proposition}
\begin{proof} See Proposition 7 in~\cite{Epstein3}.
\end{proof}

The operators $\cD_{+}^{\even}$ and $\cD_{+}^{\odd}$ are adjoint to one
another. From~\eqref{7.16.9} and the well known properties of the harmonic
oscillator, it is clear that $\cD_+^{\even}\cD_{+}^{\odd}$ is invertible. As
$\cD_+^{\even}$ has a one dimensional null space this easily implies that
$\cD_+^{\odd}$ is injective with image orthogonal to the range of $\sypr_0,$
while $\cD_+^{\even}$ is surjective. With these observations we easily invert
the model operators.

Let $[\cD_+^{\even}]^{-1}u$ denote the unique solution
to the equation
$$\cD_+^{\even} v=u,$$
orthogonal to the null space of $\cD_+^{\even}.$ We let 
\begin{equation}
\sh{u}=\left(\begin{matrix} 1-\sypr_0 & 0\\ 0& \Id\end{matrix}\right)u;
\end{equation}
this is the projection onto the range of $\cD_+^{\odd}$ and
\begin{equation}
u_0= \left(\begin{matrix} \sypr_0 & 0\\ 0& \bzero\end{matrix}\right)u,
\end{equation}
denotes the projection onto the nullspace of $\cD_+^{\even}.$ We let
$[\cD_+^{\odd}]^{-1}$ denote the unique solution to
$$\cD_+^{\odd}v=\sh{u}.$$ 
Proposition~\ref{prop12} shows that these
partial inverses are isotropic operators of order $-1.$

With this notation we find the inverse of $\esym{eH}{}(\cT^{\even}_+)(+).$
The vector $[u,v]$ satisfies
\begin{equation}
\esym{eH}{}(\cT^{\even}_+)(+)\left[\begin{matrix} u\\ v
\end{matrix}\right]=\left[\begin{matrix} a\\ b
\end{matrix}\right]
\end{equation}
if and only if
\begin{equation}
\begin{split}
u&=a_0+[\alpha\cD_+^{\even}]^{-1}(\alpha^2\ho-\alpha^2\beta)[\alpha\cD_+^{\odd}]^{-1}
\sh{a}+[\alpha\cD_+^{\even}]^{-1} b\\
v&=-[\alpha\cD_+^{\odd}]^{-1}\sh{a}.
\end{split}
\label{7.19.13}
\end{equation}
Writing out the inverse as a block matrix of operators, with appropriate
factors of $\eta_0$ included, gives:
\begin{equation}
\begin{split}
[\esym{eH}{}&(\cT^{\even}_+)(+)]^{-1}=\\
&\left[\begin{matrix} \left(\begin{matrix} \sypr_0
      &0\\0&\bzero\end{matrix}\right)
+[\cD_+^{\even}]^{-1}(\ho-\beta)[\cD_+^{\odd}]^{-1}
\left(\begin{matrix}1- \sypr_0
      &0\\0&\Id\end{matrix}\right) & \eta_0[\alpha\cD_+^{\even}]^{-1}\\
-\eta_0[\alpha\cD_+^{\odd}]^{-1}\left(\begin{matrix}1- \sypr_0
      &0\\0&\Id\end{matrix}\right) & \bzero\end{matrix}\right]
\end{split}
\end{equation}
The isotropic operators $[\alpha\cD_{+}^{\eo}]^{-1}$ are of order $-1,$ whereas
the operator,
$$[\cD_+^{\even}]^{-1}(\ho-\beta)[\cD_+^{\odd}]^{-1},$$
is of order zero. The Schwartz kernel of $\sypr_0$ is rapidly decreasing. From
this we conclude that the Heisenberg orders, as a block matrix, of the
parametrix for $[\esym{eH}{}(\cT^{\even}_+)(+)]$ are
\begin{equation}
\left(\begin{matrix} 0 & 1 \\ 1&1\end{matrix}\right).
\label{7.16.10}
\end{equation}
We get a $1$ in the lower right corner because the principal symbol of this
entry, a priori of order $2,$ vanishes. As a result, the inverses of the model
operators have Heisenberg order at most $1,$ which in turn allows us to use
this representation of the parametrix to deduce the standard subelliptic
$\ha$-estimates for these boundary value problems.

The solution for the odd case is given by
\begin{equation}
\begin{split}
u&=a_0+[\cD_+^{\even}]^{-1}(\ho+\beta)[\cD_+^{\odd}]^{-1}
\sh{a}-[\alpha\cD_+^{\even}]^{-1} b\\
v&=[\alpha\cD_+^{\odd}]^{-1}\sh{a}.
\end{split}
\end{equation}
Once again the $(2,2)$ block of  $[\esym{eH}{}(\cT^{\odd}_+)(+)]^{-1}$ vanishes,
and the principal symbol has the Heisenberg orders indicated
in~\eqref{7.16.10}.

For the case that $\pi_0'=\pi_0,$ Lemma~\ref{lem7} implies that
the model operators satisfy
\begin{equation}
[\esym{eH}{}(\cT^{\eo}_{+})(+)]^*=\esym{eH}{}(\cT^{\ooee}_{+})(+).
\end{equation}
From Proposition~\ref{prop12}, we know that these are Fredholm operators. Since
we have shown that all the operators $\esym{eH}{}(\cT^{\eo}_{+})(+)$ are
surjective, i.e., have a left inverse, it follows that all are in fact
injective and therefore invertible.  In all cases this completes the proof of
Theorem~\ref{thm1} in the special case that the principal symbol of $\cS$
equals $\pi_0.$

We now show that the parametrices for $\esym{eH}{}(\cT^{\prime\eo}_{+})(+)$ differ
from those with classical Szeg\H o projectors by operators of finite rank. The
Schwartz kernels of the correction terms are in the Hermite ideal, and so do
not affect the Heisenberg orders of the blocks in the parametrix. As before the
principal symbol in the $(2,2)$ block vanishes.

With these preliminaries and the results from the beginning of Section 7
in~\cite{Epstein3}, we can now complete the proof of Theorem~\ref{thm1}. As
noted above,  $\esym{eH}{}(\cT^{\eo}_{+})(+)$ denotes the
model operators with the projector $\pi_0,$ and
$\esym{eH}{}(\cT^{\prime\eo}_{+})(+)$  the model operators with projector $\pi_0'.$
\begin{proposition} If $\sypr_0'$  is the principal symbol of a generalized  Szeg\H o
  projection, which is a deformation of $\sypr_0,$ then
  $\esym{eH}{}(\cT^{\prime\eo}_{+})(+)$ are invertible elements of the
  isotropic algebra. The inverses satisfy
\begin{equation}
[\esym{eH}{}(\cT^{\prime\eo}_{+})(+)]^{-1}=[\esym{eH}{}(\cT^{\eo}_{+})(+)]^{-1}+
\left(\begin{matrix} c_{1} & c_2\\ c_3 & 0\end{matrix}\right).
\label{7.19.8}
\end{equation}
Here $c_1, c_2, c_3$ are finite rank operators in the Hermite ideal.
\end{proposition}
\begin{proof}
In the formul\ae\  below we let $z_0$ denote the unit vector spanning the
range of $\sypr_0$ and $z_0',$ the unit vector spanning the range of
$\sypr_0'.$

Proposition~\ref{prop12} implies that $\esym{eH}{}(\cT^{\prime\eo}_{+})(+)$ are
Fredholm operators. Since, as isotropic operators, the differences
$$\esym{eH}{}(\cT^{\prime\eo}_{+})(+)-\esym{eH}{}(\cT^{\eo}_{+})(+)$$
are finite rank operators, it follows that
$\esym{eH}{}(\cT^{\prime\eo}_{+})(+)$ have index zero. It therefore
suffices to construct a left inverse.

We begin with the $+$ even case by rewriting the equation
\begin{equation}
\esym{eH}{}(\cT^{\prime\even}_{+})(+)\left[\begin{matrix} u\\
    v\end{matrix}\right]
=\left[\begin{matrix} a\\ b\end{matrix}\right],
\end{equation}
as
\begin{equation}
\begin{split}
\left[\begin{matrix}\sypr_0' & 0\\0 &\bzero\end{matrix}\right][u+\alpha\cD_+^{\odd} v]&=
\left[\begin{matrix}\sypr_0' & 0\\0 &\bzero\end{matrix}\right]a\\
\left[\begin{matrix}1-\sypr_0' & 0\\0 &\Id\end{matrix}\right]\alpha\cD_+^{\odd} v&=
-\left[\begin{matrix}1-\sypr_0' & 0\\0 &\Id\end{matrix}\right]a\\
\alpha\cD_+^{\even} u + (\alpha^2\ho-\alpha^2\beta)v&=b.
\end{split}
\label{7.19.11}
\end{equation}
We solve the middle equation in~\eqref{7.19.11} first. Let
\begin{equation}
A_1=(\frac{z_0'\otimes z_0^t}{\langle z_0',z_0\rangle}-\sypr_0)\Pi_0 a,
\label{7.19.15}
\end{equation}
and note that $\sypr_0A_1=0.$ Corollary 2 in~\cite{Epstein3} shows that the
model operator in~\eqref{7.19.15} provides a globally defined symbol. The
section $v$ is determined as the unique solution to
\begin{equation}
\alpha\cD_+^{\odd}v=-(\sh{a}-A_1).
\end{equation}
By construction $(1-\sypr_0')(a_0+A_1)=0$ and therefore the second
equation is solved. The section $\sh{u}$ is now uniquely determined by the last
equation in~\eqref{7.19.11}:
\begin{equation}
\sh{u}=[\alpha\cD_{+}^{\even}]^{-1}(b+(\alpha^2\ho-\alpha^2\beta)[\alpha\cD_+^{\odd}]^{-1}(\sh{a}-A_1)).
\end{equation}
This leaves only the first equation, which we rewrite as
\begin{equation}
\left[\begin{matrix}\sypr_0' & 0\\0 &\bzero\end{matrix}\right]u_0=
\left[\begin{matrix}\sypr_0' & 0\\0 &\bzero\end{matrix}\right](a-\alpha\cD_+^{\odd}
v-\sh{u}).
\end{equation}
It is immediate that
\begin{equation}
u_0=\frac{z_0\otimes z_0^{\prime t}}{\langle z_0,z_0'\rangle}\Pi_0(a-\alpha\cD_+^{\odd}
v-\sh{u}).
\end{equation}
By comparing these equations to those in~\eqref{7.19.13} we see that
$[\esym{eH}{}(\cT^{\prime\even}_{+})(+)]^{-1}$ has the required form. The
finite rank operators are finite sums of terms involving $\sypr_0,$ $z_0\otimes
z_0^{\prime t}$ and $z_0^{\prime t}\otimes z_0,$ and are therefore in the
Hermite ideal.

The solution in the $+$ odd case is given by
\begin{equation}
\begin{split}
v&=[\alpha\cD_+^{\odd}]^{-1}(\sh{a}-A_1)\\
\sh{u} &=[\alpha\cD_+^{\even}]^{-1}[(\alpha^2\ho+\alpha^2\beta)v)-b]\\
u_0&=\frac{z_0\otimes z_0^{\prime t}}{\langle z_0,z_0'\rangle}\Pi_0(a+\alpha\cD_+^{\odd}
v-\sh{u})
\end{split}
\label{7.19.14}
\end{equation}
As before $A_1$ is given by~\eqref{7.19.15}. Again the inverse of
$\esym{eH}{}(\cT^{\prime\odd}_{+})(+)$ has the desired form.
\end{proof}

As noted above, the operators $\esym{eH}{}(\cT^{\prime\eo}_{+})(+)$ are
Fredholm operators of index zero. Hence, solvability of the equations
\begin{equation}
\esym{eH}{}(\cT^{\prime\eo}_{+})(+) \left[\begin{matrix} u\\
    v\end{matrix}\right]
=\left[\begin{matrix} a\\ b\end{matrix}\right],
\end{equation}
for all $[ a, b]$ implies the uniqueness and therefore the invertibility of the
model operators. This completes the proof of Theorem~\ref{thm1}. We now turn to
applications of these results.

\begin{remark} For the remainder of the paper $\cT^{\eo}_+$ is used to denote
  the comparison operator defined by $\cR_+^{\eo},$ where the rank one
  projections are given by the principal symbol of $\cS.$
\end{remark}

\section{Consequences of Ellipticity}
As in the K\"ahler case, the ellipticity of the operators $\cT^{\eo}_+$
implies that the graph closures of $(\eth^{\eo}_+,\cR_+^{\eo})$ are
Fredholm and moreover,
\begin{equation}
(\eth^{\eo}_+,\cR_+^{\eo})^*=\overline{(\eth^{\ooee}_+,\cR_+^{\ooee})}
\end{equation}
Given the ellipticity of $\cT^{\eo}_+,$ the proofs of these statements
are identical to the proofs in the K\"ahler case. For later usage, we
introduce some notation and state these results.

Let $\cU^{\eo}_{+}$ denote a 2-sided parametrix
defined so that
\begin{equation}
\begin{split}
&\cT^{\eo}_{+}\cU^{\eo}_{+}=\Id-K^{\eo}_1\\
&\cU^{\eo}_{+}\cT^{\eo}_{+}=\Id-K^{\eo}_2,
\end{split}
\label{7.22.1}
\end{equation}
with $K^{\eo}_1, K^{\eo}_2$ finite rank smoothing operators.  The principal symbol
computations show that $\cU^{\eo}_{+}$ has classical order $0$ and
Heisenberg order at most $1.$  Such an operator defines a bounded map from
$H^{\ha}(bX)$ to $L^2(bX).$
\begin{proposition}\label{prop15} The operators $\cU^{\eo}_{+}$ define bounded maps
  $$\cU^{\eo}_{+}:H^{s}(bX;F)\to H^{s-\ha}(bX;F)$$ 
for $s\in\bbR.$ Here $F$ is an appropriate vector bundle over $bX.$
\end{proposition}

The mapping properties of the boundary parametrices allow us to show that the
graph closures of the operators $(\eth^{\eo}_{+},\cR^{\eo}_{+})$ are
Fredholm. 
\begin{theorem}\label{thm2}  Let $(X,J,g,\rho)$ define a normalized strictly
  pseudoconvex $\spnc$-manifold.  The graph closures of
  $(\eth^{\eo}_{+},\cR^{\eo}_{+}),$ are Fredholm operators.
\end{theorem}
\begin{proof} The proof is exactly the same as the proof of Theorem 2
  in~\cite{Epstein3}. 
\end{proof}

We also obtain the standard subelliptic Sobolev space estimates for the
operators $(\eth^{\eo}_{+},\cR^{\eo}_{+}).$
\begin{theorem}Let $(X,J,g,\rho)$ define a normalized strictly
  pseudoconvex $\spnc$-manifold. \label{thm33} 
   For each $s\geq 0,$
  there is a positive constant $C_s$ such that if $u$ is an $L^2$-solution to
$$\eth^{\eo}_{+} u = f\in H^s(X)\text{ and }\cR^{\eo}_{+}[u]_{bX}=0$$
in the sense of distributions, then
\begin{equation}
\|u\|_{H^{s+\ha}}\leq C_s[\|\eth^{\eo}_{+} u\|_{H^s}+\|u\|_{L^2}].
\label{7.22.100}
\end{equation}
\end{theorem}
\begin{proof} Exactly as in the K\"ahler case.
\end{proof}
\begin{remark} In the case $s=0,$ there is a slightly better result: the
  Poisson kernel maps $L^2(bX)$ into $H_{\hn}(X)$ and therefore
  the argument shows that there is a constant $C_0$ such that if $u\in L^2,$
  $\eth^{\eo}_{+} u\in L^2$ and $\cR^{\eo}_{+} [u]_{bX}=0,$ then
\begin{equation}
\|u\|_{\hn}\leq C_0[\|f\|_{L^2}+\|u\|_{L^2}]
\end{equation}
This is just  the standard $\ha$-estimate for the operators
$\overline{(\eth^{\eo}_{+},\cR^{\eo}_{+})}$ 
\end{remark}

It is also possible to prove localized versions of these results. The higher
norm estimates have the same consequences as for the $\dbar$-Neumann
problem. Indeed, under certain hypotheses these estimates imply higher norm
estimates for the second order operators considered in~\cite{Epstein4}.  We
identify the adjoints:
\begin{theorem}\label{thm3} Let $(X,J,g,\rho)$ define a normalized strictly
  pseudoconvex $\spnc$-manifold, then we have
  the following relations:
\begin{equation}
(\eth^{\eo}_{+},\cR^{\eo}_{+})^*=
\overline{(\eth^{\ooee}_{+},\cR^{\ooee}_{+})}.
\label{7.21.2}
\end{equation}
\end{theorem}

As a corollary of Theorem~\ref{thm3}, we get estimates for the second order operators
$\eth^{\ooee}_{+}\eth^{\eo}_{+},$ with subelliptic boundary conditions.
\begin{corollary}\label{cor3}  Let $(X,J,g,\rho)$ define a normalized strictly
  pseudoconvex $\spnc$-manifold. For $s\geq 0$ there
  exist constants $C_s$ such that if $u\in L^2,$ $\eth^{\eo}_{+}u\in L^2,$
  $\eth^{\ooee}_{+}\eth^{\eo}_{+}u\in H^s$ and $\cR^{\eo}_{+} [u]_{bX}=0,
  \cR^{\ooee}_{+}[\eth^{\eo}_{+}u]=0$ in the sense of distributions, then
\begin{equation}
\|u\|_{H^{s+1}}\leq
C_s[\|\eth^{\ooee}_{+}\eth^{\eo}_{+}u\|_{H^s}+\|u\|_{L^2}].
\label{7.23.3}
\end{equation}
\end{corollary}

We close this section by considering $(\cP^{\eo}_+,\cR^{\eo}_+)$ as a tame
Fredholm pair, as defined in the appendix.  To apply the the functional
analytic framework set up in the appendix, we use as the family of separable
Hilbert spaces the $L^2$-Sobolev spaces $H^s(bX;F),$ where $F$ are appropriate
vector bundles. The norms on these spaces can be selected to satisfy the
conditions,~\eqref{eqn5.16.1} and~\eqref{eqna1}.  In this setting the algebra
of tame operators certainly includes the extended Heisenberg calculus. In this
setting the smoothing operators are operators in
${}^{eH}\Psi^{-\infty,-\infty,-\infty}(bX;F,G),$ i.e., operators from sections
of $F$ to sections of $G$ (two vector bundles) with a Schwartz kernel in
$\CI(bX\times bX).$

An immediate corollary of Theorem~\ref{thm1} is:
\begin{corollary} Let $(X,J,g,\rho)$ define a normalized strictly
  pseudoconvex $\spnc$-manifold, and let $\cP^{\eo}_+$ be the Calderon
  projectors for $\eth^{\eo}_+.$ If $\cR^{\eo}_+$ are projectors defining
  modified $\dbar$-Neumann boundary conditions, then
  $(\cP^{\eo}_+,\cR^{\eo}_+)$ are tame Fredholm pairs.
\end{corollary}
If $\cU^{\eo}_{+}$ are parametrices for $\cT^{\eo}_+,$ and $K^{\eo}_1,$
$K^{\eo}_2$ are smoothing operators that satisfy
\begin{equation}
\cT^{\eo}_+\cU^{\eo}_+=\Id-K^{\eo}_1,\quad \cU^{\eo}_+\cT^{\eo}_+=\Id -K^{\eo}_2,
\end{equation}
then Theorem~\ref{thma1} immediately implies:
\begin{theorem}  Let $(X,J,g,\rho)$ define a normalized strictly
  pseudoconvex $\spnc$-manifold, and let $\cP^{\eo}_+$ be  Calderon
  projectors for $\eth^{\eo}_+.$ If $\cR^{\eo}_+$ are projectors defining
  a modified $\dbar$-Neumann boundary conditions, then
\begin{equation}
\Rind(\cP^{\eo}_+,\cR^{\eo}_+)=\Tr_{L^2}(\cP^{\eo}_+K^{\eo}_2\cP^{\eo}_+)-
\Tr_{L^2}(\cR^{\eo}_+K^{\eo}_1\cR^{\eo}_+)
\label{5.16.4}
\end{equation}
\end{theorem}
If $\cP^{\eo}_+K^{\eo}_2\cP^{\eo}_+$ have Schwartz kernels $\kappa^{\eo}_2(x,y)$
and $\cR^{\eo}_+K^{\eo}_1\cR^{\eo}_+$ have Schwartz kernels
$\kappa^{\eo}_1(x,y),$ then Lidskii's theorem, see~\cite{lax}, implies that
\begin{equation}
\Rind(\cP^{\eo}_+,\cR^{\eo}_+)=\int\limits_{bX}\kappa^{\eo}_2(x,x)dS(x)-
\int\limits_{bX}\kappa^{\eo}_1(x,x)dS(x).
\label{5.16.3}
\end{equation}
This formula, coupled with~\eqref{eqn2}, is very useful for showing the
constancy of the index under smooth isotopies of the structures involved in its
definition.

\section{Invertible doubles and the Calderon Projector}\label{sec6}

In order to better understand the functorial properties of sub-elliptic
boundary value problems and prove the Atiyah-Weinstein conjecture, it is
important to be able to deform the $\spnc$-structure and projectors without
changing the indices of the operators.  We now consider the dependence of the
various operators on the geometric structures. Of particular interest is the
dependence of the Calderon projector on $(J,g,\rho).$ To examine this we need
to consider the invertible double construction from~\cite{BBW} in
greater detail.  We also want to express the indices of
$(\eth^{\eo}_+,\cR^{\eo}_+)$ as the relative indices of the tame Fredholm pairs
$(\cP^{\eo}_+,\cR^{\eo}_+).$

We now recount the invertible double construction from~\cite{BBW}. We begin
with a compact manifold $X$ with boundary, with a metric $g,$ complex spinor
bundles $\Spn^{\eo}\to X,$ and $h$ a Hermitian metric on $\Spn^{\eo}.$ Let
$Y=bX$ and suppose that an identification of a neighborhood $U$ of $bX$ with
$Y\times [-1,0]_t$ is fixed. We assume that $dt$ is a outward pointing unit
co-vector.  With respect to this collar neighborhood, we say that $X$ has a
\emph{cylindrical end} if $\Spn^{\eo},$ $h$ and $g$ are independent of the ``normal
variable,'' $t$. In this case, the invertible double of $(X,g,h,\Spn^{\eo})$ is
defined to be $\tX=X\amalg_{Y}\overline{X},$ here $\overline{X}$ is $X$ with
the opposite orientation. We denote the components of $\tX\setminus
Y\times\{0\}$ by $X_{+} (t<0),$ $X_{-} (t>0).$ The smooth structure on $\tX$ is
obtained by gluing $Y\times [-1,0]\subset X_+$ to $Y\times [0,1]\subset X_-,$
along $Y\times \{0\}.$ As $\Spn^{\eo}, h$ and $g$ are independent of $t$ is is
clear that they extend smoothly to $\tX.$

Because the orientation of $X_-$ is reversed, to get a smooth bundle of complex
spinors we glue $\Spn^{\eo}\restrictedto_{Y}$ to
$\Spn^{\ooee}\restrictedto_{Y}$ using
\begin{equation}
\bc(-dt)\cdot\sigma_+\restrictedto_{Y\times 0^-}\sim
\sigma_-\restrictedto_{Y\times 0^+}.
\end{equation}
In~\cite{BBW} it is shown that this defines a smooth Clifford module over $\tX$
and hence a $\spnc$-structure. We let
$\eth^{\eo}_X$ denote the Dirac operator, and use the notation
\begin{equation}
\eth^{\eo}_{X\pm}\overset{d}{=}\eth^{\eo}_X\restrictedto_{\CI(X_{\pm};\Spn^{\eo})}.
\end{equation}
From the construction and results in~\cite{BBW}, the following identities are obvious
\begin{equation}
\Ker\eth^{\eo}_{X\pm}=\Ker\eth^{\ooee}_{X\mp}.
\label{5.16.6}
\end{equation}
In~\cite{BBW} it is shown that $\eth^{\eo}_X$ are invertible operators, we denote
the inverses by $Q^{\eo}.$

If $(X,J,g)$ is an almost complex manifold with boundary, then it can be
included into a larger manifold $\pX$ that has a cylindrical end. It is clear
that this can be done with smooth dependence on $(J,g).$ We fix an
identification of a neighborhood $U$ of $bX$ with $[-3,-2]\times Y.$ Using this
identification we smoothly glue $Y\times[-2,0]$ to $X.$ Denote this manifold by
$\pX.$ Using Lemma~\ref{lem8} below we easily show that the almost complex structure
can be extended to $\pX$ so that by the time we reach $t=-1$ it is independent
of $t.$ Hence we can also extend $\Spn^{\eo}$ to $\pX.$ Using the Seeley
extension theorem we can extend $(g,h)$ to $Y\times [-2,-1]$ in such a way that
the extended metric tensors depend continuously, in the $\CI$-topology, on
$(g,h)\restrictedto_{Y\times[-3,-2]},$ and $(g,h)$ also have a product
structure by the time we reach $Y\times \{-1\}.$ Everything can be further
extended to $Y\times [-1,0]$ so that it is independent of $t,$ and hence $\pX$,
with this hermitian spin-structure, has a cylindrical end.  
Compatible connections can be chosen on $\Spn^{\eo}$ so that both the metric
and spin geometries of $\thX=\pX\amalg_{Y\times \{0\}}\overline{\pX}$ depends
smoothly on the geometry of $(X,J,g,h).$ In particular the symbols of
$\eth^{\eo}_{\thX}$ depend smoothly on the symbols of $\eth^{\eo}_{X}.$

Fix a collar neighborhood, $U,$ of $bX$ so that $TX\restrictedto _U$ is
independent of $t.$ We can also normalize so that the 1-jet of $t$ along $t=0$
equals that of $\rho.$ We first homotope the Hermitian metric through a family
$\{g_s:\: s\in [0,\ha]\}$ so that the $g_0=g,$ and $g_{\ha}$ has a product
structure in $U.$ Moreover we can fix the metric on $bX$ throughout this
homotopy. For each $s$ there is a unique positive definite endomorphism $A_s$
of $TX$ so that, for all vector fields $V,W,$ we have
\begin{equation}
g(V,W)=g_s(A_sV,A_sW).
\end{equation}
If we set $J_s=A_sJ A_s^{-1},$ then this is a smooth family of almost complex
structures $\{J_s\}$ compatible with $g_s,$ and $J_s\restrictedto_{t=0}$
remaining fixed. Finally, with the metric in $U$ fixed to equal $g_{\ha},$ we can
deform $J_{\ha}$ through a family $\{J_s:\: s\in[\ha,1]\}$ so that:
\begin{enumerate}
\item $J_s=J_{\ha}$ outside of a small neighborhood of $bX.$
\item $J_1$ has a product structure within a smaller neighborhood of $bX.$
\item $J_s$ is compatible with $g_{\ha}$ for $s\in [\ha,1].$
\item $J_s\restrictedto_{t=0}=J\restrictedto_{t=0}.$
\end{enumerate}
That this is possible follows from the fact that the space of almost complex
structures compatible with $g_{\ha}$ can be represented as sections of a smooth
fiber bundle $\cJ$ with fiber equal to $SO(2n)/U(n).$ This representation is
obtained by using $J\restrictedto_{t=0},$ pulled back to the collar
neighborhood, to define a reference structure.  By compactness, there is an
$\epsilon<0$ so that, the section of $\cJ\restrictedto_{\epsilon\leq t\leq 0}$
defined by $J_{\ha}$ lies in a neighborhood of $\cJ$ retractable onto the
``zero section,'' defined by the reference structure. Hence we can perform the
homotopy described. We let $g_s=g_{\ha}$ for $s\in [\ha,1]$ and $\rho_s=t$ for
all $s.$ For later application we summarize the results of this discussion as a
lemma. We refer to this process as \emph{flattening the end}.
\begin{lemma}\label{lem8} Let $(X,J,g,\rho)$ be a normalized strictly pseudoconvex
  $\spnc$-manifold. There exists a smooth family
  $\{(X,J_s,g_s,\rho_s):\: s\in [0,1]\}$ of normalized strictly pseudoconvex
  $\spnc$-manifolds with
\begin{enumerate}
\item The structure at $s=0$ is equal to the given structure.
\item Throughout the homotopy, the data remains fixed along $bX.$
\item The space $(X,J_1,g_1,\rho_1)$ has a cylindrical end.
\end{enumerate}
\end{lemma}

We let $Q^{\eo}_{\thX}$ be the inverses of $\eth^{\eo}_{\thX}.$ These are
classical pseudodifferential operators of order $-1,$ whose symbols depend
smoothly on the symbols of $\eth^{\eo}_{\thX}$ and therefore, in turn on the
geometric data on $X.$ Throughout the discussion below we use the fact that the
operator norms of a pseudodifferential operator depend continuously on finite
semi-norms of the ``full'' symbol of the operator, see~\cite{Hormander3}.

We state a general result:
\begin{proposition}\label{prp10} Let $M$ be a compact manifold and $E,F$
  complex vector bundles over $M.$ Let $\{A_\tau\in\Psi^{1}(M;E,F):\: \tau\in
  \fT\}$ be a compact smooth family of invertible elliptic pseudodifferential
  operators. For any $s\in\bbR$ the family of inverses $A_{\tau}^{-1}$ is a
  norm continuous family of operators from $H^s(M;F)$ to $H^{s+1}(M;E).$
\end{proposition}
\begin{proof} As $\{A_\tau\}$ is a smooth family of pseudodifferential
  operators of order 1, for each $s\in\bbR,$ $\tau\mapsto \left[A_{\tau}:H^s\to
H^{s-1}\right]$ is norm continuous.  We can write
\begin{equation}
A_{\tau}=A_{\tau_0}(I-E_{\tau,\tau_0})
\label{1.10.6.1}
\end{equation}
where
\begin{equation}
E_{\tau,\tau_0}=A_{\tau_0}^{-1}(A_{\tau_0}-A_{\tau}).
\end{equation}
For any fixed $\tau_0,$ the operators $E_{\tau,\tau_0}$ are a smooth family of
pseudodifferential operators of order zero, with $E_{\tau_0,\tau_0}=0.$ Hence
for any fixed $s,$ there is a $\delta_s$ such that $|\tau-\tau_0|<\delta_s$
implies that operator norm of $E_{\tau,\tau_0}:H^s\to H^s$ is less than $\frac
12.$ Hence~\eqref{1.10.6.1} implies that $A^{-1}_\tau$ is given by the series:
\begin{equation}
A^{-1}_\tau=A^{-1}_{\tau_0}+\sum_{k=1}^{\infty}E_{\tau,\tau_0}^k,
A^{-1}_{\tau_0},
\label{5.17.5} 
\end{equation}
which converges in the norm topology of operators from $H^s$ to $H^{s+1}.$
This completing the proof of the proposition.
\end{proof}

Once the collar neighborhood is fixed, we can define the Calderon projectors
for the hypersurface $t=0:$  The Calderon projectors
$\cP^{\eo}_{\pm}$ are defined by
\begin{equation}
\cP^{\eo}_{\pm}g=\lim_{\epsilon\to 0^+} \gamma_{\mp\epsilon}
Q^{\eo}\gamma_0^{*}\bc(\mp dt)g.
\label{5.18.3}
\end{equation}
Here $\gamma_\epsilon$ is the operation of restriction to the submanifold
$t=\epsilon.$ In fact, we can define a pair of Calderon projectors for any
hypersurface $t=t_0$ lying in the collared part of the manifold. In~\cite{BBW}
it is shown that,
\begin{equation}
\cP^{\eo}_{+}+\cP^{\eo}_{-}=\Id.
\label{5.17.9}
\end{equation}
We can extend $t$ smoothly to all of $\thX$ so that it is negative on the
original manifold $X,$ and positive on the interior of $\thX\setminus\pX.$ The
following result is very useful in our analysis.
\begin{proposition}\label{prp11} Let $X$ be a $\spnc$-manifold with boundary and let $\thX$ be
  an invertible double for $X.$ If $\cP^{\eo}_{\pm}$ are the Calderon
  projectors defined by the above prescription, then the adjoints satisfy:
\begin{equation}
\cP^{\eo*}_{\pm}=\bc(\pm dt)\cP^{\ooee}_{\mp}\bc(\pm dt)^{-1}.
\label{5.17.10}
\end{equation}
\end{proposition}
\begin{proof} We give the proof for the $+$ even case; the odd and $-$ cases are
  identical. Let $f$ and $g$ be smooth sections of
  $\Spn^{\even}\restrictedto_{bX},$ and let $G$ denote the extension of $g$ to
  the collar neighborhood of $bX$ that is constant in $t.$ As
  $[Q^{\eo}]^{*}=Q^{\ooee},$ the definition implies that
\begin{equation}
\begin{split}
\langle f,\cP^{\even*}_{+}g\rangle=
\langle \cP^{\even}_{+}f,g\rangle&=\lim_{\epsilon\to 0^-}
\langle\gamma_{\epsilon}
Q^{\even}\gamma_0^{*}\bc(-dt)f,G\rangle\\
&=\lim_{\epsilon\to 0^-}\langle 
\bc(-dt)f,\gamma_0 Q^{\odd}\gamma_{\epsilon}^* G\rangle
\end{split}
\label{5.18.2}
\end{equation}
The proof of the proposition follows from the observations that $\epsilon<0$
in~\eqref{5.18.2}, $\bc(-dt)^*=\bc(dt),$ and
\begin{equation}
\lim_{\epsilon\to 0^-}\gamma_0 Q^{\odd}\gamma_{\epsilon}^* G=
\cP^{\odd}_-\bc(dt)^{-1} g.
\label{5.17.11}
\end{equation}
This follows because $u_{\epsilon}=Q^{\odd}\gamma_{\epsilon}^* G$ solves
$\eth^{\odd}u_{\epsilon}=0$ in the subset
 $$\thX_{\epsilon}=\{x:\:t(x)>\epsilon\}.$$ 
It is easy to see that, as $\epsilon\uparrow 0,$ this is a
uniformly bounded family of solutions, which converges uniformly on $\thX_-$ to
$u_0=Q^{\odd}\gamma_{0}^* G.$ Clearly the restrictions to $t=0$ converge to
$\cP^{\odd}_-\bc(dt)^{-1} g.$ 
\end{proof}
This proposition has an interesting and useful corollary
\begin{corollary} Suppose that $X$ is $\spnc$-manifold, with a cylindrical end and 
  $\cP^{\eo}_{\pm}$ are defined using the inverse of $\eth^{\eo}_{\tX}$ on the
  invertible double $\tX.$ Then these are self adjoint projection operators.
\end{corollary}
\begin{proof} We do the even-$+$ case, the others are
  identical. Proposition~\ref{prp11} shows that
  $\cP^{\even*}_{+}=\bc(dt)\cP^{\odd}_-\bc(dt)^{-1}.$ On the other hand, because
  $\tX$ is obtained by doubling across $bX,$ we have~\eqref{5.16.6}, implying
  that
\begin{equation}
\range\cP^{\even}_+=\range\cP^{\even *}_+.
\end{equation}
Generally we have that
\begin{equation}
[\range\cP^{\even}_+]^{\bot}=\Ker\cP^{\even *}_+=\range(\Id-\cP^{\even *}_+)
=\range(\Id-\cP^{\even }_+),
\end{equation}
and therefore $\langle\cP^{\even}f,(\Id-\cP^{\even})g\rangle=0,$ for all pairs,
$f,g.$ 
These relations imply that
\begin{equation}
\cP^{\even *}_+\cP^{\even}_+=\cP^{\even}\text{ and }
\cP^{\even *}_+(\Id-\cP^{\even}_+)=0,
\end{equation}
from which the conclusion is immediate.
\end{proof}

The symbols of $\cP^{\eo}_{\pm}$ are smooth functions of the
symbols of $Q^{\eo}_{\thX}.$ Using the norm continuity of $Q^{\eo}_{\thX}$ we
conclude that the Calderon projectors also depend continuously, in the
uniform norm topology, on the geometric data on $X.$ 
\begin{proposition}\label{prp12} Suppose that $\{(X,J_{\tau},g_{\tau},\rho_{\tau}):\:
  \tau\in \fT\}$ is a compact smooth family of normalized strictly pseudoconvex
  $\spnc$-manifolds. The Calderon projectors, $\cP^{\eo}_{\pm\tau},$
  defined by the invertible double construction are smooth families of
  pseudodifferential operators of order zero, and
\begin{equation}
\tau\mapsto \left[\cP^{\eo}_{\pm\tau}:L^2(bX;\Spn^{\eo})\to L^2(bX;\Spn^{\eo})\right]
\end{equation}
are continuous in the uniform norm topology.
\end{proposition}
\begin{proof} First we show that $\cP^{\eo}_{\pm\tau}$ is a norm continuous
  family of operators on $L^2.$ Let $Q^{\eo}_{\tau}$ denote the inverse of
  $\eth^{\eo}_{\tau},$ the $\spnc$-Dirac operator on the invertible double
  defined by the data $(X,J_{\tau},g_{\tau},\rho_{\tau}),$ and
  $\tQ^{\eo}_{\tau}$ a parametrix with
\begin{equation}
\eth^{\eo}_{\tau}\tQ^{\eo}_{\tau}=\Id-K_{1\tau},\quad
\tQ^{\eo}_{\tau}\eth^{\eo}_{\tau}=\Id-K_{2\tau}.
\end{equation}
Proposition~\ref{prp10} shows that $Q^{\eo}_{\tau}:L^2(X;\Spn^{\ooee})\to
  H^1(X;\Spn^{\eo})$ are norm continuous families. 
The inverse and the parametrix are related by
\begin{equation}
Q^{\eo}_{\tau}=\tQ^{\eo}_{\tau}+Q^{\eo}_{\tau}K_{1\tau}\text{ and }
Q^{\eo}_{\tau}=\tQ^{\eo}_{\tau}+K_{2\tau}Q^{\eo}_{\tau}.
\label{5.17.6}
\end{equation}
Recall that the restriction maps from $H_{\hn}(X_{\pm};\Spn^{\eo})\to
L^2(bX_{\pm};\Spn^{\eo})$ are continuous. The second statement of the
proposition is an easy consequence of this fact, the relation~\eqref{5.17.6},
and the observation that
\begin{equation}
\tQ^{\eo}_{\tau}:L^2(bX;\Spn^{\eo})\to H_{\hn}(X_{\pm};\Spn^{\eo})\text{ and }
K_{1\tau}:L^2(bX;\Spn^{\eo})\to H^1(X;\Spn^{\eo})
\label{5.17.8}
\end{equation}
are norm continuous families. The proof that the maps in~\eqref{5.17.8} define
norm continuous families is a simple adaptation of the argument showing that
the norm of a pseudodifferential operator is bounded by a finite semi-norm of
its full symbol, which we leave to the interested reader.

To see that $\cP^{\eo}_{\pm\tau}$ is a smooth family of pseudodifferential
operators we use the relation,~\eqref{5.17.6} to conclude that
\begin{equation}
\cP^{\eo}_{\pm\tau}=\tcP^{\eo}_{\pm\tau}(\Id+k'_{1\tau})
+k_{2\tau}\cP^{\eo}_{\pm\tau}k_{1\tau},
\end{equation}
where $\tcP^{\eo}_{\pm\tau}$ is the smooth family of pseudodifferential
operators defined by using the parametrices, $\tQ^{\eo}_{\tau},$ in the
definition of the Calderon projector,~\eqref{5.18.3}, and $k'_{1\tau},k_{1\tau},
k_{2\tau}$ are smooth families of smoothing operators. This relation, combined
with the $L^2$-norm continuity of $\tau\mapsto\cP^{\eo}_{\pm\tau},$ show that
$\cP^{\eo}_{\pm\tau}$ is a smooth family of pseudodifferential operators.
\end{proof}

We can suppose that the contact structure induced on the boundary of the family
$(X,J_\tau,g_\tau,\rho_\tau)$ is fixed and we let $\cR^{\eo}_{+\tau}$ denote a
smooth family of modified $\dbar$-Neumann conditions. The non-trivial part of
such a family is a smooth family of generalized Szeg\H o projectors
$\tau\mapsto\cS_{\tau}.$ In~\cite{EpsteinMelrose} it is shown that such a
family is norm continuous as a family of maps
$\tau\mapsto[\cS_{\tau}:L^2(bX)\to L^2(bX)].$
\begin{theorem}\label{thm6} If $\{(X,J_{\tau},g_{\tau},\rho_{\tau}):\: \tau\in
  \fT\}$ is a compact smooth family of normalized strictly pseudoconvex
  $\spnc$-manifolds and $\cR^{\eo}_{+\tau}$ is a smooth family of modified
  $\dbar$-Neumann conditions, then $\Rind(\cP^{\eo}_{+\tau},\cR^{\eo}_{+\tau})$
  is constant.
\end{theorem}
\begin{proof} Proposition~\ref{prp12} shows that the operators
  $\cT^{\eo}_{+\tau}$ are a smooth family of extended Heisenberg operators and
  therefore so are the parametrices $\cU^{\eo}_{+\tau}.$ Hence the residual
  terms
\begin{equation}
K^{\eo}_{1\tau}=\Id-\cT^{\eo}_{+\tau}\cU^{\eo}_{+\tau},\quad
K^{\eo}_{2\tau}=\Id-\cU^{\eo}_{+\tau}\cT^{\eo}_{+\tau}
\end{equation}
are smooth families of smoothing operators. As $\tau\mapsto\cP^{\eo}_{+\tau}$
and $\tau\mapsto\cR^{\eo}_{+\tau}$ are norm continuous as maps from $L^2$ to
itself, the operators $\cP^{\eo}_{+\tau}K^{\eo}_{2\tau}\cP^{\eo}_{+\tau}$ and
$\cR^{\eo}_{+\tau}K^{\eo}_{1\tau}\cR^{\eo}_{+\tau}$ are continuous in the trace
norm. It follows from~\eqref{5.16.4} that
$\Rind(\cP^{\eo}_{+\tau},\cR^{\eo}_{+\tau})$ depends continuously on $\tau.$ As
it is integer valued, it is constant.
\end{proof}

\section{The relative index formula}
In this section we prove the formula in~\eqref{eqn2}, expressing
$\Ind(\eth^{\eo}_{+},\cR^{\eo}_+)$ as the relative index of the tame
Fredholm pair $(\cP^{\eo}_{+},\cR^{\eo}_+),$ and derive several
consequences of this formula. 
\begin{theorem}\label{thm7} Let $(X,J,g,\rho)$ be a normalized strictly pseudoconvex
  $\spnc$-manifold, and $\cR^{\eo}_+$  projections defining a modified
  $\dbar$-Neumann problems, then
\begin{equation}
\Ind(\eth^{\eo}_+,\cR^{\eo}_+)=\Rind(\cP^{\eo}_+,\cR^{\eo}_+).
\label{5.17.15}
\end{equation}
\end{theorem}
\begin{proof} We give the proof for the even case, the odd case is
  identical. The kernel of $(\eth^{\even}_+,\cR^{\even}_+)$ consists of smooth
  forms $\sigma$ that satisfy
\begin{equation}
\eth^{\even}_+\sigma=0\text{ and }\cR^{\even}_+[\sigma]_{bX}=0.
\end{equation}
The first condition implies that $\cP^{\even}_+[\sigma]_{bX}=[\sigma]_{bX},$
and therefore $\sigma\in\Ker\cR^{\even}_+\cP^{\even}_+.$ Conversely, by the
unique continuation theorem for $\Ker\eth^{\eo}_{\pm},$ any form in the range
of $\cP^{\even}_+$ that lies in the kernel of $\cR^{\even}_+$ defines a unique
form in $\Ker(\eth^{\even}_+,\cR^{\even}_+).$ Thus
$$\Ker(\eth^{\even}_+,\cR^{\even}_+)\simeq
\Ker\cR^{\even}_+\cP^{\even}_+\restrictedto_{\range\cP^{\even}_+}.$$

The cokernel of $\cR^{\even}_+\cP^{\even}_+$ is isomorphic to the null space
of
\begin{equation}
\cP^{\even*}_+:\range\cR^{\even}_+\to\range\cP^{\even*}_+.
\end{equation}
Proposition~\ref{prp11} and equation~\eqref{5.17.9} show that 
\begin{equation}
\cP^{\even*}_+=\bc(dt)(\Id-\cP^{\odd}_+)\bc(dt)^{-1}
\end{equation}
This identity, along with~\eqref{5.17.17} show that the cokernel of
$\cR^{\even}_+\cP^{\even}_+$ is isomorphic to the null space of $\cR^{\odd}_+$
acting on $\range\cP^{\odd}_+,$ which, by the first part of this argument, is
isomorphic to $\Ker(\eth^{\odd}_+,\cR^{\odd}_+).$ Applying Theorem~\ref{thm3},
we complete the proof of the theorem.
\end{proof}

As a corollary we have
\begin{corollary}\label{cor5} Let $(X,J,g,\rho)$ be a normalized strictly
  pseudoconvex $\spnc$-manifold, and $\cR^{\eo}_+$ projection operators
  defining modified $\dbar$-Neumann conditions, then
\begin{equation}
\Rind(\cP^{\eo}_+,\cR^{\eo}_+)=-\Rind(\Id-\cP^{\eo}_+,\Id-\cR^{\eo}_+)
\label{5.17.18}
\end{equation}
\end{corollary}
\begin{proof} We give the proof for the even case, the odd case is
  identical.  First suppose that $(X,J,g,\rho)$ has a cylindrical end. In
this case $\cP^{\even}_+=\cP^{\even*}_+$ and therefore~\eqref{5.17.9} and
Proposition~\ref{prp11} imply that
\begin{equation}
\Id-\cP^{\even}_+=\bc(dt)\cP^{\odd}_{+}\bc(dt)^{-1}
\end{equation}
As $\Id-\cR^{\even}_+=\bc(dt)\cR^{\odd}_{+}\bc(dt)^{-1},$  the
relation~\eqref{5.17.15} implies that
\begin{equation}
\Rind(\Id-\cP^{\even}_+,\Id-\cR^{\even}_+)=\Rind(\cP^{\odd}_+,\cR^{\odd}_+)=
\Ind(\eth^{\odd}_+,\cR^{\odd}_+)=-\Ind(\eth^{\even}_+,\cR^{\even}_+).
\end{equation}
The corollary, in this case, follows with one further application
of~\eqref{5.17.15}.

We treat the general case by deformation. Lemma~\ref{lem8} gives a smooth
family $\{(X,J_s,g_s,\rho_s):\: s\in [0,1]\}$ such that $(X,J_1,g_1,\rho_1)$
has a cylindrical end and the data along $bX$ is fixed.  As the metric, almost
complex structure and defining function remain constant along $bX$, the family
$\{(X,J_s,g_s,\rho_s):\: s\in [0,1]\}$ satisfies the hypotheses of
Theorem~\ref{thm6} and therefore $\Rind(\Id-\cP^{\eo}_{+s},\Id-\cR^{\eo}_{+s})$
and $\Rind(\cP^{\eo}_{+s},\cR^{\eo}_{+s})$ are independent of $s.$ The argument
above shows that
\begin{equation}
\Rind(\Id-\cP^{\eo}_{+1},\Id-\cR^{\eo}_{+1})=
-\Rind(\cP^{\eo}_{+1},\cR^{\eo}_{+1}),
\end{equation}
completing the proof of the theorem.
\end{proof}

As a corollary of the corollary we have the following result.
\begin{corollary} The operators  $\cT^{\eo}_+$   are tame Fredholm operators of
  index   zero.
\end{corollary}
\begin{proof} The first statement follows from the ellipticity of $\cT^{\eo}_+$
  in the extended Heisenberg calculus.  The indices of $\cT^{\eo}_+$ are computed
  using the trace formula:
\begin{equation}
\Ind(\cT^{\eo}_+)=\Tr K^{\eo}_2-\Tr K^{\eo}_1.
\label{5.20.1}
\end{equation}
Using formula~\eqref{5.16.4} and its analogue for
$\Rind(\Id-\cP^{\eo}_+,\Id-\cR^{\eo}_+)$, we see that
\begin{multline}
\Rind(\cP^{\eo}_+,\cR^{\eo}_+)+\Rind(\Id-\cP^{\eo}_+,\Id-\cR^{\eo}_+)=\\
\Tr \cP^{\eo}_+ K_2^{\eo} \cP^{\eo}_+-\Tr \cR^{\eo}_+ K_1^{\eo} \cR^{\eo}_+
+\\
\Tr(\Id- \cP^{\eo}_+) K_2^{\eo}(\Id- \cP^{\eo}_+)-\Tr (\Id-\cR^{\eo}_+)
K_1^{\eo}(\Id- \cR^{\eo}_+)
\label{5.20.111}
\end{multline}
It is a simple computation to show that the right hand side
of~\eqref{5.20.111} equals
\begin{multline}
\Tr K^{\eo}_2-\Tr K^{\eo}_1+\\\Tr [\cP^{\eo}_+ K_2^{\eo},\cP^{\eo}_+]+
\Tr [\cP^{\eo}_+, K_2^{\eo}\cP^{\eo}_+]-
\Tr[\cR^{\eo}_+ K_1^{\eo}, \cR^{\eo}_+]-\Tr[\cR^{\eo}_+ ,K_1^{\eo}
  \cR^{\eo}_+].
\label{5.20.1111}
\end{multline}
The commutators on the right hand side of~\eqref{5.20.1111} are of the form
$\Tr[K,A]$ where $K$ is a smoothing operator and $A$ is bounded. As such terms
vanish, this corollary follows from~\eqref{5.20.1}--\eqref{5.20.1111} and Corollary~\ref{cor5}.
\end{proof}

\section{The Agranovich-Dynin formula}
In~\cite{Epstein4} we proved a generalization of the Agranovich-Dynin formula
for subelliptic boundary conditions, assuming that the $\spnc$-structure arises
from an integrable almost complex structure. In this section we show that the
integrability is not necessary.
\begin{theorem}\label{thm8} Let $(X,J,h,\rho)$ be a normalized strictly pseudoconvex
  $\spnc$-manifold. Let $\cS_1$ and $\cS_2$ be two generalized Szeg\H o
  projections and $\cR^{\even}_{+1}, \cR^{\even}_{+2}$ the modified $\dbar$-Neumann
  conditions they define, then
\begin{equation}
\Rind(\cS_1,\cS_2)=\Ind(\eth^{\even},\cR^{\even}_{+2})-\Ind(\eth^{\even},\cR^{\even}_{+1}).
\label{5.18.4}
\end{equation}
\end{theorem}
\begin{proof} We apply Theorem~\ref{thm7} to conclude that
\begin{equation}
\Ind(\eth^{\even}_+,\cR^{\even}_{+2})-\Ind(\eth^{\even}_+,\cR^{\even}_{+1})=
\Rind(\cP^{\even}_+,\cR^{\even}_{+2})-\Rind(\cP^{\even}_+,\cR^{\even}_{+1}).
\label{5.18.5}
\end{equation}
By relabeling we may assume that $\Rind(\cS_1,\cS_2)\geq 0.$
In~\cite{EpsteinMelrose} it is shown that $\cS_2$ can be deformed, through a
smooth family of generalized Szeg\H o projections to a projection $\cS_3$ which
is a sub-projection of $\cS_1.$ That is $\cS_1=\cS_3+P,$ where $P$ is a finite
rank orthogonal projection, with a smooth kernel and
\begin{equation}
\cS_3 P= P\cS_3=0.
\end{equation}
If $\{\cR^{\even}_{+s}:\: s\in [2,3]\}$ is the associated family of modified
$\dbar$-Neumann conditions, then the proof of Theorem~\ref{thm6} applies
equally well to show that $\Rind(\cP^{\even}_+,\cR^{\even}_{+s})$ is
constant. So we are reduced to showing that
\begin{equation}
\Rind(\cS_1,\cS_2)=\Rind(\cS_1,\cS_3)=
\Rind(\cP^{\even}_{+},\cR^{\even}_{+3})-\Rind(\cP^{\even}_{+},\cR^{\even}_{+1}),
\label{5.18.6}
\end{equation}
with $\cS_3$ a finite corank sub-projection of $\cS_1.$ For convenience we apply
Lemma~\ref{lem8} to deform to a structure with a cylindrical end, so that we can 
assume that $\cP^{\even}_+$ is self adjoint. Theorem~\ref{thm11} in the
appendix applies to show that 
$\Rind(\cP^{\even}_{+},\cR^{\even}_{+3})-\Rind(\cP^{\even}_{+},\cR^{\even}_{+1})$ is
the index of the operator 
\begin{equation}
\cR^{\even}_{+3}\cP^{\even}_{+}\cR^{\even}_{+1}:L^2\cap\range\cR^{\even}_{+1}
\longrightarrow\cK\cap\range\cR^{\even}_{+3},
\label{5.18.9}
\end{equation}
where $\cK$ is an appropriately defined Hilbert space. Let $\cU^{\even}_{+j}$
$j=1,3$ denote parametrices for $\cT^{\even}_{+j}.$ 
The space $\cK$ is defined as the closure of $\CI$ with respect to the inner
product defined by
\begin{equation}
\|u\|^2_{\cK}=\|u\|^2_{L^2}+\|\cU^{\even*}_{+1}\cU^{\even}_{+3}u\|^2_{L^2}.
\label{5.18.15}
\end{equation}

We can rewrite
\begin{equation}
\cR^{\even}_{+3}\cP^{\even}_+\cR^{\even}_{+1}=
\cR^{\even}_{+3}\cT^{\even}_{+3}\cT^{\even *}_{+1}\cR^{\even}_{+1}
\label{5.18.10}
\end{equation}
We recall that the difference $\cR^{\even}_{+1}-\cR^{\even}_{+3}$ is a smoothing
operator and therefore so is $\cT^{\even}_{+1}-\cT^{\even}_{+3}.$ Hence
the operator on the right hand side of~\eqref{5.18.10} has the same index as
\begin{equation}
\cR^{\even}_{+3}\cT^{\even}_{+3}\cT^{\even *}_{+3}\cR^{\even}_{+1}
\label{5.18.11}
\end{equation}

It is a simple matter to show that
$[\cR^{\even}_{+3},\cT^{\even}_{+3}\cT^{\even *}_{+3}]=0$ and therefore, as
$[\cR^{\even}_{+3}]^2=\cR^{\even}_{+3},$ we are reduced to computing the index of
\begin{equation}
\cR^{\even}_{+3}\cT^{\even}_{+3}\cT^{\even *}_{+3}\cR^{\even}_{+3}\cR^{\even}_{+1}.
\label{5.18.12}
\end{equation}
We think of this as a composition
\begin{equation}
\cR^{\even}_{+3}\cR^{\even}_{+1}:L^2\cap\range\cR^{\even}_{+1}\to
L^2\cap\range\cR^{\even}_{+3}
\label{5.18.13}
\end{equation}
with
\begin{equation}
\cW=\cR^{\even}_{+3}\cT^{\even}_{+3}\cT^{\even
  *}_{+3}\cR^{\even}_{+3}:L^2\cap\range\cR^{\even}_{+3}\to
\range\cR^{\even}_{+3}\cap\cK.
\end{equation}
The index of the operator in~\eqref{5.18.13} is
$\Rind(\cS_1,\cS_3)=\Rind(\cS_1,\cS_2);$ hence we are left to show that $\cW$
has index zero. 

As $\cT^{\even}_{+1}-\cT^{\even}_{+3}$ is a smoothing operator,
the space defined by the inner product in~\eqref{5.18.15} is unchanged if
$\cU^{\even *}_{+1}$ is replaced by $\cU^{\even *}_{+3};$ thus $\cW$ is a
bounded operator. The operator
\begin{equation}
\cV=\cR^{\even}_{+3}\cU^{\even
  *}_{+3}\cU^{\even}_{+3}\cR^{\even}_{+3}:\range\cR^{\even}_{+3}\cap\cK\to
L^2\cap\range\cR^{\even}_{+3}
\end{equation}
is a parametrix for $\cR^{\even}_{+3}\cT^{\even}_{+3}\cT^{\even
  *}_{+3}\cR^{\even}_{+3}.$ If
\begin{equation}
\cV\cW=
\cR^{\even}_{+3}-\cR^{\even}_{+3}K\cR^{\even}_{+3},
\end{equation}
where $K$ is smoothing, then
\begin{equation}
\cW\cV=
\cR^{\even}_{+3}-\cR^{\even}_{+3}K^*\cR^{\even}_{+3}.
\end{equation}
Hence the $\Ind(\cW)$ is given by
\begin{equation}
\Ind(\cW)=\Tr(\cR^{\even}_{+3}K^*\cR^{\even}_{+3})-
\Tr(\cR^{\even}_{+3}K\cR^{\even}_{+3}).
\label{5.18.16}
\end{equation}
But this means it must be zero, because the right hand side of~\eqref{5.18.16}
is a purely imaginary number. This completes the proof of Theorem~\ref{thm8}.
\end{proof}
\begin{remark} The proof of the Atiyah-Weinstein conjecture is a small
  modification of the proof of the Agranovich-Dynin formula.
\end{remark}

\section{The Atiyah-Weinstein conjecture}
In~\cite{Weinstein} Weinstein considers the following situation: let $X_0, X_1$
be strictly pseudoconvex Stein manifolds with boundary. The CR-structures on
$bX_0,$ and $bX_1$ define Szeg\H o projectors $\cS_0, \cS_1$ as  projectors
onto the nullspaces of $\dbarb$-operators acting on functions. Suppose that
there is a contact diffeomorphism $\phi: bX_1\to bX_0.$  Weinstein describes a
construction, using stable almost complex structures, for gluing $X_0$ to $X_1$
via $\phi,$ to obtain a compact manifold $X$ with a well defined
$\spnc$-structure.  Weinstein conjectures that
\begin{equation}
\Rind(\cS_0,[\phi^{-1}]^*\cS_1\phi^{*})=\Ind(\eth^{\even}_{X}).
\end{equation}
He also gives a conjecture for a formula when $X_0,$ or $X_1$ is a Stein space
and not a Stein manifold. As described in the introduction, these conjectures
evolved from conjectures, made jointly with Michael Atiyah in the 1970s, for the
indices of elliptic Fourier integral operators defined by contact
transformations between co-sphere bundles of compact manifolds.

In this section we prove a more general formula, covering all these cases. As
noted earlier, we do not use the stable almost complex structure construction
to build a $\spnc$-structure on $X,$ but rather a simple extension of the
invertible double construction. It seems clear that these two constructions
lead to the same $\spnc$-structure on the glued space.

Our set-up is the following: we have two normalized strictly pseudoconvex
$\spnc$-manifolds, $(X_k,J_k,g_k,\rho_k),$ $k=0,1,$ and a co-orientation
preserving contact diffeomorphism of their boundaries, $\phi:bX_1\to bX_0.$
On each of these boundaries we choose a generalized Szeg\H o projector,
$\cS_0,\cS_1.$ Using the contact diffeomorphism we obtain a second
Szeg\H o projector on $bX_0,$ where
\begin{equation}
\cS'_1=\phi^{-1*}\cS_1\phi^*.
\label{5.19.9}
\end{equation}

We can extend the contact diffeomorphism to a diffeomorphism of collar
neighborhoods of the boundaries. Let us heretofore suppose that such an
identification is fixed. The two compatible almost complex structures can now
be regarded as being defined on a neighborhood of one and the same contact
manifold. We show below that there is a homotopy, $\{J_s:\: s\in [0,1]\},$
through compatible almost complex structures joining $J_0$ to $J_1.$ Using this
and our collar neighborhood construction we can add collars to $X_0,$ $X_1$
obtaining $\spnc$-manifolds $\hX_0,\hX_1$ with cylindrical ends and
identical structures on a collar neighborhood of their boundaries. We describe
this construction more precisely below.

Using the obvious extension of the invertible double construction, we can now
build a compact space
\begin{equation}
\tX_{01}=\hX_0\amalg_{b\hX_j}\overline{\hX_1},
\label{5.19.10}
\end{equation}
with a well defined isotopy class of $\spnc$-structures and Dirac operator
$\eth^{\even}_{\tX_{01}}.$ The operators $\eth^{\eo}_{\tX_{01}}$ need not be
invertible, nor have index zero. In the sequel we refer to this as the
\emph{extended}  double construction.
\begin{theorem}\label{thm9} Let $(X_k,j_k,g_k,\rho_k), k=0,1$ be normalized
  strictly pseudoconvex $\spnc$-manifolds, and suppose that $\phi:bX_1\to bX_1$
  is a co-orientation preserving contact diffeomorphism. Suppose that $\cS_0,
  \cS_1$ are generalized Szeg\H o projections on $bX_0, bX_1,$ respectively,
  which define modified $\dbar$-Neumann conditions $\cR^{\even}_{+j}$ on $X_j.$
  With $\cS_1'$ as defined in~\eqref{5.19.9} and $\tX_{01}$ as defined
  in~\eqref{5.19.10} we have
\begin{equation}
\Rind(\cS_0,\cS'_1)=\Ind(\eth^{\even}_{\tX_{01}})-\Ind(\eth^{\even}_{X_0},\cR^{\even}_{+0})+
\Ind(\eth^{\even}_{X_1},\cR^{\even}_{+1}).
\label{5.18.20}
\end{equation}
\end{theorem}
\noindent
In the sequel, we refer to the indices of the boundary value problems
in~\eqref{5.18.20} as the \emph{boundary terms}.
\begin{remark}
This includes Weinstein's conjectured formula: we take the CR-structures $J_0,
J_1$ to be integrable, and $X_0, X_1$ to be compact complex manifolds with the
given boundary. The classical Szeg\H o projectors are used for $\cS_0,\cS_1,$
respectively. In this case, it is shown in Section 7 of~\cite{Epstein4} that
the indices of the boundary value problems in~\eqref{5.18.20} are renormalized
holomorphic Euler characteristics:
\begin{equation}
\Ind(\eth^{\even}_{X_j},\cR^{\even}_{j})=\chi_{\cO}'(X_j)=
\sum_{q=1}^{n}\dim H^{0,q}(X_j)(-1)^q.
\end{equation}
Hence~\eqref{5.18.20} gives
\begin{equation}
\Rind(\cS_0,\cS_1)=\Ind(\eth^{\even}_{\tX_{01}})-\chi'_{\cO}(X_0)+
\chi'_{\cO}(X_1).
\label{5.18.21}
\end{equation}
If $X_0$ and $X_1$ are Stein manifolds then $\chi'_{\cO}(X_j)=0, j=0,1.$
\end{remark}

Applying the Atiyah-Singer theorem for a Dirac operator we obtain a
(partially) cohomological formula for the index.
\begin{corollary} With the hypotheses of Theorem~\ref{thm9}, we have
 \begin{equation}
\Rind(\cS_0,\cS'_1)=\langle e^{\ha c_1}\hat{\bf A}(\tX_{01}),\tX_{01}\rangle
-\Ind(\eth^{\even}_{X_0},\cR^{\even}_{+0})+
\Ind(\eth^{\even}_{X_1},\cR^{\even}_{+1}),
\label{5.18.201}
\end{equation}
with $c_1=c_1(\Spn^{\even})$ the canonical class of the $\spnc$-structure on
$\tX_{01}.$
\end{corollary}
\begin{remark} The terms $\Ind(\eth^{\even}_{X_0},\cR^{\even}_{+0})$ and
$\Ind(\eth^{\even}_{X_1},\cR^{\even}_{+1})$ are essential parts of this formula
  as they capture the non-symbolic nature of the relative index. In the case
  that $X_0, X_1$ are the co-ball bundles of compact manifolds, an equivalent
  formula appears in~\cite{LeichtNestTsy}. A related formula for a contact self
  map is given in~\cite{EpsteinMelrose}.
\end{remark}

If $\phi'$ is a different choice of contact diffeomorphism, then
$\psi=\phi'\circ\phi^{-1}$ is a contact automorphism of $bX_0.$ The projector
\begin{equation}
\cS''_1=[\phi']^{-1*}\cS_1\phi^{\prime *}=\psi^{-1*}\cS'_1\psi^*.
\end{equation}
The cocycle formula, proved in~\cite{EpsteinMelrose} shows that
\begin{equation}
\Rind(\cS_0,\cS''_1)=\Rind(\cS_0,\cS'_1)+\cdeg(\psi),
\end{equation}
where $\cdeg(\psi)$ is the contact degree. This is shown to be a topological
invariant of the isotopy class of $\psi$ in the contact mapping class group. A
formula for $\cdeg(\psi)$ as the index of a Dirac operator on the mapping torus
of $\psi,$ $\eth_{Z_\psi}$ is also provided, i.e.
\begin{equation}
\cdeg(\psi)=\Ind(\eth_{Z_\psi}).
\end{equation} 
Using the cohomological expression for $\Ind(\eth_{Z_\psi}),$ one can easily
show that the contact degree alway vanishes if $\dim Y=3.$ Hence, if
$\dim_{\bbR} X_j=4,$ then the $\Rind(\cS_0,\cS'_1)$ does not depend on the
choice of contact diffeomorphism.

\begin{proof}[Proof of Theorem~\ref{thm9}]
We now turn to the details of the proof of Theorem~\ref{thm9}. We suppose that the
extension of $\phi$ has been applied to identify a neighborhood of $bX_0$ with
a neighborhood of $bX_1.$ We then apply Lemma~\ref{lem8} to reduce to the
situation that $(X_k,J_k,g_k,\rho_k),$ $k=0,1$ have cylindrical ends.  This
deformation does not change the index of the operators
$(\eth^{\even}_{X_k},\cR^{\even}_{k}).$ 

We use $j_0,j_1$ to denote the compatible almost CR-structures defined on $H$
by $J_0,J_1,$ respectively. The set of almost CR-structures compatible with a
given contact form is contractible. Let $\{j_s:\: s\in [0,1]\}$ denote a
deformation of the compatible almost CR-structure $j_0$ to the compatible
almost CR-structure $j_1$ through compatible almost CR-structures. The contact
form, $\theta$ is fixed throughout this deformation and hence, so is the Reeb
vector field $T.$ Let $\pi_H$ denote the projection of $TY\times [0,1]$ to $H$
(pulled back to $Y\times [0,1]$) along $\Span\{T,\pa_t\}.$ We extend the almost
CR-structure to an almost complex structure, $J$ on $T Y\times[0,1]$ by setting
$J T=\pa_t.$ With this definition, the function $t$ satisfies
\begin{equation}
\theta\restrictedto_{Y\times\{t\}}=\frac{i}{2}\dbar
t\restrictedto_{Y\times\{t\}}
\label{5.19.1}
\end{equation}
We extend $\theta$ to $\Theta$ defined on $TY\times[0,1]$ by setting
\begin{equation}
\Theta(\pa_t)=0.
\end{equation}
The metric on the collar is given by
\begin{equation}
ds^2_{Y\times\{s\}}=dt^2+\Theta\cdot\Theta +d\theta(\pi_H\cdot, j_s\pi_H\cdot).
\end{equation}

Observe that  we can reparametrize the family $\{j_s\}$ so that both
ends of $Y\times [0,1]$ are cylindrical. To do this we choose a smooth non-decreasing function
$\varphi:[0,1]\to [0,1]$ with $\varphi(s)=0$ for $s\in [0,\oqu]$ and
$\varphi(s)=1$ for $s\in [\frac 34,1].$ If we replace $\{j_s\}$ with
$\{j_{\varphi(s)}\},$ then the forgoing construction defines an almost complex
structure on $Y\times [0,1]$ with both ends cylindrical and agreeing with the given
structures. We summarize this construction in a lemma.
\begin{lemma}\label{lem9} Suppose that $(Y,H)$ is a contact manifold with contact form $\theta$
  and two compatible almost CR-structures $j_0,j_1.$ Then there is an almost
  complex structure on $Y\times [0,1]$ with both ends cylindrical. The structure
  induced on $Y\times\{1\}$ is strictly pseudoconvex and agrees with $j_1,$
  while that on $Y\times\{0\}$ is strictly pseudoconcave and agrees with $j_0,$
  with its co-orientation reversed. For all members of the family the
  relation~\eqref{5.19.1} holds.
\end{lemma}

We use Lemma~\ref{lem9} to define an almost complex structure on $Y\times
[0,1].$ For each $0\leq\tau\leq 1$ we set $\hX^\tau_0=X_0\amalg_{Y\times\{\tau\}}
(Y\times[0,\tau]).$  The relative index formula and Theorem~\ref{thm6} imply
that $\Ind(\eth^{\even}_{\hX^{\tau}_0},\cR^{\even}_{+0})$ is independent of $\tau.$
The boundary of $\hX^1_0$ is cylindrical and isomorphic to a neighborhood of the
(flattened) boundary of $X_1,$ hence we can glue $\hX^1_0$ to $\overline{X_1}$ to
obtain
\begin{equation}
\tX_{01}=\hX^1_0\amalg_{Y\times\{1\}} \overline{X_1}.
\end{equation}
This is a manifold with a $\spnc$-structure and Dirac operator
$\eth^{\even}_{\tX_{01}}.$ Let $\cP^{\even}_{+0}$ denote Calderon projector on
$b\hX^1_0$ defined by the invertible double construction and $\cP^{\even}_{+1}$
that defined on $bX_1$ via the invertible double construction.  Since these two
manifolds agree in neighborhoods of their respective boundaries, it is clear
that $\cP^{\even}_{+0}-\cP^{\even}_{+1}$ is a smoothing operator.

We can use the boundary projector $\cR^{\even\prime}_{+1}$ to define a boundary
condition on $\hX^{1}_0.$ Because
$\Ind(\eth^{\even}_{\hX^1_0},\cR^{\even}_{+0})=
\Ind(\eth^{\even}_{X_0},\cR^{\even}_{+0}),$ the Agranovich--Dynin formula (Theorem~\ref{thm8})
implies that
\begin{equation}
\begin{split}
\Rind(\cS_0,\cS_1')&=-\Ind(\eth^{\even}_{\hX^1_0},\cR^{\even}_{+0})+
\Ind(\eth^{\even}_{\hX^1_0},\cR^{\even\prime}_{+1})\\
&=-\Ind(\eth^{\even}_{X_0},\cR^{\even}_{+0})+
\Ind(\eth^{\even}_{\hX^1_0},\cR^{\even\prime}_{+1}),
\end{split}
\label{5.18.211}
\end{equation}
and therefore
\begin{multline}
\Rind(\cS_0,\cS_1')+
\Ind(\eth^{\even}_{X_0},\cR^{\even}_{+0})-
\Ind(\eth^{\even}_{X_1},\cR^{\even\prime}_{+1})=\\
\Ind(\eth^{\even}_{\hX^1_0},\cR^{\even\prime}_{+1})-
\Ind(\eth^{\even}_{X_1},\cR^{\even\prime}_{+1}).
\label{5.18.22}
\end{multline}
The relative index formula (Theorem~\ref{thm7}) implies that
\begin{equation}
\begin{split}
\Ind(\eth^{\even}_{\hX^1_0},\cR^{\even\prime}_{+1})&=
\Rind(\cP^{\even}_{+0},\cR^{\even\prime}_{+1})\\
\Ind(\eth^{\even}_{X_1},\cR^{\even\prime}_{+1})&=
\Rind(\cP^{\even}_{+1},\cR^{\even\prime}_{+1}).
\end{split}
\end{equation}
To complete the proof we need to show that
\begin{equation}
\Rind(\cP^{\even}_{0+},\cR^{\even\prime}_{+1})-
\Rind(\cP^{\even}_{+1},\cR^{\even\prime}_{+1})=
\Rind(\cP^{\even}_{+0},\cP^{\even}_{+1})
\label{6.26.1}
\end{equation}
This is done in Proposition~\ref{prp13}. 
As $\cP^{\even}_{+1}=\Id-\cP^{\even}_{-1}$ it follows that
\begin{equation}
\Rind(\cP^{\even}_{+0},\cP^{\even}_{+1})=\Rind(\cP^{\even}_{+0},\Id-\cP^{\even}_{-1}).
\end{equation}
Applying Bojarski's formula, (see Theorem 24.1 in~\cite{BBW}), we  therefore conclude that:
\begin{equation}
\Rind(\cP^{\even}_{+0},\cP^{\even}_{+1})=\Ind(\eth^{\even}_{\tX_{01}}).
\label{5.18.23}
\end{equation}
Theorem~\ref{thm9} follows from~\eqref{5.18.22}--\eqref{5.18.23}.
The proof of~\eqref{6.26.1} is essentially the same as the proof of the
Agranovich--Dynin formula.

\begin{proposition}\label{prp13} Let $X_0, X_1$ be normalized strictly
  pseudoconvex $\spnc$-manifolds with cylindrical ends. Suppose that a collar
  neighborhood of $bX_0$ is isomorphic to a collar neighborhood of $bX_1.$ Let
  $\cR^{\even}_+$ be a modified $\dbar$-Neumann condition defined on
  $bX_0\simeq bX_1$ and $\cP^{\even}_{+j}, j=0,1$ are the Calderon projectors
  defined via the invertible double construction on $X_j, j=0,1.$ We have that
\begin{equation}
\Rind(\cP^{\even}_{0+},\cR^{\even}_+)-\Rind(\cP^{\even}_{+1},\cR^{\even}_+)=\\
\Rind(\cP^{\even}_{+0},\cP^{\even}_{+1}).
\label{5.19.2}
\end{equation}
\end{proposition}
\begin{proof}
We give an outline for the proof. As in the proof of Theorem~\ref{thm8}, we
consider the map
\begin{equation}
A=\cP^{\even}_{+1}\cR_+^{\even}\cP^{\even}_{+0}=\cP^{\even}_{+1}\cT^{\even
  *}_{1+}
\cT^{\even}_{+0}\cP^{\even}_{+0}.
\end{equation}
\begin{enumerate}
\item[Step 1]
We show that $A$ is a Fredholm map from $L^2\cap\range\cP^{\even}_{+0}$ to
$\cK\cap\range\cP^{\even}_{+1},$ with index 
$$\Rind(\cP^{\even}_{0+},\cR_+^{\even})-\Rind(\cP^{\even}_{+1},\cR_+^{\even}).$$
Here $\cK$ is the closure of $\CI$ with respect to the inner product defined by
\begin{equation}
\|u\|^2_{\cK}=\|u\|^2_{L^2}+\|\cU^{\even}_{+0}\cU^{\even *}_{+1}u\|^2_{L^2}.
\end{equation}
\item[Step 2]
We then observe that, as $\cP^{\even}_{+0}-\cP^{\even}_{+1}$ is a smoothing
operator, so is $\cT^{\even}_{+0}-\cT^{\even}_{+1}.$ Hence the index of $A$
equal to that of
\begin{equation}
B=\cP^{\even}_{+1}\cT^{\even *}_{1+}\cT^{\even}_{+1}\cP^{\even}_{+0}:
L^2\range\cP^{\even}_{+0}\to
\cK\cap\range\cP^{\even}_{+1}.
\end{equation}

\item[Step 3] As before, the commutator $[\cP^{\even}_{+1},\cT^{\even
*}_{1+}\cT^{\even}_{+1}]=0,$ hence we can replace $B$ with
\begin{equation}
B=\cP^{\even}_{+1}\cT^{\even
  *}_{1+}\cT^{\even}_{+1}\cP^{\even}_{+1}\cP^{\even}_{+0},
\end{equation}
which we think of as a composition of Fredholm maps
\begin{equation}
\begin{split}
&\cP^{\even}_{+1}\cP^{\even}_{+0}:L^2\cap\range\cP^{\even}_{+0}\to
L^2\cap\range\cP^{\even}_{+1}
\text{ and }\\
\cW=&\cP^{\even}_{+1}\cT^{\even *}_{1+}\cT^{\even}_{+1}\cP^{\even}_{+1}:
L^2\cap\range\cP^{\even}_{+1}\to\cK\cap\range\cP^{\even}_{+1}.
\end{split}
\end{equation}
The index of the first term is $\Rind(\cP^{\even}_{+0},\cP^{\even}_{+1}).$ 
\item[Step 4] To complete the proof we need to
show that $\Ind(\cW)=0.$  Again, this is a formally self adjoint operator, so
the vanishing of the index follows from the trace formula.
\end{enumerate}

Step 1 essentially follows from Theorem~\ref{thm11} in the Appendix. For
clarity we outline the argument.  To show that $A$ is a Fredholm map, as
indicated, follows easily from the commutation relations
\begin{equation}
\begin{split}
(\Id+K)\cP^{\even}_{+0}\cU^{\even}_{+0}&=\cU^{\even}_{+0}\cR_+^{\even}(\Id+K)\\
(\Id+K)\cR_+^{\even}\cU^{\even*}_{+1}&=\cU^{\even *}_{+1}\cP^{\even}_{+1}(\Id+K),
\end{split}
\label{5.19.6}
\end{equation}
where we use $K$ to denote a variety  of smoothing operators. To compute its
index we factor it as the composition of
\begin{equation}
\begin{split}
&\cR_+^{\even}\cP^{\even}_{+0}:L^2\cap\range\cP^{\even}_{+0}
\to\cH\cap\range\cR_+^{\even}\text{
    and }\\
&\cP^{\even}_{+1}\cR_+^{\even} :\cH\cap\range\cR_+^{\even}
\to\cK\cap\range\cP^{\even}_{+1}.
\end{split}
\label{5.19.5}
\end{equation}
Here $\cH$ is the closure of $\CI$ with respect to the inner product defined by
\begin{equation}
\|u\|^2_{\cH}=\|u\|^2_{L^2}+\|\cU^{\even}_{+0}u\|^2_{L^2},
\end{equation}
and  $\cK$ is the closure of $\CI$ with respect to the inner product defined by
\begin{equation}
\|u\|^2_{\cK}=\|u\|^2_{L^2}+\|\cU^{\even}_{+0}\cU^{\even *}_{+1}u\|^2_{L^2}.
\end{equation}
That $\cR_+^{\even}\cP^{\even}_{+0}$ is a Fredholm map with respect to these
spaces and has index 
$$\Rind(\cP^{\even}_{+0},\cR_+^{\even})$$ 
follows immediately from the results in the appendix.

We now show that the second map in~\eqref{5.19.5} is Fredholm and has index
$$-\Rind(\cP^{\even}_{+1},\cR_+^{\even}).$$
Using the commutation relations,~\eqref{5.19.6}, it follows easily that this
map is bounded.  The commutation relations imply that the map
$\cR_+^{\even}\cU^{\even *}_{+1}\cP^{\even}_{+1}:\cK\to\cH$ is bounded and satisfies
\begin{equation}
\begin{split}
&\cP^{\even}_{+1}\cR_+^{\even}[\cR_+^{\even}\cU^{\even *}_{+1}\cP^{\even}_{+1}]=
\cP^{\even}_{+1}(\Id-K_1)\cP^{\even}_{+1}\text{ and }\\ &[\cR_+^{\even}\cU^{\even
*}_{+1}\cP^{\even}_{+1}]\cP^{\even}_{+1}\cR_+^{\even}=
\cR_+^{\even}(\Id-K_2)\cR_+^{\even}
\end{split}
\end{equation}
for $K_1, K_2$ smoothing operators. This shows that
$\cP^{\even}_{+1}\cR_+^{\even}$ is Fredholm.  As follows from the results in
the appendix, the index of this operator is given by
\begin{equation}
\Ind(\cP^{\even}_{+1}\cR_+^{\even})=\Tr\cP^{\even}_{+1}K_1\cP^{\even}_{+1}-
\Tr\cR_+^{\even}K_2\cR_+^{\even},
\end{equation}
where as usual, we can use the $L^2$-topology to compute the traces. Comparing
this to~\eqref{eqna14} we see that
\begin{equation}
\Ind(\cP^{\even}_{+1}\cR_+^{\even})=-\Rind(\cP^{\even}_{+1},\cR_+^{\even}).
\end{equation}
As $\cP^{\even}_{+1}\cR_+^{\even}\cP^{\even}_{+0}$ is the composition of the
Fredholm maps in~\eqref{5.19.5}, this completes step 1. 

Step 2 is obvious as
\begin{equation}
A-B=\cP^{\even}_{+1}K\cP^{\even}_{+0},
\end{equation}
for $K$  a smoothing operator. We now turn to step 3. Using the commutation
relations in~\eqref{5.19.6}, we easily show that $\cW$ is a Fredholm map with
parametrix $\cV=\cP^{\even}_{+1}\cU^{\even}_{+1}\cU^{\even
  *}_{+1}\cP^{\even}_{+1}.$ As before, if
\begin{equation}
\cV\cW=\cP^{\even}_{+1}-\cP^{\even}_{+1}K\cP^{\even}_{+1},
\end{equation}
for $K$ a smoothing operator, then
\begin{equation}
\cW\cV=\cP^{\even}_{+1}-\cP^{\even}_{+1}K^*\cP^{\even}_{+1}.
\end{equation}
Step 4 is completed, as before, by using the trace formula
\begin{equation}
\Ind(\cW)=\Tr\cP^{\even}_{+1}K\cP^{\even}_{+1}-
\Tr\cP^{\even}_{+1}K^*\cP^{\even}_{+1}.
\end{equation}
Because the right had side is a purely imaginary number, this relation implies
that the index is zero. This completes the proof of the proposition.
\end{proof}
With this proposition the proof of Theorem~\ref{thm9} is complete.  As noted
above, the Atiyah--Weinstein conjecture is an immediate consequence of
Theorem~\ref{thm9}.
\end{proof}

\section{Vector Bundle Coefficients}\label{sec11}
The foregoing analysis applies equally well if we consider the $\spnc$-Dirac
operators acting on sections of $E\otimes\Spn,$ where $E\to X$ is a complex
vector bundle. The results in this paper rest entirely on the properties of the
principal symbols of the comparison operators, $\cT^{\eo}_+.$ The computations
of these principal symbols follow from
equations~\eqref{5.23.1}--\eqref{5.23.2}.  These, in turn, are consequences of
the geometric statements in equations~\eqref{5.23.3}--\eqref{5.23.4}, as well
as~\eqref{7.2.1}. These geometric normalizations continue to be possible if we
twist the spin bundle with a Hermitian vector bundle.

Let $\nabla^{\Spn}$ denote the compatible connection on the $\spnc$-bundle and
$\nabla^E,$ an Hermitian connection on $E.$ The connection
\begin{equation}
\nabla^{\Spn\otimes E}=\nabla^{\Spn}\otimes\Id_E+\Id_{\Spn}\otimes \nabla^{E},
\end{equation}
is a compatible connection on $\Spn\otimes E.$ We fix a point $p\in bX,$ and
let $(x_1,\dots,x_{2n})$ be normal coordinates at $p.$ Let $\{\sigma_J\}$ be a
local frame field for $\Spn,$ satisfying~\eqref{5.23.4}, and $\{e_l\}$ a local
framing for $E$ with
\begin{equation}
\nabla^E e_l=O(|x|).
\label{5.23.8}
\end{equation}
In this case $\{\sigma_J\otimes e_l\}$ is a local framing for $\Spn\otimes E$
that satisfies
\begin{equation}
\nabla^{\Spn\otimes E}\sigma_J\otimes e_l=O(|x|).
\label{5.23.9}
\end{equation}

Let $\eth_E$ be the $\spnc$-Dirac operator acting on sections of $\Spn\otimes
E.$ Because the coordinates are normal,~\eqref{5.23.9} implies that
\begin{equation}
\begin{split}
\eth_E\sum_{J,l} a_{J,l}\sigma_J\otimes e_l=&
\ha\sum_{J,l}\sum_{j=1}^{2n}\bc(dx_j)\cdot
\left[(\pa_{x_j}a_{J,l})\sigma_J\otimes e_l+
a_{J,l}\nabla^{\Spn\otimes E}\sigma_J\otimes e_l\right]\\
=&\sum_{J,l}\sum_{j=1}^{2n}[\dbar+\dbar^*]_{\bbC^n}[a_{J,l}d\bz^J\otimes
  e_l]+\sO_1(|x|)+\sO_0(|x|).
\end{split}
\end{equation}
It is a general result about Dirac operators that $[\eth^E]^2=\Delta+R,$ where
$\Delta$ is the Laplace operator and $R$ is an operator of order zero. We
compute the action of the Laplace operator in the normal coordinates at $p:$
\begin{equation}
\begin{split}
{\Delta } &\sum_{J,l} a_{J,l}\sigma_J\otimes e_l=\\
&\sum_{j,k=1}^{2n}(\delta_{jk}+O(|x|^2))[\nabla^{\Spn\otimes E}_{\pa_{x_j}}
\nabla^{\Spn\otimes E}_{\pa_{x_k}} a_{J,l}\sigma_J\otimes e_l-
\nabla^{\Spn\otimes E}_{\nabla^g_{\pa_{x_j}}\pa_{x_k}}
a_{J,l}\sigma_J\otimes e_l] \\
=&\sum_{J,l} \sum_{j=1}^{2n}\pa_{x_j}^2a_{J,l}\sigma_J\otimes
e_l+\sO_2(|x|^2)+
2\pa_{x_j}a_{J,l}\nabla^{\Spn\otimes E}_{\pa_{x_j}}
\sigma_J\otimes_l+\\&\sO_1(|x|)+\sO_0(1).
\end{split}
\end{equation}
Using once again that $\nabla^{\Spn\otimes E}_{\pa_{x_j}}\sigma_J\otimes
e_l=O(|x|)$ we see that 
\begin{equation}
{\Delta } \sum_{J,l} a_{J,l}\sigma_J\otimes e_l=
\sum_{J,l} \sum_{j=1}^{2n}\pa_{x_j}^2a_{J,l}\sigma_J\otimes
e_l+\sO_2(|x|^2)+\sO_1(|x|)+\sO_0(1).
\end{equation}

These formul{\ae} demonstrate that the necessary symbolic conditions are
satisfied by $\eth_E$ and $\eth^2_E.$ As described in~\cite{Epstein3}, complex
vector bundle coefficients are easily incorporated into the generalized 
Szeg\H o projector formalism,  and therefore the results proved in the previous
sections apply equally well when vector bundle coefficients are included. We do
not wish to exhaustively enumerate these generalizations, but simply list a few
of these results.

We let $\cP^{\eo}_{\pm E}$ denote the Calderon projectors with bundle
coefficients and $\cR^{\eo}_{+E}$ the modified $\dbar$-Neumann conditions
defined by a generalized Szeg\H o projector $\cS_{E}.$ As before we set
\begin{equation}
\cT^{\eo}_{+E}=\cR^{\eo}_{+E}\cP^{\eo}_{+E}+(\Id-\cR^{\eo}_{+E})(\Id-\cP^{\eo}_{+E}).
\end{equation}
The basic analytic result is:
\begin{theorem}\label{thm1E} Let $(X,J,g,\rho)$ be a normalized strictly pseudoconvex
  $\spnc$-mani-fold, $(E,h)\to X$ a Hermitian vector bundle and $\cS_E$ a
  generalized Szeg\H o projector acting on sections of $E.$ The principal
  symbol of $\cS_E$ is defined by a compatible deformation of the almost
  complex structure on $H$ induced by the embedding of $bX$ as the boundary of
  $X.$ Then the operators $\cT^{\eo}_{+E}$ are elliptic, in the extended
  Heisenberg calculus, with parametrices having Heisenberg orders
\begin{equation}
\left(\begin{matrix} 0 & 1\\ 1 & 1\end{matrix}\right).
\end{equation}
\end{theorem}

As before, this result shows that the graph closures of
$(\eth^{\eo}_{+E},\cR^{\eo}_{+E})$ are Fredholm and
\begin{equation}
(\eth^{\eo}_{+E},\cR^{\eo}_{+E})^*=\overline{(\eth^{\ooee}_{+E},\cR^{\ooee}_{+E})}.
\end{equation}
Moreover, the $(\cP^{\eo}_{+E},\cR^{\eo}_{+E})$ are tame Fredholm pairs with
\begin{equation}
\Ind(\eth^{\eo}_{+E},\cR^{\eo}_{+E})=\Rind(\cP^{\eo}_{+E},\cR^{\eo}_{+E}).
\end{equation}
We have the Agranovich-Dynin Formula:
\begin{theorem}\label{thm88} Let $(X,J,h,\rho)$ be a normalized strictly pseudoconvex
  $\spnc$-mani-fold, and $E\to X$ a Hermitian vector bundle. Let
  $\cS_{E1}$ and $\cS_{E2}$ be two generalized Szeg\H o projections and
  $\cR^{\even}_{+E1}, \cR^{\even}_{+E2}$ the modified $\dbar$-Neumann conditions
  they define, then
\begin{equation}
\Rind(\cS_{E1},\cS_{E2})=\Ind(\eth^{\even}_{+}E,\cR^{\even}_{+E2})-
\Ind(\eth^{\even}_{+E},\cR^{\even}_{+E1}). 
\label{5.18.44}
\end{equation}
\end{theorem}

Finally, we have the Atiyah-Weinstein conjecture for this case. Over $X_0, X_1$
we have bundles $E_0\to X_0, E_1\to X_1.$ If $\Phi$ denotes the extension of
the contact diffeomorphism to a neighborhood of $bX_j,$ then we need to assume
that, in the collar neighborhood, $\Phi^*E_1$ is isomorphic to $E_0,$ via an
bundle map $\Psi.$  Altogether we get a vector bundle $E\to\tX_{01},$ which may
depend on the choice of $\Psi.$

\begin{theorem}\label{thm99} Let $(X_k,j_k,g_k,\rho_k), k=0,1$ be normalized
  strictly pseudoconvex $\spnc$-manifolds, with $E_j\to X_j, j=0,1$ Hermitian
  vector bundles. Suppose that $\phi:bX_1\to bX_1$ is a co-orientation
  preserving contact diffeomorphism, and $\Psi:E_0\to \Phi^*E_1,$ is a bundle
  equivalence, covering $\phi.$ Suppose that $\cS_{E_0}, \cS_{E_1}$ are
  generalized Szeg\H o projections on $bX_0, bX_1,$ respectively, which define
  modified $\dbar$-Neumann conditions $\cR^{\even}_{+E_j}$ on $X_j.$ If we let
  $\cS_{E_1}'=\Psi^{-1}\cS_{E_1}\Psi,$ $\tX_{01}$ as defined
  in~\eqref{5.19.10}, and $E$ the bundle over $\tX_{01}$ defined by gluing
  $E_0$ to $\Phi^* E_1$ via $\Psi,$ then we have
\begin{equation}
\Rind(\cS^{\phantom{\prime}}_{E_0},\cS^{\prime}_{E_1})=
\Ind(\eth^{\even}_{\tX_{01}E})-
\Ind(\eth^{\even}_{X_0E_0},\cR^{\even}_{+E_0})+
\Ind(\eth^{\even}_{X_1E_1},\cR^{\even}_{+E_1}).
\label{5.18.200}
\end{equation}
\end{theorem}

\section{The Relative Index Conjecture}\label{s.3d}

In~\cite{Epstein} we introduced the relative index for pairs of embeddable
CR-structures on 3-dimensional manifolds with the same underlying contact
structure. In those papers we used the opposite convention to that employed in
the current series of papers and therefore the relative index defined there is
minus that defined here.  With the present convention, Proposition 8.1
in~\cite{Epstein} implies that if $\cS_0$ is the Szeg\H o projector defined by
a ``reference'' embeddable CR-structure ${}^0T^{0,1}Y$ on $(Y,H),$ and $\cS_1$
is the Szeg\H o projector defined by a \emph{sufficiently small}, embeddable
deformation, ${}^1T^{0,1}Y,$ of this CR-structure, then
\begin{equation}
\Rind(\cS_0,\cS_1)\geq 0.
\end{equation}
In~\cite{Epstein} we showed that, for $n\in\bbN,$ the set of embeddable
deformations of ${}^0T^{0,1}Y$ that satisfy $\Rind(\cS_0,\cS_1)\leq n$ is
closed in the $\CI$-topology. This motivated our relative index conjecture,
which asserts (with our current sign convention) that
$\cS_1\to\Rind(\cS_0,\cS_1)$ is bounded from above, among sufficiently
small, embeddable deformations, ${}^1T^{0,1}Y,$ of ${}^0T^{0,1}Y.$
In~\cite{Epstein5} we establish this conjecture for CR-structures that bound
pseudo\emph{concave} manifolds $X_-,$ satisfying either
\begin{equation}
\begin{split}
H^2_c(X_-;\Theta\otimes &[-Z])=0\text{ or, }\\
H^2_c(X_-;\Theta)=0&\text{ and }H^1(Z;N_Z)=0,
\end{split}
\end{equation}
here $Z\subsubset X_-$ is a smooth, compact holomorphic curve with positive
normal bundle, and $\Theta$ is the tangent sheaf of $X_-.$

Suppose that $j_0, j_1$ define CR-structures on a 3-dimensional contact
manifold $(Y,H),$ which bound strictly pseudoconvex complex manifolds
$(X_0,J_0),$ $(X_1,J_1).$ If $\cS_0,\cS_1$ are the classical Szeg\H o
projectors defined by $j_0,j_1,$ respectively, then formula~\eqref{5.18.21}
gives:
\begin{equation}
\Rind(\cS_0,\cS_1)=\Ind(\eth^{\even}_{\tX_{01}})+\dim H^{0,1}(X_0,J_0)-\dim
H^{0,1}(X_1,J_1).
\label{7.14.104}
\end{equation}
The Atiyah-Singer index theorem provides a cohomological formula for
$\Ind(\eth^{\even}_{\tX_{01}}):$
\begin{equation}
\Ind(\eth^{\even}_{\tX_{01}})=\frac{c_1^2(\Spn^{\even})[\tX_{01}]-\sign[\tX_{01}]}{8}.
\label{7.14.100}
\end{equation}
Here $\Spn^{\even}$ is the $+$-spinor bundle defined on $\tX_{01}$ by the
extended  double construction, and $\sign[M]$ is the signature of the
oriented 4-manifold $M.$

In~\cite{OzbagciStipsicz2} a general formula is given relating the
characteristic numbers on a compact $\spnc$ 4-manifold, $M$:
\begin{equation}
c_2(\Spn^{\even})[M]=\frac{c_1^2(\Spn^{\even})[M]-3\sign[M]-2\chi[M]}{4},
\label{7.14.40}
\end{equation}
Here $[M]$ denotes the fundamental class of the oriented manifold $M.$ This
formula is stated as an exercise, whose solution we briefly explain: One first
shows that, if $L$ is a Hermitian line bundle, then the Chern-Weil
representative of $c_2(\Spn^{\even}\otimes L)-c_1^2(\Spn^{\even}\otimes L)/4$ does not depend
$L.$ Hence we can \emph{locally} represent $\Spn^{\even}\otimes L$ as $\Spn^{\even}_{0},$ where
$\Spn^{\even}_{0}$ is the $+$-spinor bundle coming from a \emph{locally defined}
Spin-bundle. Using the expression for the curvature of $\Spn^{\even}\otimes L$ that
arises from such a local representation, we show that the Chern-Weil
representative of $c_2(\Spn^{\even})-c_1^2(\Spn^{\even})/4$ agrees with that of
$-(p_1(M)+2e(M)),$ where $p_1$ is the first Pontryagin class and $e$ is the
Euler class.

Putting~\eqref{7.14.40} into~\eqref{7.14.100} gives:
\begin{equation}
\Ind(\eth^{\even}_{\tX_{01}})=\frac{2c_2(\Spn^{\even})[\tX_{01}]+\sign[\tX_{01}]+
\chi[\tX_{01}]}{4}. 
\label{7.14.101}
\end{equation}
Note that on $X_0,$ $\Spn^{\even}\simeq \Lambda^{0,0}X_0\oplus\Lambda^{0,2}X_0,$ and on
$X_1,$ $\Spn^{\even}\simeq\Lambda^{0,1}X_1.$ Since $\Spn^{\even}$ has a global section over $X_0$
and over the part of $\tX_{01}$ coming from the neck joining $X_0$ to $X_1,$
and a neighborhood of the boundary of $X_1,$ we can choose a metric for $\Spn^{\even}$
so that the Chern-Weil representative of $c_2(\Spn^{\even})$ is supported in the
interior of $X_1.$ Over $X_1$
$$\Spn^{\even}=\Lambda^{0,1} X_1=[T^{0,1}X_1]^*\simeq T^{1,0}X_1,$$ 
and therefore it follows that $c_2(\Spn^{\even})\restrictedto_{X_1}=e(X_1).$  Recalling
that the orientation of $X_1$ is reversed in $\tX_{01},$ and using the
additivity of the signature and Euler characteristic we obtain
\begin{equation}
\begin{split}
\Ind(\eth^{\even}_{\tX_{01}})&=\frac{2e(X_1)[-X_1]+\sign[X_0]-
\sign[X_1]+\chi[X_0]+\chi[X_1]}{4}\\
&=\frac{\sign[X_0]-\sign[X_1]+\chi[X_0]-\chi[X_1]}{4}.
\label{7.14.102}
\end{split}
\end{equation}
Using this formula in~\eqref{7.14.104} completes the proof of the following
theorem.
\begin{theorem}\label{thm13} Let $(Y,H)$ be a compact 3-dimensional co-oriented, contact
  manifold, and let $j_0,j_1$ be CR-structures with underlying plane field $H.$
  Suppose that $(X_0,J_0),$ $(X_1,J_1)$ are strictly pseudoconvex complex
  manifolds with boundary $(Y,H,j_0),$ $(Y,H,j_1),$ respectively. If $\cS_0,$
  $\cS_1$ are the classical Szeg\H o projectors defined by these CR-structures
  then
\begin{equation}
\begin{split}
\Rind(\cS_0,\cS_1)=&\dim H^{0,1}(X_0,J_0)-\dim
H^{0,1}(X_1,J_1)+\\
&\frac{\sign[X_0]-\sign[X_1]+\chi[X_0]-\chi[X_1]}{4}.
\end{split}
\label{7.14.106}
\end{equation}
\end{theorem}

\begin{remark} If $(M,\Spn^{\even},\Spn^{\odd})$ is a compact $\spnc$ 4-manifold, then
\begin{equation}
d_{\SW}(M)=\frac{c_1^2(\Spn^{\even})[M]-3\sign[M]-2\chi[M]}{4}
\end{equation}
is the formal dimension of the moduli space of solutions to the Seiberg-Witten
equations, see~\cite{Morgan}. Our calculations show that for manifolds
$\tX_{01}\simeq X_0\amalg\overline{X_1},$ with $\spnc$-structure defined by the
invertible double construction,
\begin{equation}
d_{\SW}(\tX_{01})=-\chi[X_1].
\end{equation}
Reversing the orientation of $\tX_{01}$ gives $\tX_{10}$ and interchanges $\Spn^{\even}$
with $\Spn^{\odd},$ so that
\begin{equation}
d_{\SW}(\tX_{10})=-\chi[X_0].
\end{equation}
Note finally that, under the hypotheses of Theorem~\ref{thm13},
equation~\eqref{7.14.106} implies that
\begin{equation}
\frac{\sign[X_0]-\sign[X_1]+\chi[X_0]-\chi[X_1]}{4}\in\bbZ.
\end{equation}
\end{remark}
Formula~\eqref{7.14.106} has a direct bearing on the relative index conjecture.
\begin{corollary}\label{cor8}
 Let $(Y,H)$ be a compact 3-dimensional, co-oriented contact manifold. Suppose
  that among Stein manifolds $(X,J)$ with pseudoconvex boundary $(Y,H)$ the
  signature, $\sign(X)$ and Euler characteristic $\chi(X)$ assume only finite
  many values, then, among embeddable deformations, $\cS_1$ of a given
  embeddable reference CR-structure, $\cS_0,$ the relative index
  $\Rind(\cS_0,\cS_1)$ is bounded from above.
\end{corollary}
\begin{proof} Suppose that $X$ is a strictly pseudoconvex complex manifold with
  boundary $(Y,H).$ We can assume that $X$ is minimal, i.e. all inessential
  compact varieties are blown down. Bogomolov and DeOliveira proved that one
  can deform the complex structure on $X$ to obtain a Stein manifold,
  see~\cite{BogomolovDeOliveira}. Such a deformation does not change the
  topological invariants $\sign[X],$ $\chi[X].$ Thus among minimal strictly
  pseudoconvex complex manifolds $X,$ with boundary $(Y,H),$ the numbers
  $\sign[X],$ $\chi[X]$ assume only finitely many values, and the corollary
  follows from~\eqref{7.14.106}.
\end{proof}

If $X_0$ is diffeomorphic to $X_1,$  then~\eqref{7.14.106} implies that
\begin{equation}
\Rind(\cS_0,\cS_1)=\dim H^{0,1}(X_0,J_0)-\dim H^{0,1}(X_1,J_1).
\label{5.19.21}
\end{equation}
For such deformations of the CR-structure  we see that
\begin{equation}
\Rind(\cS_0,\cS_1)\leq \dim H^{0,1}(X_0,J_0).
\end{equation}
This becomes an equality if $(X,J_1)$ is a Stein manifold. It says that a
singular Stein surface, with $H^{0,1}(X_0,J_0)\neq 0,$ has a larger algebra of
holomorphic functions than its smooth (Stein) deformations.  If both
$(X_0,J_0)$ and $(X_1,J_1)$ are Stein manifolds, then the relative index is
zero.

For some time it was believed that a given compact, contact 3-manifold should
have, at most, finitely many diffeomorphism classes of Stein fillings.  Indeed
this expectation has been established for some classes of contact
3-manifolds. Eliashberg proved that the 3-sphere, with its standard contact
structure, has a unique Stein filling. McDuff extended this to the lens spaces
$L(p,1)$ for $p\neq 4,$ with contact structure induced from $S^3.$ Stipsicz
showed the uniqueness of the Stein filling for the only fillable contact
structure on the 3-torus, see~\cite{Stipsicz2}. Lisca showed
that other lens spaces have finitely many diffeomorphism classes of fillings,
see~\cite{Lisca}. Ohta and Ono proved the finiteness statement for simple and
simple elliptic isolated singularities, see~\cite{OhtaOno1,OhtaOno2}. For most
of these cases, the local form of the relative index conjecture was proved by
other means in~\cite{Epstein5}. 

Recently, I. Smith and Ozbagci and Stipsicz have produced examples of
compact, contact 3-manifolds that arise as the boundaries of infinitely many
diffeomorphically distinct Stein manifolds, see~\cite{OzbagciStipsicz}.  In the
examples of Ozbagci and Stipsicz, the signatures and Euler characteristics of
the 4-manifolds are all equal. In a later paper Stipsicz conjectured that,
given a 3-dimensional compact, contact manifold, $(Y,H),$ the signatures and
Euler characteristics for Stein manifolds with boundary $(Y,H)$ should assume
only finitely many distinct values, see~\cite{stipsicz}. Using
Corollary~\ref{cor8}, this conjecture would clearly imply the strengthened form
of the relative index conjecture.

Stipsicz has established his conjecture in a variety of cases, among them, the
boundaries of circle bundles in line bundles of degree $n$ over surfaces of
genus $g,$ provided
\begin{equation}
|n|>2g-2.
\end{equation}
This allows us to extend the result proved in~\cite{Epstein5}, where the
relative index conjecture was established for circle bundles under the
assumption that $|n|\geq 4g-3.$ Stipsicz also proved his conjecture for the
Seifert fibered 3-manifold $\Sigma(2,3,11),$ and so the relative index
conjecture holds in this case as well. 

Even if the Stipsicz conjecture is false, it would not necessarily invalidate
our conjecture. The relative index conjecture is a conjecture for
\emph{sufficiently small} embeddable deformations of the CR-structure. It could
well be that as one moves through the infinitely many diffeomorphism types,
arising among the fillings, the deformations of the CR-structure on the
boundary become large. Conversely, when such bounds do exist, a global upper
bound for the relative index follows from~\eqref{7.14.106} and
Corollary~\ref{cor8}.

\section{Further Interesting Special Cases}
We consider some other special cases of the index formul{\ae} proved above. 

\subsection{Co-ball Bundles}
First we
consider the original Atiyah-Weinstein conjecture which concerns pairs of
co-ball bundles, $X_0=B^*M_0, X_1=B^*M_1$ and a contact diffeomorphism of their
boundaries, $\phi:S^* M_1\to S^* M_0.$ In this case there is a complex structure
on each of the manifolds, well defined up to isotopy. In these structures,
$X_0, X_1$ are Stein manifolds and therefore
\begin{equation}
\chi'_{\cO}(X_0)=\chi'_{\cO}(X_1)=0.
\label{6.28.1}
\end{equation}
If $\cS_0,$ $\cS_1$ are the classical Szeg\H o projectors defined by the
complex structures on $X_0, X_1,$ respectively, then equations~\eqref{5.18.21}
and~\eqref{6.28.1} imply that, with $\cS_1'=[\phi^{-1}]^*\cS_1\phi^*,$ 
\begin{equation}
\Rind(\cS_0,\cS_1')=\Ind(\eth^{\even}_{X_{\phi}}),
\end{equation}
where
\begin{equation}
X_\phi=B^*M_0\amalg_{\phi} \overline{B^* M_1}.
\end{equation}

As noted in the introduction, the mapping $\phi$ and the complex structures on
the ball bundles define a class of elliptic Fourier integral operators. Let
$F^\phi:\CI(M_0)\to\CI(M_1)$ be an element of this class.  From the definition
of $F^{\phi}$ it is clear that its index equals $\Rind(\cS_0,\cS_1'),$ and
therefore
\begin{equation}
\Ind(F^\phi)=\Ind(\eth^{\even}_{X_{\phi}}).
\end{equation}
Let $\Phi$ denote the homogeneous extension of $\phi$ to
$B^*M_1\setminus\{0\}.$ If $\Phi$ extends smoothly across the zero section,
i.e., $\phi$ is defined as the differential of a diffeomorphism  $f:M_0\to
M_1,$ then the glued space $X_\phi$ is essentially an invertible double. Hence
$\Ind(\eth^{\even}_{X_{\phi}})=0,$ and we have:
\begin{proposition} If $\Phi$ extends smoothly across the zero section of
  $B^*M_1,$ then $\Ind F^{\phi}=0.$
\end{proposition}
\noindent
The case of a pair of co-ball bundles is treated in~\cite{LeichtNestTsy}, where
equivalent results are proved.

In the 3-dimensional case, we can use formula~\eqref{7.14.102} to obtain:
\begin{equation}
\Ind(F^\phi)=\frac{\sign[X_0]-\sign[X_1]+\chi[X_0]-\chi[X_1]}{4},
\label{7.25.05.2}
\end{equation}
where $X_0=B^*M_0$ and $X_1=B^*M_1.$ Note that $X_j$ retracts onto $M_j.$  Hence
$H_2(X_j,\bbZ)$ is one dimensional, and generated by the zero section, $[M_j].$
The self intersection of the zero section in a vector bundle equals the Euler
charteristic of the bundle.  In this case it is clear that, for $j=0,1,$
\begin{equation}
\chi[X_j]=\chi[M_j]\text{ and } 
\sign[X_j]=\begin{cases}-1&\text{ if }\chi[M_j]<0\\
0&\text{ if }\chi[M_j]=0\\
1&\text{ if }\chi[M_j]>0.
\end{cases}
\label{7.25.05.1}
\end{equation}

Using~\eqref{7.25.05.2} and~\eqref{7.25.05.1} we easy prove that, in the
3-dimensional case, $\Ind(F^\phi)$ vanishes.
\begin{theorem}
Let $M_0, M_1$ be compact, oriented 2-manifolds and suppose that $\phi $ is a
co-orientation preserving contact diffeomorphism, $\phi :S^*M_1\to S^* M_0,$
then
$$\Ind(F^\phi)=0.$$
\end{theorem}
\begin{proof}
The theorem follows from the formul{\ae} above and the following lemma:
\begin{lemma}\label{lem10} If $Y$ is the co-sphere bundle of an oriented
  compact surface, $M,$ then
\begin{equation}
H_1(M;\bbQ)\simeq H_1(Y;\bbQ).
\end{equation}
\end{lemma}
\begin{proof}[Proof of lemma] The Leray spectral sequence implies that there is
  a short exact sequence:
\begin{equation}
\begin{CD}0@>>> \bbZ/\chi[M]\bbZ @>>> H_1(Y;\bbZ)@>>>H_1(M;\bbZ)@>>>0\end{CD}.
\end{equation}
Taking tensor product with $\bbQ$ leaves the sequence exact and proves the lemma.
\end{proof}
The theorem follows from the lemma as it implies that $\chi[M_0]=\chi[M_1].$
\end{proof}

\begin{remark} A similar vanishing theorem in dimension 3 is proved in~\cite{EpsteinMelrose} for
  $\phi: Y\to Y$ a contact self map.
\end{remark}

\subsection{The Atiyah-Singer Index Theorem}
Let $M$ be a compact manifold without boundary, and let $E, F$ be complex
vector bundles over $M.$ Let $P$ be an elliptic pseudodifferential operator of
order zero, $P:\CI(M;E)\to\CI(M;F).$ We can use our theorem to give an
analytic proof of the K-theoretic step in the original proof of the
Atiyah-Singer, which states that the index of $P$ equals
that of a Dirac operator on the glued space:
\begin{equation}
\bbT M=B^* M\amalg_{S^* M}\overline{B^* M}.
\end{equation}
Because $B^*_{\epsilon}M$ is a Stein manifold, Oka's principle implies that the
lifted bundles $\pi^* E, \pi^* F$ have well defined complex structures. Let
$G_{b\epsilon}^E, G_{b\epsilon}^F$ be the maps from
$\cO_b(B^*_{\epsilon}M;E),$ $\cO_b(B^*_{\epsilon}M;F)$ to
$\CI(M;E),$ $\CI(M;F),$ respectively. defined by pushforward. As before these maps are
isomorphisms for small enough $\epsilon.$ 

Let $\cS^E_\epsilon, \cS^F_{\epsilon}$ be classical Szeg\H o projectors onto
the boundary values of holomorphic sections of $\pi^* E, \pi^* F,$
respectively. If
$$\sigma_P\in \CI(S^* M;\Hom(\pi^* E,\pi^* F))$$ 
is the principal symbol of $P,$ then, for sufficiently small $\epsilon>0,$ the
composition,
\begin{equation}
F_{P} s= G_{b\epsilon}^F\cS_\epsilon^F\sigma_P\cS_\epsilon^E
[G_{b\epsilon}^{E}]^{-1}s
\end{equation}
is a pseudodifferential operator with principal symbol $\sigma_P,$
see~\cite{BoutetdeMonvel4}. The operators $G_{b\epsilon}^F, G_{b\epsilon}^E$
are invertible and therefore
\begin{equation}
\Ind(P)=\Ind(F_P)=\Ind(\cS_\epsilon^F\sigma_P\cS_\epsilon^E),
\label{6.28.2}
\end{equation}
where the last term in~\eqref{6.28.2} is the index of a Toeplitz operator. 

The Toeplitz index in~\eqref{6.28.2} is easily seen to equal the relative index
$$\Ind(\cS_\epsilon^F\sigma_P\cS_\epsilon^E)=
\Rind(\cS_\epsilon^E,[\sigma_P]^{-1}\cS_\epsilon^F\sigma_P).$$ 
This relative index is computed in Theorem~\ref{thm99}. Equation (137)
in~\cite{Epstein4} and the fact that $B^*_\epsilon M$ is a Stein manifold imply
that the boundary terms in~\eqref{5.18.200} vanish, and therefore
\begin{equation}
\Ind(P)=\Ind(\eth^{\even}_{\bbT M,\, V_P}).
\end{equation}
Here $V_P$ is the bundle obtained by gluing $\pi^* E$ to $\pi^* F$ along $S^*M$
via the symbol of $P.$ This assertion is an important step, proved using
K-theory, in the original proof of the Atiyah-Singer theorem. The relative
index formalism in this paper, along with results from~\cite{Epstein4} provide
a completely analytic proof of this statement.

\subsection{Higher Dimensional Complex Manifolds}

Suppose that $(X,J_0)$ is a strictly pseudoconvex complex $n$-dimensional
manifold, and $J_1$ is a deformation of the complex structure on $X$ that is
again strictly pseudoconvex. For our purposes it is sufficient if the
deformation takes place through a smooth family, $\{J_s:\: s\in [0,1]\},$ of
strictly pseudoconvex almost complex structures. Let $\tX_{01}$ denote the
$\spnc$-manifold obtained by gluing $(X,J_0)$ to $\overline{(X,J_1)},$ via the
extended double construction.  Because $J_0$ is homotopic to $J_1,$
it follows that the space $\tX_{01}$ is isotopic, as a $\spnc$-manifold, to the
invertible double, $\tX_0,$ of $(X,J_0).$ Using the results of Chapter 9
in~\cite{BBW} we conclude that
\begin{equation}
\Ind(\eth^{\even}_{\tX_{01}})=\Ind(\eth^{\even}_{\tX_{0}})=0.
\end{equation}
In this case,~\eqref{5.18.20} implies that
\begin{equation}
\Rind(\cS_0,\cS_1)=\chi'_{\cO}(X,J_1)-\chi'_{\cO}(X,J_0).
\label{5.19.20}
\end{equation} 
Thus for CR-structures that can be obtained by deformation of the complex
structure through almost complex structures on a fixed manifold, the relative
index is simply the change in the renormalized holomorphic Euler
characteristic. The relative index is again non-negative for small
integrable deformations. Hence, if the deformation arises
from a deformation of the almost complex structure on $X,$ then we get the
semi-continuity result for the renormalized Euler characteristic:
\begin{equation}
\chi'_{\cO}(X, J_1)\geq \chi'_{\cO}(X,J_0).
\end{equation}

\section{Appendix A: Tame Fredholm Pairs}
In this appendix we present a generalization of the theory of Fredholm pairs
and the index theory for such pairs. We give this discussion in a fairly
general functional analytic setting. We suppose that we have a nested family of
separable Hilbert spaces $(H_s,\|\cdot\|_s),$ labeled by $s\in\bbR.$ If $s<t,$
then $H_t\subset H_s$ and
\begin{equation}
\|x\|_s\leq\|x\|_t \text{ for all }x\in H_t.
\label{eqn5.16.1}
\end{equation}
The intersection
$$H_\infty=\bigcap\limits_{s=-\infty}^{\infty} H_s$$
is assumed to be dense in $H_s$ for all $s\in\bbR.$ We also define
$$H_{-\infty}=\bigcup\limits_{s=-\infty}^{\infty} H_s.$$
The inner product on $H_0$ is assumed to satisfy the
generalized H\"older inequality, for all $x,y\in H_{\infty}$ we have:
\begin{equation}
|\langle x,y\rangle|\leq \|x\|_s\|y\|_{-s}.
\label{eqna1}
\end{equation}
Indeed we suppose that $H_s'\simeq H_{-s}.$

We consider several classes of operators generalizing notions from the theory
of pseudodifferential operators.
\begin{definition}
A tame operator $T$ is an operator defined on $H_{-\infty}$  for which there is
a fixed $m\in\bbR$ such that, for all $s\in\bbR,$  
\begin{equation}
T H_s\subset H_{s-m},
\end{equation}
and the map $T:H_s\to H_{s-m}$ is bounded.
In this case we say that $T$ has order $m.$
\end{definition}
For $x,y$ in $H_{\infty}$ we define the \emph{formal adjoint}, $T^*$ of the tame
operator $T$ by duality:
\begin{equation}
\langle Tx,y\rangle=\langle x,T^*y\rangle.
\end{equation}
In fact, this extends by continuity to $x\in H_{\infty}$ and $y\in H_{-\infty}$
or $x\in H_{-\infty}$ and $y\in H_{\infty}.$
The definition of tameness implies that $Tx\in H_0,$ for $x\in H_{\infty};$ so
this notion of adjoint is consistent with the $L^2$-adjoint.  In this appendix
the notation $T^*$ always refers to the formal adjoint. 
\begin{lemma} If $T$ is a tame operator of order $m,$ then its formal adjoint
  $T^*$ is as well.
\end{lemma}
\begin{proof} We use the fact that $H_s'\simeq H_{-s}.$ Let $x,y\in H_{\infty}$
  and fix a value of $s.$ For $x,y\in H_{\infty}$ we have $\langle
  Tx,y\rangle=\langle x,T^*y\rangle.$ The generalized Cauchy-Schwarz inequality
  implies that
\begin{equation}
|\langle x,T^*y\rangle|\leq \|Tx\|_{-s}\|y\|_s\leq C\|x\|_{-s+m}\|y\|_{s}.
\end{equation}
This inequality implies that for $y\in H_{\infty}$ we have the estimate
\begin{equation}
\|T^* y\|_{s-m}\leq C\|y\|_s,
\end{equation}
as $H_{\infty}$ is dense in $H_s$ this proves the lemma.
\end{proof}
It is clear that, under composition, the set of tame operators defines a star
algebra. 

\begin{definition} A tame operator $K$ that maps $H_{-\infty}$ to $H_{\infty}$
  is called a smoothing operator. For any $s\in\bbR,$ we suppose that, when
  acting on $H_s,$ a smoothing operator is a trace class operator.
\end{definition}
A tame operator $K$ is smoothing if and only if it is a tame operator of order
$m$ for every $m\in (0,-\infty).$ It is clear that the class of smoothing
operators is closed under adjoints and defines a two sided ideal in the algebra
of tame operators.
\begin{definition} An tame operator $T$ is tamely elliptic if there is a tame
  operator $U$ so that
\begin{equation}
TU=\Id-K_1,\quad UT=\Id-K_2,
\end{equation}
where $K_1, K_2$ are smoothing operators. The operator $U$ is called a
parametrix for $T.$
\end{definition}
As the smoothing operators are closed under taking adjoints, if $T$ is tamely
elliptic, then so is $T^*.$ If $T$ is invertible in any reasonable sense and
tamely elliptic, then it follows easily that $U-T^{-1}$ is a smoothing
operator.

If $T$ is a tame operator of non-positive order, then the operator $\Id+T^*T$
is self adjoint and invertible on $H_0.$ We make the following additional
assumption: 
\smallskip
\centerline{If $T$ has non-positive order then $\Id+T^*T$ is tamely
elliptic.} 
\smallskip
This is true of any operator calculus with a good symbol map.  Tame
smoothing operators are compact operators and therefore, acting on
$H_{-\infty},$ an elliptic operator has a finite dimensional kernel contained
in $H_{\infty}.$ If $T$ is elliptic then so is $T^*$ and it also has a finite
dimensional kernel, which is contained in $H_{\infty}.$ We define the tame
index of a tamely elliptic operator to be
\begin{equation}
\tind(T)=\dim\ker T-\dim\ker T^*.
\end{equation}

We define the $H_0$-domain of tame operator $T$ to be the graph closure, in the
$H_0$-norm, of $T$ acting on $H_{\infty}.$ 
\begin{equation}
\Dom_0(T)=\{x\in H_{0} :\: \exists <x_n>\subset H_{\infty}\text{ with
}
x_n\to x\text{ and } Tx_n\to Tx \text{ in }H_0\}.
\end{equation}
As $H_{\infty}\subset \Dom_0(T)$ it is clear that $\Dom_0(T)$ is dense
in $H_0.$ 
\begin{lemma}
The operator $(T,\Dom_0(T))$  is a densely defined, closed operator. 
\end{lemma}
\begin{proof} If the order of $T$ is non-positive, then $T$ it is bounded on $H_0$ and
  therefore $\Dom_0(T)=H_0.$ Now assume that the order of $T$ is $m>0.$ Let
  $<x_n>\subset \Dom_0(T).$ Suppose that $<x_n>$ converges to $x$ in $H_0$ and
  $<Tx_n>$ converges to $y$ in $H_0.$ For each $n$ we can choose a sequence
  $<x_{nj}>\subset H_{\infty}$ so that $x_{nj}\to x_n$ and $Tx_{nj}\to Tx_n$ in
  $H_0.$ We can easily extract a subsequence $<x_{n_kj_k}>$ that converges to
  $x$ and so that $<Tx_{n_kj_k}>$ converges to $y.$ As $T$ has order $m\geq 0,$
  $<Tx_{n_kj_k}>$ converges to $Tx$ in $H_{-m}.$ By duality and the density of
  $H_{\infty},$ it is immediate that $y=Tx,$ and therefore $x\in\Dom_0(T).$
\end{proof} 

If $U$ is an operator of positive order, $\Dom_0(U)$ is a Hilbert space with
the inner product
\begin{equation}
\langle x,y\rangle_{U}=\langle x,y\rangle+\langle U x,Uy\rangle.
\end{equation}
For $x\in H_{\infty}$ this can be rewritten,
\begin{equation}
\langle x,(\Id+U^*U)y\rangle,
\end{equation}
where $(\Id+U^*U)y$ is an element of $H_{-\infty}.$ This identifies the
Hilbert space dual to $\Dom_0(U)$ with $(\Id+U^*U)\Dom_0(U)\subset H_{-\infty}.$

The following observation
is useful:
\begin{lemma}\label{lem15} Let $U$ be a tame operator of non-negative order and $K$ a
  smoothing operator. The map $K:\Dom_0(U)\to H_0$ is compact.
\end{lemma}
\begin{proof} This follows because the unit ball in $\Dom_0(U)$ is contained in
  the unit ball in $H_0.$ Thus its image under $K$ is a precompact subset of
  $H_0.$
\end{proof}

\begin{proposition}\label{prp144} Let $T$ be a tame elliptic operator of non-positive
  order. Let $U$ denote a parametrix for $T,$  then
  $T:H_0\to \Dom_0(U)$ is a Fredholm operator. Moreover we have
\begin{equation}
\Ind(T)=\tind(T).
\end{equation}
\end{proposition}
\begin{proof} We first observe that $TH_0\subset \Dom_0(U).$ Let $x\in H_0$ and
  let $<x_n>\subset H_{\infty}$ converge to $x$ in $H_0.$ Note that
  $<Tx_n>\subset H_{\infty}.$ As $T$ is of non-positive order $Tx_n\to Tx$ in
  $H_0$ as well. As $UTx_n=x_n-K_2x_n$ and $UTx=x-K_2x,$ it is clear that
  $UTx_n$ converges to $UTx$ in $H_0;$ thus verifying that $TH_0\subset
  \Dom_0(U).$ Hence, the operator $U$ maps $\Dom_0(U)$ boundedly onto $H_0.$ It
  follows from Lemma~\ref{lem15} that $T:H_0\to \Dom_0(U)$ has a left and right
  inverse, up to a compact error, and is therefore a Fredholm operator. The
  null space of $T$ acting on $H_{-\infty}$ is contained in $H_{\infty};$ hence
  $\ker T$ does not depend on the topology. To complete the proof we need to
  show that the $\coker T$ is isomorphic to the kernel of the formal adjoint
  $T^*.$

As $T^*T(\Id+U^*U)=I+T^*T+K_1,$ and $(\Id+U^*U)T^*T=I+T^*T+K_2,$ for smoothing
operators $K_1, K_2$ it follows from the assumption that $\Id+T^*T$ is tamely
elliptic that $\Id+U^*U$ is tamely elliptic as well.  Hence the null-space of
$(\Id+U^*U)$ is contained in $H_{\infty}$ and is therefore trivial.  Indeed, as
an operator on $H_0$ we can identify $(\Id+U^*U)$ as the self adjoint operator
defined by Friedrichs' extension from the symmetric quadratic form
$\langle\cdot,\cdot\rangle_U.$ This operator is self adjoint and therefore
invertible.  The cokernel of $T$ is isomorphic to the set of $y$ in $\Dom_0(U)$
with
\begin{equation}
\langle Tx,y\rangle_U=0
\end{equation}
for all $x\in H_{\infty}.$ Using the extension of the pairing
$\langle\cdot,\cdot\rangle$ to $H_{\infty}\times H_{-\infty},$ this implies
that $T^*(\Id+U^*U)y=0.$ As $T^*(\Id+U^*U)$ is tamely elliptic we see that
$y\in H_{\infty}.$ Thus the $\coker T$ is isomorphic to $(\Id+U^*U)^{-1}\ker
T^*,$ completing the proof of the proposition.
\end{proof}

A bounded operator $P$ is a projection if $P^2=P.$ The
following elementary fact about bounded projections is very useful:
\begin{lemma}\label{lem16} Suppose that $H$ is a Hilbert space and $P:H\to H$ is a bounded
  projection operator, then $\range P$ is a closed subspace of $H.$
\end{lemma}
\begin{proof} This follows immediately from the observation that $\range P=\Ker(\Id-P).$
\end{proof}

\begin{definition} A tame projection $P$ is a projection that is also a tame
  operator of order $0.$ A pair of tame projections $(P,R)$ defines a
tame Fredholm pair if the operator
\begin{equation}
T=RP+(I-R)(I-P)
\end{equation}
is tamely elliptic. We let $U$ denote a parametrix for $T.$
\end{definition}

\begin{proposition}\label{prp14} If $(P,R)$ are a tame Fredholm pair, then
\begin{equation}
RP:H_0\cap\range P\longrightarrow \Dom_0(U)\cap\range R
\label{eqna10}
\end{equation}
is a Fredholm operator
\end{proposition}
\begin{proof} We first need to show that $RP$ is bounded from $H_0$ to
  $\Dom_0(U).$ Let $x\in H_{\infty},$ as $URP=UTP$ we have
\begin{equation}
\begin{split}
\|RP x\|_U^2&=\|RP x\|^2_0+\|UTP x\|^2_0\\
&=\|RP x\|^2_0+\|(\Id-K_2)P x\|^2_0\\
&\leq C\|x\|^2_0.
\end{split}
\end{equation}
This establishes the boundedness. Moreover $R:\Dom_0(U)\to\Dom_0(U)$ is a
bounded map. This follows from the identity $RT=TP,$ which implies that
\begin{equation}
UR(\Id-K_1)=(\Id-K_2)PU
\end{equation}
As $URK_1$ is a smoothing operator and $P$ is order $0,$ for $x\in\Dom_0(U)$ we have
\begin{equation}
\|URx\|=\|(PU+URK_1-K_2PU)x\|\leq C[\|Ux\|+\|x\|].
\end{equation}
Hence Lemma~\ref{lem16} implies that $H_0\cap\range P$ and
$\Dom_0(U)\cap\range R$ are closed subspaces of their respective Hilbert spaces
and are therefore Hilbert spaces in their own rights.  Similarly we can show
that $PUR:\Dom_0(U)\to H_0$ is bounded. That the map is Fredholm follows from
the identities:
\begin{equation}
\begin{split}
&(PUR)(RP)=PUTP=P(\Id-K_2)P\\
&(RP)(PUR)=RTUR=R(\Id-K_1)R,
\end{split}
\label{eqna15}
\end{equation}
which imply that the map in~\eqref{eqna10}  is a bounded map
between Hilbert spaces,  invertible up to a compact error.
\end{proof}

\begin{definition} For a $(P,R)$ a tame Fredholm pair we let $\Rind(P,R)$
  denote the Fredholm index of the operator in~\eqref{eqna10}. The number
  $\Rind(P,R)$ is also called the \emph{relative index} of $P$ and $R.$
\end{definition}

The relative index can be identified with  difference of the dimensions of
null spaces. This is a relative index analogue of Proposition~\ref{prp144}.
\begin{proposition}\label{reltr} Let $(P,R)$ be a tame Fredholm pair, then
\begin{equation}
\Rind(P,R)=\dim[\Ker RP\restrictedto_{PH_{-\infty}}]-\dim[\Ker P^*R^*\restrictedto_{R^* H_{-\infty}}]
\label{5.20.2}
\end{equation}
\end{proposition}
\begin{proof}
By definition 
\begin{equation}
\Rind(P,R)=\dim[\Ker RP\restrictedto_{PH_{0}}]-\dim[\coker PR\restrictedto_{P
    H_{0}}].
\label{5.20.3}
\end{equation}
Ellipticity of $T$ implies that
\begin{equation}
\Ker RP\restrictedto_{PH_{-\infty}}\subset P H_{\infty},
\end{equation}
so the  first terms in~\eqref{5.20.2} and~\eqref{5.20.3} agree. Using the
representation of $\Dom_{0}(U)'$ as $(\Id+U^*U)\Dom_{0}(U)$ we see that
\begin{equation}
\coker PR\restrictedto_{P H_{0}}\simeq \{v\in\Dom_0(U)'\cap\range R^*:\:
\langle RP u,v\rangle=0\text{ for all }u\in H_0\}.
\label{5.20.4}
\end{equation}
The right hand side of~\eqref{5.20.4} is certainly contained in $\Ker
P^*R^*\restrictedto_{R^* H_{-\infty}}.$  Again ellipticity implies that this
subspace is contained in $H_{\infty}$ and therefore is also contained in
$\Dom_0(U)'\cap\range R^*.$ This proves~\eqref{5.20.2}
\end{proof}

\begin{remark} If $K:H\to H$ is a trace class operator, then we use $\Tr_H(K)$
  to denote its trace \emph{as an operator on $H.$} A priori the definition of
  the trace class and the value of the trace appear to depend on the inner product
  on $H.$ In many cases one can show that the value of the trace is independent
  of the inner product.
\end{remark}

The operators $RK_1R$ and $PK_2P$ are smoothing operators and therefore trace
class operators on $H_0.$ In order to use the standard trace formula for the
index we need to show that $RK_1R:\Dom_0(U)\to\Dom_0(U)$ is trace class and
that
\begin{equation}
\Tr_{\Dom_0(U)}(RK_1R)=\Tr_{H_0}(RK_1R).
\label{eqna11}
\end{equation}

Theorem VI.2.23 in~\cite{Kato} implies that $\sqrt{\Id+U^*U}$ is a (possibly
unbounded) self adjoint operator with domain $\Dom_0(U).$ Moreover
$\sqrt{\Id+U^*U}$ has a bounded inverse, $(\Id+U^*U)^{-\frac 12}$ and
\begin{equation}
\langle x,y\rangle_U=\langle\sqrt{\Id+U^*U}x,\sqrt{\Id+U^*U}y\rangle,\quad
\forall x,y\in \Dom_0(U).
\label{eqna12}
\end{equation}
To show that $RK_1R$ is a trace class operator on $\Dom_0(U)$ we need to verify
that there is an $M$ so that, for every pair of orthonormal bases $\{f_j^1\},
\{f_j^2\}$ of $\Dom_0(U),$ we have
\begin{equation}
\sum_{j=1}^{\infty}|\langle RK_1R f_j^1,f_j^2\rangle_U|\leq M.
\end{equation}
As $RK_1R$ is a smoothing operator, $RK_1Rf_j^1\in
H_{\infty}\subset\Dom_0(\Id+U^*U),$ hence the Friedrichs extension theorem
and~\eqref{eqna12} imply that
\begin{equation}
\sum_{j=1}^{\infty}|\langle RK_1R f_j^1,f_j^2\rangle_U|=
\sum_{j=1}^{\infty}|\langle (\Id+U^*U)RK_1R f_j^1,f_j^2\rangle|.
\end{equation}
For $i=1$ or $2$ we let
\begin{equation}
e_j^i=(\Id+U^*U)^{\frac 12}f^i_j,
\end{equation}
where $\{e^1_j\},$ $\{e^2_j\}$ are orthonormal bases for $H_0.$ As
  $(\Id+U^*U)^{-\frac 12}$ is a bounded self adjoint operator, we obtain:
\begin{multline}
\sum_{j=1}^{\infty}|\langle RK_1R f_j^1,f_j^2\rangle_U|=\\
\sum_{j=1}^{\infty}|\langle(\Id+U^*U)^{-\frac
  12}(\Id+U^*U)RK_1R(\Id+U^*U)^{-\frac 12} e_j^1,e_j^2\rangle|.
\label{eqna13}
\end{multline}
The operator $(\Id+U^*U)RK_1R$ is a smoothing operator and therefore a trace
class operator on $H_0.$ The operator $(\Id+U^*U)^{-\frac 12}$ is bounded on
$H_0,$ which shows that $(\Id+U^*U)^{-\frac
  12}(\Id+U^*U)RK_1R(\Id+U^*U)^{-\frac 12}$ is also a trace class operator on
$H_0.$ Hence there exists an upper bound $M$ for the sum on the right hand side
of~\eqref{eqna13} valid for any  pair of orthonormal bases $\{e_j^1\},
\{e_j^2\}.$ This completes the proof  that 
$$RK_1R:\Dom_0(U)\to\Dom_0(U)$$ 
is  a trace class operator. 

To compute the trace we select an orthonormal basis
$\{f_j\}$ of $\Dom_0(U)$ and let $\{e_j\}$ be the corresponding orthonormal
basis of $H_0.$ Arguing as above we conclude that
\begin{equation}
\begin{split}
\Tr_{\Dom_0(U)}(RK_1R)&=\sum_{j=1}^{\infty}
\langle(\Id+U^*U)^{-\frac
  12}(\Id+U^*U)RK_1R(\Id+U^*U)^{-\frac 12} e_j,e_j\rangle\\
&=\Tr_{H_0}\left(\Id+U^*U)^{-\frac 12}(\Id+U^*U)RK_1R(\Id+U^*U)^{-\frac 12}\right).
\end{split}
\end{equation}
The operator $(\Id+U^*U)^{-\frac 12}(\Id+U^*U)RK_1R$ is $H_0$-trace class and
the operator $(\Id+U^*U)^{-\frac 12}$ is $H_0$ bounded. Therefore
\begin{multline}
\Tr_{H_0}\left((\Id+U^*U)^{-\frac 12}(\Id+U^*U)RK_1R(\Id+U^*U)^{-\frac 12}\right)=\\
\Tr_{H_0}[(\Id+U^*U)^{-\frac 12}]^2(\Id+U^*U)RK_1R=\Tr_{H_0}RK_1R.
\end{multline}
This completes the proof of the following result:
\begin{proposition}\label{prp17} If $(P, R)$ is a tame Fredholm pair then
  $RK_1R:\Dom_0(U)\to\Dom_0(U)$ is a trace class operator and
\begin{equation}
\Tr_{\Dom_0(U)}(RK_1R)=\Tr_{H_0}(RK_1R).
\end{equation}
\end{proposition}

To obtain a trace formula for $\Rind(P,R)$ we need to compute the traces of
$RK_1R$ and $PK_2P$ restricted to the ranges of $R$ and $P$ respectively. Let
$\{e_j^1\}$ be an orthonormal basis for the range of $P$ and $\{e_j^2\}$ be a
orthonormal basis for the orthocomplement of the range. Note that $P^*e_j^2=0$
for all $j$ and therefore
\begin{equation}
\langle PK_2P e_j^2,e_j^2\rangle=\langle PK_2P e_j^2,P^*e_j^2\rangle=0.
\end{equation}
This shows that
\begin{equation}
\Tr_{H_0}(PK_2P)=\sum_{j=1}^{\infty}\langle PK_2P e_j^1,e_j^1\rangle.
\end{equation}
Similar considerations apply to $\Tr_{\Dom_0(U)}(RK_1R)$ and
$\Tr_{H_0}(RK_1R).$ From these observations we obtain:
\begin{theorem}\label{thma1} If $(P,R)$ is a tame Fredholm pair, then
\begin{equation}
\Rind(P,R)=\Tr_{H_0}(PK_2P)-\Tr_{H_0}(RK_1R).
\label{eqna14}
\end{equation}
\end{theorem}
\begin{proof} This follows from~\eqref{eqna15}, Proposition~\ref{prp17} and
  Theorem 15 in Chapter 30 of~\cite{lax}.
\end{proof}
In applications, $H_0$ is often $L^2(M;E)$ where $M$ is a compact manifold and
$E\to M$ is a vector bundle. Theorem~\ref{thma1} and Lidskii's Theorem imply
that $\Rind(P,R)$ can be computed by integrating the Schwartz kernels of
$PK_2P$ and $RK_1R$ along the diagonal. This is a very useful fact.

Finally we have a logarithmic property for relative indices.
\begin{theorem}\label{thm11} Let $(P,Q)$ and $(Q,R)$ be tame Fredholm
  pairs. For an appropriately defined Hilbert space, $H_{UV},$ $RQP$ is a tame
  Fredholm operator from $\range P\cap H_0$ to $\range Q\cap H_{UV}$ and
\begin{equation}
\Ind(RQP)=\Rind(P,Q)+\Rind(Q,R).
\end{equation}
\end{theorem}
\begin{proof} For $S=QP+(\Id-Q)(\Id-P)$ and $T=RQ+(\Id-R)(\Id-Q),$ let $U,V$ be
  parametrices, with
\begin{equation}
\begin{split}
&SU=\Id-K_1\quad US=\Id-K_2\\
&TV=\Id-K_3\quad VT=\Id-K_4.
\end{split}
\label{5.20.6}
\end{equation}
We define the Hilbert spaces $H_U, H_V$ and $H_{UV}$ as the closures of $H_{\infty}$
with respect to the inner products
\begin{equation}
\|x\|^2_{U}=\|x\|^2+\|Ux\|^2,\quad
\|x\|^2_{V}=\|x\|^2+\|Vx\|^2,\quad
\|x\|^2_{UV}=\|x\|^2+\|UVx\|^2.
\end{equation}
By definition $QP:\range P\cap H_0\to \range Q\cap H_{U}$ and $RQ:\range Q\cap
H_0\to \range R\cap H_V$ are Fredholm with indices $\Rind(P,Q), \Rind(Q,R)$
respectively. To prove the theorem we need to show that 
\begin{equation}
RQ:\range Q\cap H_U\longrightarrow\range R\cap H_{UV} 
\label{5.20.7}
\end{equation}
is Fredholm and has index $\Rind(Q,R).$ The proofs of these statements make use
of the commutation relations $QS=SP, RT=TQ,$ which, along with~\eqref{5.20.6}
imply that
\begin{equation}
UQ(\Id-K_1)=(\Id-K_2)PU\text{ and }VR(\Id-K_3)=(\Id-K_4)QV.
\label{5.20.8}
\end{equation}

First we show that $Q$ acts boundedly on $H_U$ and $R$ acts boundedly on
$H_{UV},$ so that we can apply Lemma~\ref{lem16}. Using~\eqref{5.20.8} and the
fact that the smoothing operators are a two sided ideal we see that for $x\in
H_{\infty},$ we have
\begin{equation}
\|UQx\|=\|(UQK_1+PU-K_2 PU)x\|\leq C_1[\|Ux\|+\|x\|.
\end{equation}
In the following computation $K$ denotes a variety of smoothing operators:
\begin{multline}
\|UV Rx\|=\|(UQV+K)x\|\\=\|(PUV+K)x\|\leq C_2[\|UV
  x\|+\|x\|].
\end{multline}
The constants $C_1, C_2$ are independent of $x.$ This shows that $\range Q\cap
H_U$ and $\range R\cap H_{UV}$ are closed subspaces. 

The operator $QVR$ is a parametrix for the operator in~\eqref{5.20.7}. First we
see that $QVR:\range R\cap H_{UV}\to\range Q\cap H_U$ is bounded:
\begin{equation}
\| UQVR x\|=\|(UQK_1+PU-K_2PU)VR\|\leq C\|x\|_{UV}.
\end{equation}
That it is a parametrix follows from
\begin{equation}
(QVR)(RQ)=Q(\Id- K_4)Q\text{ and }(RQ)(QVR)=R(\Id-K_3)R.
\end{equation}
The null-space of $RQ\restrictedto_{Q H_U}$ agrees with the null-space of $RQ$
acting on $QH_{-\infty}.$ We use the operator $(\Id+V^*U^*U V)$ to identify the
dual space $H'_{UV}$ of $H_{UV}$ as a subspace of $H_{-\infty}$ via the pairing
$\langle\cdot,\cdot\rangle.$ With this identification, the cokernel of the
operator in~\eqref{5.20.7}, is isomorphic to the null-space of $Q^*R^*$ acting
on $R^* H'_{UV}.$ As in the proof of Proposition~\ref{reltr}, we see that the
cokernel is therefore isomorphic to the null-space of $Q^*R^*$ acting on
$R^*H_{-\infty}.$ This, along with Proposition~\ref{reltr},  shows that the
index of the operator in~\eqref{5.20.7} equals $\Rind(Q,R).$ The theorem now
follows from the standard logarithmic property for the Fredholm indices of
operators acting on Hilbert spaces.
 
\end{proof}


\begin{thebibliography}{10}

\bibitem{BogomolovDeOliveira}
{\sc F.~A.Bogomolov and B.~de~Oliveira}, {\em Stein small deformations of
  strictly pseudoconvex surfaces}, in Birational Algebraic Geometry (Baltimore,
  MD, 1996), vol.~207 of Contemp. Math., Providence, RI, 1997, Amer. Math.
  Soc., pp.~25--41.

\bibitem{ABKLR}
{\sc B.~Aebischer, M.~Borer, M.~K{\"a}lin, C.~Leuenberger, and H.~Reimann},
  {\em Symplectic Geometry}, vol.~124 of Progress in Mathematics,
  Birkh{\"a}user, Basel, Boston and Berlin, 1994.

\bibitem{BBW}
{\sc B.~Booss-Bavnbek and K.~P. Wojciechowsi}, {\em Elliptic Boundary Problems
  for the Dirac Operator}, Birkh\"auser, Boston, 1996.

\bibitem{BoutetdeMonvel4}
{\sc L.~{Boutet de Monvel}}, {\em On the index of {T}oeplitz operators of
  several complex variables}, Invent. Math., 50 (1979), pp.~249--272.

\bibitem{BoutetdeMonvel-Guillemin1}
{\sc L.~{Boutet de Monvel} and V.~Guillemin}, {\em The spectral theory of
  {T}oeplitz operators}, vol.~99 of Ann. of Math. Studies, Princeton University
  Press, 1981.

\bibitem{duistermaat}
{\sc J.~Duistermaat}, {\em The heat kernel Lefschetz fixed point theorem
  formula for the Spin-c Dirac Operator}, Birkh\"auser, Boston, 1996.

\bibitem{Epstein}
{\sc C.~L. Epstein}, {\em A relative index on the space of embeddable
  {C}{R}-structures, {I}, {I}{I}}, Annals of Math., 147 (1998), pp.~1--59,
  61--91.

\bibitem{Epstein5}
\leavevmode\vrule height 2pt depth -1.6pt width 23pt, {\em Geometric bounds on
  the relative index}, Jour. Inst. Math. Jussieu, 1 (2002), pp.~441--465.

\bibitem{Epstein4}
\leavevmode\vrule height 2pt depth -1.6pt width 23pt, {\em Subelliptic
  {S}pin${}_{\bbC}$ {D}irac operators,{I}}, to appear Annals of Math.,  (2006),
  pp.~1--36.

\bibitem{Epstein3}
\leavevmode\vrule height 2pt depth -1.6pt width 23pt, {\em Subelliptic
  {S}pin${}_{\bbC}$ {D}irac operators,{II}}, to appear Annals of Math.,
  (2006), pp.~1--54.

\bibitem{EpsteinMelrose2}
{\sc C.~L. Epstein and R.~Melrose}, {\em Shrinking tubes and the
  {$\bar{\partial}$}-{N}eumann problem}, preprint, 1990.

\bibitem{EpsteinMelrose}
\leavevmode\vrule height 2pt depth -1.6pt width 23pt, {\em Contact degree and
  the index of {F}ourier integral operators}, Math. Res. Letters, 5 (1998),
  pp.~363--381.

\bibitem{EpsteinMelrose3}
\leavevmode\vrule height 2pt depth -1.6pt width 23pt, {\em The {H}eisenberg
  algebra, index theory and homology}, preprint, 2004.

\bibitem{FollandKohn1}
{\sc G.~Folland and J.~Kohn}, {\em The {N}eumann problem for the
  {C}auchy-{R}iemann complex}, no.~75 in Ann. of Math. Studies, Princeton Univ.
  Press, 1972.

\bibitem{GGK}
{\sc V.~Guillemin, V.~Ginzburg, and Y.~Karshon}, {\em Moment maps, cobordisms,
  and {H}amiltonian group actions}, vol.~98 of Mathematical Surveys and
  Monographs, American Mathematical Society, Providence, RI, 2002.

\bibitem{GuilleminStenzel1}
{\sc V.~Guillemin and M.~Stenzel}, {\em Grauert tubes and the homogeneous
  {M}onge-{A}mpere equation {I}}, J. Differential Geom., 34 (1991),
  pp.~561--570.

\bibitem{GuilleminStenzel2}
\leavevmode\vrule height 2pt depth -1.6pt width 23pt, {\em Grauert tubes and
  the homogeneous {M}onge-{A}mpere equation {II}}, J. Differential Geom., 35
  (1992), pp.~627--641.

\bibitem{Hormander3}
{\sc L.~H{\"o}rmander}, {\em The Analysis of Linear Partial Differential
  Operators}, vol.~3, Springer-Verlag, Berlin{,} Heidelberg{,} New York{,}
  Tokyo, 1985.

\bibitem{Kang}
{\sc H.~Kang}, {\em Polynomial Hulls and Relative Indices in Dimension One},
  Thesis, University of PA, Philadelphia, PA, 2003.

\bibitem{Kato}
{\sc T.~Kato}, {\em Perturbation Theory for Linear Operators, corrected 2nd
  printing}, vol.~132 of Grundlehren der mathematischen Wissenschaften,
  Springer Verlag, Berlin Heidelberg New York, 1980.

\bibitem{LawsonMichelsohn}
{\sc H.~B. {Lawson Jr.} and M.-L. Michelsohn}, {\em Spin Geometry}, vol.~38 of
  Princeton Mathematical Series, Princeton University Press, 1989.

\bibitem{lax}
{\sc P.~D. Lax}, {\em Functional Analysis}, Pure and Applied Mathematics, John
  Wiley and Sons, New York, 2002.

\bibitem{LeichtNestTsy}
{\sc E.~Leichtnam, R.~Nest, and B.~Tsygan}, {\em Local formula for the index of
  a {F}ourier integral operator}, Jour. of Diff. Geo., 59 (2001), pp.~269--300.

\bibitem{LS}
{\sc L.~Lempert and R.~Szoke}, {\em Global solutions of the homogeneous complex
  {M}onge-{A}mpere equation and complex structures on the tangent bundle of
  {R}iemannian manifolds}, Math. Ann., 2904 (1991), pp.~689--712.

\bibitem{Lisca}
{\sc P.~Lisca}, {\em On lens spaces and their symplectic fillings},
  arXiv:math.SG/0203006,  (2002).

\bibitem{Morgan}
{\sc J.~W. Morgan}, {\em The {S}eiberg-{W}itten Equations and Applicatgions to
  the Topology of Smooth Four-Manifolds}, Princeton University Press,
  Princeton, NJ, 1996.

\bibitem{NestTsygan}
{\sc R.~Nest and B.~Tsygan}, {\em Deformations of symplectic {L}ie algebroids,
  {D}eformations of holomorphic symplectic structures, and index theorems},
  Asian Jour. of Math., 5 (2001), pp.~599--633.

\bibitem{OhtaOno1}
{\sc H.~Ohta and K.~Ono}, {\em Simple singularities and topology of
  symplectically filling 4-manifolds}, Comm. Math. Helv., 74 (1999),
  pp.~575--590.

\bibitem{OhtaOno2}
\leavevmode\vrule height 2pt depth -1.6pt width 23pt, {\em Symplectic fillings
  of the link of simple elliptic singularities}, Jour. Reine Angew. Math., 565
  (2003), pp.~183--205.

\bibitem{OzbagciStipsicz}
{\sc B.~Ozbagci and A.~I. Stipsicz}, {\em Contact 3-manifolds with infinitely
  many {S}tein fillings}, Proceedings of the American Mathematical Society, 132
  (2003), pp.~1549--1558.

\bibitem{OzbagciStipsicz2}
\leavevmode\vrule height 2pt depth -1.6pt width 23pt, {\em Surgery on Contact
  3-Manifolds and Stein Surfaces}, Springer Verlag and J\'anos Bolyai
  Mathematical Society, Budapest, 2004.

\bibitem{Seeley0}
{\sc R.~Seeley}, {\em Singular integrals and boundary value problems}, Amer.
  Jour. of Math., 88 (1966), pp.~781--809.

\bibitem{SegalWilson}
{\sc G.~Segal and G.~Wilson}, {\em Loop groups and equations of {K}d{V} type},
  Publ. Math., 61 (1985), pp.~5--65.

\bibitem{Stipsicz2}
{\sc A.~I. Stipsicz}, {\em Gauge theory and {S}tein fillings of certain
  3-manifolds}, in Proceedings of the G{\"o}kova Geometry-Topology Conference,
  2001, pp.~115--131.

\bibitem{stipsicz}
\leavevmode\vrule height 2pt depth -1.6pt width 23pt, {\em On the geography of
  {S}tein fillings of certain 3-manifolds}, Michigan Math. J., 51 (2003),
  pp.~327--337.

\bibitem{tataru}
{\sc G.~R. Tataru}, {\em Adiabatic limit and Szeg\H o Projections}, Thesis,
  MIT, Cambridge, MA, 2003.

\bibitem{Weinstein}
{\sc A.~Weinstein}, {\em {Some questions about the index of quantized contact
  transformations}}, RIMS Kokyuroku, 1014 (1997), pp.~1--14.

\end{thebibliography}
\end{document}